\documentclass{amsart}
\usepackage{amsfonts,amssymb,amsmath,amsthm}
\usepackage{url}
\usepackage{enumerate}

\newtheorem{theorem}{Theorem}[section]

\theoremstyle{definition}
\newtheorem{definition}[theorem]{Definition}

\numberwithin{equation}{section}

\author{U. Bekbaev}

\address{
	Department of Science in Engineering, Faculty of Engineering,\\
	IIUM, 50728, KL, Malaysia\\
}
\email{bekbaev@iium.edu.my}
\begin{document}
	\begin{center}\small{In the name of Allah, Most Gracious, Most Merciful.}\end{center}
	\vspace{1cm}
	\title{Complete classification of a class of $3$-dimensional algebras }

	\begin{abstract} 		
		A complete classification of a class of $3$-dimensional algebras is provided. In algebraically closed field $\mathbb{F}$ case this class is an open, dense (in Zariski topology) subset of $\mathbb{F}^{27}$. \end{abstract}
		\maketitle
	
	%%%%% END OF TITLE PAGE %%%%%%%%%%%%%%%%%%%%%%%%%%%%%%%%%%%%%%%%%%%%%%
	
	\section{Introduction}
The classification problem of different mathematical structures, in particular finite dimensional algebras, is one of the central problems in algebra. One can hope to get a complete classification of small dimensional algebras though it may be too difficult for higher dimensional cases. There are many investigations on special classes of small dimensional algebras, for example, \cite{A, A1, Ay, J, Bu, C, D, M, M1, MA, P2000, R}. In \cite{B} a complete classification of all $2$-dimensional algebras, over algebraically closed fields, is presented in terms of their matrix of structure constants and in \cite{B3} it is done over the field of real numbers. These classifications can be used for classification of any class of $2$-dimensional algebras defined by a system of polynomial identities such as commutative, associative, Jordan, alternative and etcetera \cite{B1, B3}. The results of \cite{B} are used in \cite{B2} to describe their automorphism groups and derivation algebras. The present paper deals with complete classification of a huge part of $3$-dimensional algebras. By huge we mean density of that set, in Zariski topology, in the space of all $3$-dimensional algebras over algebraically closed fields. For the classification problem of special classes of $3$-dimensional algebras one can see \cite{Az, E, F}. 

\section{Preliminaries}

Let $\mathbb{F}$ be any field and $A\otimes B$ stand for the Kronecker product of matrices $A$,$B$, that is for the matrix with blocks $(a_{ij}B),$ where $A=(a_{ij})$, $B$ are
matrices over $\mathbb{F}$.

Let ($\mathbb{A}, \cdot)$ be a $m$-dimensional algebra over $\mathbb{F}$ and $e=(e_1,e_2,...,e_m)$ be its basis. Then the
bilinear map $\cdot$ is represented by a matrix $A=(A^k_{ij})\in M(m\times m^2;\mathbb{F})$ as follows \[\mathbf{u}\cdot
\mathbf{v}=eA(u\otimes v),\] for $\mathbf{u}=eu,\mathbf{v}=ev,$
where $u = (u^1, u^2, ..., u^m),$ and $v = (v^1, v^2, ..., v^m)$ are column coordinate vectors of $\mathbf{u}$ and
$\mathbf{v},$ respectively.
The matrix $A\in M(m\times m^2;\mathbb{F})$ defined above is called the matrix of structural constants (MSC) of
$\mathbb{A}$ with respect to the basis $e$. 

If $e'=(e'_1,e'_2,...,e'_m)$ is another basis of $\mathbb{A}$, $e'g=e$ with $g\in G=GL(m;\mathbb{F})$, and  $A'$ is MSC
of $\mathbb{A}$ with respect to $e'$ then it is known that
\[A'=gA(g^{-1})^{\otimes 2}\] is valid. Thus, we can reformulate the isomorphism of algebras as follows.

\begin{definition} Two $m$-dimensional algebras $\mathbb{A}$, $\mathbb{B}$ over $\mathbb{F}$, given by
	their matrices of structure constants $A$, $B$, are said to be isomorphic if $B=gA(g^{-1})^{\otimes 2}$ holds true
	for some $g\in GL(m;\mathbb{F})$.\end{definition}

Further we consider only $m=3$ case, assume that a basis $e$ is fixed and we do not make a difference between the algebra
$\mathbb{A}$ and its MSC 
\[A=\left(
\begin{array}{ccccccccc}
A^{1}_{1,1} & A^{1}_{1,2}& A^{1}_{1,3}& A^{1}_{2,1} & A^{1}_{2,2}& A^{1}_{2,3}& A^{1}_{3,1}& A^{1}_{3,2}& A^{1}_{3,3}\\
A^{2}_{1,1} & A^{2}_{1,2}& A^{2}_{1,3}& A^{2}_{2,1} & A^{2}_{2,2}& A^{2}_{2,3}& A^{2}_{3,1}& A^{2}_{3,2}& A^{2}_{3,3}\\
A^{3}_{1,1} & A^{3}_{1,2}& A^{3}_{1,3}& A^{3}_{2,1} & A^{3}_{2,2}& A^{3}_{2,3}& A^{3}_{3,1}& A^{3}_{3,2}& A^{3}_{3,3}\\
\end{array}\right).\] We use notations $\mathcal{A}_1=\left(
\begin{array}{ccc}
A^{1}_{1,1} & A^{1}_{1,2}& A^{1}_{1,3}\\
A^{2}_{1,1} & A^{2}_{1,2}& A^{2}_{1,3}\\
A^{3}_{1,1} & A^{3}_{1,2}& A^{3}_{1,3}
\end{array}\right),$ \[ 
\mathcal{A}_2=\left(\begin{array}{ccc}
A^{1}_{2,1} & A^{1}_{2,2}& A^{1}_{2,3}\\
A^{2}_{2,1} & A^{2}_{2,2}& A^{2}_{2,3}\\
A^{3}_{2,1} & A^{3}_{2,2}& A^{3}_{2,3}	\end{array}\right),\
\ \mathcal{A}_3=\left(\begin{array}{ccc}
A^{1}_{3,1}& A^{1}_{3,2}& A^{1}_{3,3}\\
A^{2}_{3,1}& A^{2}_{3,2}& A^{2}_{3,3}\\
A^{3}_{3,1}& A^{3}_{3,2}& A^{3}_{3,3}	\end{array}\right),
\]
\[\mathbf{Tr_1}(A)=(A^{1}_{1,1} + A^{2}_{2,1}+ A^{3}_{3,1}, A^{1}_{1,2} + A^{2}_{2,2}+ A^{3}_{3,2}, A^{1}_{1,3} + A^{2}_{2,3}+ A^{3}_{3,3}),\]
\[\mathbf{Tr_2}(A)=(A^{1}_{1,1} + A^{2}_{1,2}+ A^{3}_{1,3}, A^{1}_{2,1} + A^{2}_{2,2}+ A^{3}_{2,3}, A^{1}_{3,1} + A^{2}_{3,2}+ A^{3}_{3,3})\] and note that the equalities 
\begin{equation}\label{E1}\mathbf{Tr_1}(A')=\mathbf{Tr_1}(A)g^{-1},\ \ \mathbf{Tr_2}(A')=\mathbf{Tr_2}(A)g^{-1}\end{equation} hold true, whenever $A\in M(3\times 3^2;\mathbb{F})$, $g\in 
GL(3;\mathbb{F})$.

 Further we deal with only $3$-dimensional algebras $\mathbb{A}$ for which the system $\{\mathbf{Tr_1}(A),\mathbf{Tr_2}(A)\}$ is linear independent and
 for simplicity we use notations \[A=\left(
 \begin{array}{ccccccccc}
 \alpha_1 & \alpha_2& \alpha_3&\alpha_4&\alpha_5&\alpha_6&\alpha_7&\alpha_8&\alpha_9\\
 \beta_1 & \beta_2& \beta_3&\beta_4&\beta_5&\beta_6&\beta_7&\beta_8&\beta_9\\
 \gamma_1 & \gamma_2& \gamma_3&\gamma_4&\gamma_5&\gamma_6&\gamma_7&\gamma_8&\gamma_9\end{array}\right),\]
 $\mathbf{Tr_1}(A)=(\alpha_1 + \beta_4+ \gamma_7, \alpha_2 + \beta_5+ \gamma_8, \alpha_3 +\beta_6+ \gamma_9),$\\
 $\mathbf{Tr_2}(A)=(\alpha_1 +\beta_2+ \gamma_3, \alpha_4 + \beta_5+ \gamma_6, \alpha_7 + \beta_8+ \gamma_9)$.\\
 Due to (\ref{E1}) and the linear independence of $\{\mathbf{Tr_1}(A),\mathbf{Tr_2}(A)\}$ it is enough to consider only
 $\mathbf{Tr_1}(A)=(1,0,0),\mathbf{Tr_2}=(0,1,0)$ case so \[A=\left(
 \begin{array}{ccccccccc}
 \alpha_1 & \alpha_2& \alpha_3&\alpha_4&\alpha_5&\alpha_6&\alpha_7&\alpha_8&\alpha_9\\
 \beta_1 & \beta_2& \beta_3&\beta_4&\beta_5&\alpha_7+\beta_8-\alpha_3&\beta_7&\beta_8&\beta_9\\
 \gamma_1 & \gamma_2& -\alpha_1-\beta_2&\gamma_4&\gamma_5&1-\alpha_4-\beta_5&1-\alpha_1-\beta_4&-\alpha_2-\beta_5&-\alpha_7-\beta_8\end{array}\right).\]
 
 It is clear that $(1,0,0)g^{-1}=(1,0,0), (0,1,0)g^{-1}=(0,1,0)$ if and only if $g^{-1}=\left(
 \begin{array}{ccc}
 1 & 0& 0\\ 0 & 1& 0\\ a & b& c\end{array}\right)$ and therefore further we consider $A'=gA(g^{-1})^{\otimes 2}$ only with respect to such $g\in 
 GL(3;\mathbb{F})$. 
 
 One has $A'=gA(g^{-1})^{\otimes 2}=(g(\mathcal{A}_1+a\mathcal{A}_3)g^{-1},g(\mathcal{A}_2+b\mathcal{A}_3)g^{-1},cg\mathcal{A}_3g^{-1})$ 
 that is \[\mathcal{A}'_1=g(\mathcal{A}_1+a\mathcal{A}_3)g^{-1}, \mathcal{A}'_2=g(\mathcal{A}_2+b\mathcal{A}_3)g^{-1}, \mathcal{A}'_3=cg\mathcal{A}_3g^{-1},\ \mbox{ where}\] 
 \[\mathcal{A}_1+a\mathcal{A}_3=\left(
 \begin{array}{ccc}
 \alpha_1+a\alpha_7 & \alpha_2+a\alpha_8& \alpha_3+a\alpha_9\\
 \beta_1+a\beta_7 & \beta_2+a\beta_8& \beta_3+a\beta_9\\
 \gamma_1+a(1-\alpha_1-\beta_4) & \gamma_2-a(\alpha_2+\beta_5)& -\alpha_1-\beta_2-a(\alpha_7+\beta_8)\end{array}\right),\] 
 \[\mathcal{A}_2+b\mathcal{A}_3=\left(
 \begin{array}{ccc}
 \alpha_4+b\alpha_7&\alpha_5+b\alpha_8&\alpha_6+b\alpha_9\\
 \beta_4+b\beta_7&\beta_5+b\beta_8&\alpha_7+\beta_8-\alpha_3+b\beta_9\\
 \gamma_4+b(1-\alpha_1-\beta_4)&\gamma_5-b(\alpha_2+\beta_5)&1-\alpha_4-\beta_5-b(\alpha_7+\beta_8)\end{array}\right).\]
 
 Note that  columns of \[g\left(
 \begin{array}{ccc}
 x_1 & x_2& x_3\\ y_1 & y_2& y_3\\ z_1 & z_2& z_3\end{array}\right)g^{-1},\ \mbox{where}\ g^{-1}= \left(
 \begin{array}{ccc}
 1 & 0& 0\\ 0 & 1& 0\\ a & b& c\end{array}\right),\ \mbox{are} \] 
 \[1. 
 \begin{array}{c}
 x_1+ax_3 \\ y_1+ay_3 \\  -c^{-1}(a(x_1+ax_3)+b( y_1+ay_3)-(z_1+az_3)) \end{array}\]
 \[2. \begin{array}{c}
 x_2+bx_3\\
 y_2+by_3\\
 -c^{-1}(a(x_2+bx_3)+b( y_2+by_3)-(z_2+bz_3))\end{array}\]
 \[3.
 \begin{array}{c} cx_3\\ cy_3\\ -ax_3-by_3+z_3.\end{array}\]
  
 Therefore the columns of $\mathcal{A}'_1=g(\mathcal{A}_1+a\mathcal{A}_3)g^{-1}$ ($\mathcal{A}'_2=g(\mathcal{A}_2+b\mathcal{A}_3)g^{-1}$,\\ $\mathcal{A}'_3=cg\mathcal{A}_3g^{-1}$) are
 \begin{equation}\label{E2} \begin{array}{c}
 1. \begin{array}{c}
 \alpha_1+a\alpha_7+a(\alpha_3+a\alpha_9)\\
 \beta_1+a\beta_7+a(\beta_3+a\beta_9)\\
 c^{-1}[-a(\alpha_1+a\alpha_7+a(\alpha_3+a\alpha_9))-b(\beta_1+a\beta_7+a(\beta_3+a\beta_9))+\\ \gamma_1+ a(1-\alpha_1-\beta_4)+a(-\alpha_1-\beta_2-a(\alpha_7+\beta_8))]\end{array}\\ 
 \hspace*{\fill}\\
 2.  \begin{array}{c}
 \alpha_2+a\alpha_8+a(\alpha_3+a\alpha_9)\\
 \beta_2+a\beta_8+a(\beta_3+a\beta_9)\\
 c^{-1}[-a(\alpha_2+a\alpha_8+a(\alpha_3+a\alpha_9))-b(\beta_2+a\beta_8+a(\beta_3+a\beta_9))+\\ \gamma_2- a(\alpha_2+\beta_5)+a(-\alpha_1-\beta_2-a(\alpha_7+\beta_8))]\end{array}\\
  \hspace*{\fill}\\
 3. 
 \begin{array}{c}
 c(\alpha_3+a\alpha_9)\\
 c(\beta_3+a\beta_9)\\
 -a(\alpha_3+a\alpha_9)-b(\beta_3+a\beta_9)-\alpha_1-\beta_2-a(\alpha_7+\beta_8)\end{array}\end{array}\end{equation}
 ( respectively, \begin{equation}\label{E3} \begin{array}{c}
 1. 
 \begin{array}{c}
 \alpha_4+b\alpha_7+a(\alpha_6+b\alpha_9)\\
 \beta_4+b\beta_7+a(\alpha_7+\beta_8-\alpha_3+b\beta_9)\\
 c^{-1}[-a(\alpha_4+b\alpha_7+a(\alpha_6+b\alpha_9))-b(\beta_4+b\beta_7+a(\alpha_7+\beta_8-\alpha_3+\\ b\beta_9))+\gamma_4+ b(1-\alpha_1-\beta_4)+a(1-\alpha_4-\beta_5-b(\alpha_7+\beta_8))]\end{array}\\
  \hspace*{\fill}\\
 2. 
 \begin{array}{c}
 \alpha_5+b\alpha_8+b(\alpha_6+b\alpha_9)\\
 \beta_5+b\beta_8+b(\alpha_7+\beta_8-\alpha_3+b\beta_9)\\
 c^{-1}[-a(\alpha_5+b\alpha_8+b(\alpha_6+b\alpha_9))-b(\beta_5+b\beta_8+b(\alpha_7+\beta_8-\alpha_3+\\b\beta_9))+\gamma_5- b(\alpha_2+\beta_5)+b(1-\alpha_4-\beta_5-b(\alpha_7+\beta_8))]\end{array}\\
  \hspace*{\fill}\\
 3. 
 \begin{array}{c}
 c(\alpha_6+b\alpha_9)\\
 c(\alpha_7+\beta_8-\alpha_3+b\beta_9)\\
 -a(\alpha_6+b\alpha_9)-b(\alpha_7+\beta_8-\alpha_3+b\beta_9)+1-\alpha_4-\beta_5-b(\alpha_7+\beta_8)\end{array}\end{array}\end{equation}
 \begin{equation}\label{E4} \begin{array}{c}
 1.
 \begin{array}{c}
 c(\alpha_7+a\alpha_9)\\
 c(\beta_7+a\beta_9)\\
 -a(\alpha_7+a\alpha_9)-b(\beta_7+a\beta_9)+1-\alpha_1-\beta_4-a(\alpha_7+\beta_8)\end{array}\\
  \hspace*{\fill}\\
 2.
 \begin{array}{c}
 c(\alpha_8+b\alpha_9)\\
 c(\beta_8+b\beta_9)\\
 -a(\alpha_8+b\alpha_9)-b(\beta_8+b\beta_9)-\alpha_2-\beta_5-b(\alpha_7+\beta_8)\end{array}\\
  \hspace*{\fill}\\
 3.
 \begin{array}{c}
 c^2\alpha_9\\\
 c^2\beta_9\\
 c(-a\alpha_9-b\beta_9-\alpha_7-\beta_8).\end{array}\end{array}\end{equation}
 
To simplify the formulation of our main result we agree that all occurring in the result variables $\alpha_i, \beta_j, \gamma_k$  may have any values from $\mathbb{F}$ provided that the corresponding conditions are satisfied.
Further we have to consider matrices constructed by the use of rows listed as\\ \[\begin{array}{l} 1.\  (3\alpha_1+\beta_2+\beta_4-1,-\beta_1), \\
2.\ (2\alpha_2+\alpha_1+\beta_5+\beta_2,\\
3.\ (2\alpha_4+\beta_5-1, 1-\alpha_1-2\beta_4),\\
4.\ (\alpha_5, -\alpha_2-\alpha_4+1-3\beta_5),\end{array}\]
\\ for example, we use
$M_{(1,2,4)}(A)=M_{(1,2,4)}$ for the matrix consisting of the rows $1.,2.,4.$, and notation $\lambda\times i+ \mu \times j$, where $\lambda, \mu \in \mathbb{F}$, stands for the linear combination of rows $i.,j.$, $rk(A)$ stands for the rank of the matrix $A$. For the canonical forms of MSC we use notation $A_i$ without showing variables, to keep space,  on which it depends as far as it is clear from the context. Of course due to the notation
$A'=gA(g^{-1})^{\otimes 2}$ one expects, for obtained canonical MSCs, expressions in terms of $\alpha'_i, \beta'_j, \gamma'_k$ as far as $\alpha_i, \beta_j, \gamma_k$ are used for the entries of $A$. In the paper, for simplicity reasons, the obtained canonical MSCs also are given in terms of $\alpha_i, \beta_j, \gamma_k$. In some expressions we use * for an expression if there is no need for it or it is to big to present it there. 
 
Due to the big differences in formulations of results in $Ch.(\mathbb{F})\neq 2$ and $Ch.(\mathbb{F})= 2$ cases we have to formulate them separately.
 
\section{\bf Classification of $3$-dimensional algebras in $Ch.(\mathbb{F})\neq 2$ case}

\begin{theorem}. Let $Ch.(\mathbb{F})\neq 2$ and $Ch.(\mathbb{F})\neq 5$ ($Ch.(\mathbb{F})=5$) be a field over which any quadratic equation is solvable. In this case any $3$-dimensional algebra $\mathbb{A}$ over $\mathbb{F}$, for which $\{\mathbf{Tr_1}(A),\mathbf{Tr_2}(A)\}$ is linear independent, is isomorphic to only one algebra listed below by their (canonical) MSCs (respectively, is isomorphic to only one algebra listed below by their (canonical) MSCs except for the last two ones): 
\[A_1=\left(
\begin{array}{ccccccccc}
\alpha_1 & \alpha_2& \alpha_3&\alpha_4&\alpha_5&\alpha_6&0&0&1\\
\beta_1 & \beta_2& \beta_3&\beta_4&\beta_5&\beta_8-\alpha_3&\beta_7&\beta_8&\beta_9\\
\gamma_1 & \gamma_2& -\alpha_1-\beta_2&\gamma_4&\gamma_5&1-\alpha_4-\beta_5&1-\alpha_1-\beta_4&-\alpha_2-\beta_5&-\beta_8\end{array}\right)\cong \]
\[\left(
\begin{array}{ccccccccc}
\alpha_1 & \alpha_2& -\alpha_3&\alpha_4&\alpha_5&-\alpha_6&0&0&1\\
\beta_1 & \beta_2& -\beta_3&\beta_4&\beta_5&-\beta_8+\alpha_3&-\beta_7&-\beta_8&\beta_9\\
-\gamma_1 & -\gamma_2& -\alpha_1-\beta_2&-\gamma_4&-\gamma_5&1-\alpha_4-\beta_5&1-\alpha_1-\beta_4&-\alpha_2-\beta_5&\beta_8\end{array}\right),\] 
\[A_2=\left(
\begin{array}{ccccccccc}
\alpha_1 & \alpha_2& \alpha_3&\alpha_4&\alpha_5&\alpha_6&\alpha_7&\alpha_8&0\\
\beta_1 & \beta_2&\beta_3&\beta_4&\beta_5&\alpha_7-\alpha_3&0&0&1\\
\gamma_1 & \gamma_2& -\alpha_1-\beta_2&\gamma_4&\gamma_5&1-\alpha_4-\beta_5&1-\alpha_1-\beta_4&-\alpha_2-\beta_5&-\alpha_7\end{array}\right)\cong\]
\[\left(
\begin{array}{ccccccccc}
\alpha_1 & \alpha_2& -\alpha_3&\alpha_4&\alpha_5&-\alpha_6&-\alpha_7&-\alpha_8&0\\
\beta_1 & \beta_2&-\beta_3&\beta_4&\beta_5&-\alpha_7+\alpha_3&0&0&1\\
-\gamma_1 & -\gamma_2& -\alpha_1-\beta_2&-\gamma_4&-\gamma_5&1-\alpha_4-\beta_5&1-\alpha_1-\beta_4&-\alpha_2-\beta_5&\alpha_7\end{array}\right),\]
\[A_3=\left(
\begin{array}{ccccccccc}
\alpha_1 & \alpha_2& \alpha_3&0&0&1&\alpha_7&\alpha_8&0\\
\beta_1 & \beta_2& \beta_3&\beta_4&\beta_5&\alpha_7+\beta_8-\alpha_3&\beta_7&\beta_8&0\\
\gamma_1 & \gamma_2& -\alpha_1-\beta_2&\gamma_4&\gamma_5&1-\beta_5&1-\alpha_1-\beta_4&-\alpha_2-\beta_5&-\alpha_7-\beta_8\end{array}\right),\] where $\alpha_8\neq -1$,
\[A_4=\left(
\begin{array}{ccccccccc}
\alpha_1 & \alpha_2& \alpha_3&0&\alpha_5&1&\alpha_7&-1&0\\
\beta_1 & \beta_2& \beta_3&\beta_4&0&\alpha_7+\beta_8-\alpha_3&\beta_7&\beta_8&0\\
\gamma_1 & \gamma_2& -\alpha_1-\beta_2&\gamma_4&\gamma_5&1&1-\alpha_1-\beta_4&-\alpha_2&-\alpha_7-\beta_8\end{array}\right),\] where $\alpha_7+2\beta_8-\alpha_3\neq 0$,
\[A_5=\left(
\begin{array}{ccccccccc}
\alpha_1 & \alpha_2& \alpha_3&0&\alpha_5&1&\alpha_3-2\beta_8&-1&0\\
\beta_1 & \beta_2& \beta_3&0&\beta_5&-\beta_8&\beta_7&\beta_8&0\\
\gamma_1 & \gamma_2& -\alpha_1-\beta_2&\gamma_4&\gamma_5&1-\beta_5&1-\alpha_1&-\alpha_2-\beta_5&-\alpha_3+\beta_8\end{array}\right),\]
where $\beta_7\neq (2\beta_8-\alpha_3)\beta_8$,
\[A_6=\left(
\begin{array}{ccccccccc}
0 & \alpha_2& \alpha_3&0&\alpha_5&1&\alpha_3-2\beta_8&-1&0\\
\beta_1 & \beta_2& \beta_3&\beta_4&\beta_5&-\beta_8&(2\beta_8-\alpha_3)\beta_8&\beta_8&0\\
\gamma_1 & \gamma_2& -\beta_2&\gamma_4&\gamma_5&1-\beta_5&1-\beta_4&-\alpha_2-\beta_5&-\alpha_3+\beta_8\end{array}\right),\] where 
$(\alpha_3-2\beta_8)(\alpha_3-\beta_8)\neq 0$,
\[A_7=\left(
\begin{array}{ccccccccc}
\alpha_1 & \alpha_2& \alpha_3&0&\alpha_5&1&-\alpha_3&-1&0\\
0 & \beta_2& \beta_3&\beta_4&\beta_5&-\alpha_3&\alpha^2_3&\alpha_3&0\\
\gamma_1 & \gamma_2& -\beta_2&\gamma_4&\gamma_5&1-\beta_5&1-\alpha_1-\beta_4&-\alpha_2-\beta_5&0\end{array}\right),\] where 
$\alpha_3(\alpha^2_3+\beta_3)\neq 0$,
\[A_8=\left(
\begin{array}{ccccccccc}
\alpha_1 & 0& \alpha_3&0&\alpha_5&1&-\alpha_3&-1&0\\
\beta_1 & \beta_2& -\alpha^2_3&\beta_4&\beta_5&-\alpha_3&\alpha^2_3&\alpha_3&0\\
\gamma_1 & \gamma_2& -\beta_2&\gamma_4&\gamma_5&1-\beta_5&1-\alpha_1-\beta_4&-\beta_5&0\end{array}\right),\] where 
$\alpha_3\neq 0,1$,
\[A_{9}=\left(
\begin{array}{ccccccccc}
\alpha_1 & \alpha_2& 2\beta_8&0&\alpha_5&1&0&-1&0\\
\beta_1 & -\alpha_1& \beta_3&\beta_4&\beta_5&-\beta_8&0&\beta_8&0\\
\gamma_1 & \gamma_2& 0&\gamma_4&\gamma_5&1-\beta_5&1-\alpha_1-\beta_4&-\alpha_2-\beta_5&-\beta_8\end{array}\right),\] where 
$\beta_3\neq 0$,
\[A_{10}=\left(
\begin{array}{ccccccccc}
\alpha_1 & \alpha_2& 2\beta_8&0&\alpha_5&1&0&-1&0\\
\beta_1 & \beta_2& 0&\beta_4&\beta_5&-\beta_8&0&\beta_8&0\\
\gamma_1 & 0& -\alpha_1-\beta_2&\gamma_4&\gamma_5&1-\beta_5&1-\alpha_1-\beta_4&-\alpha_2-\beta_5&-\beta_8\end{array}\right),\] where 
$\beta_2\neq 0$,
\[A_{11}=\left(
\begin{array}{ccccccccc}
\alpha_1 & \alpha_2& 2\beta_8&0&\alpha_5&1&0&-1&0\\
\beta_1 & 0& 0&\beta_4&\beta_5&-\beta_8&0&\beta_8&0\\
0 & \gamma_2& -\alpha_1&\gamma_4&\gamma_5&1-\beta_5&1-\alpha_1-\beta_4&-\alpha_2-\beta_5&-\beta_8\end{array}\right),\] where 
$\beta_1\neq 0$,
\[A_{12}=\left(
\begin{array}{ccccccccc}
\alpha_1 & \alpha_2& 2\beta_8&0&\alpha_5&1&0&-1&0\\
0 & 0& 0&\beta_4&-\alpha_2&-\beta_8&0&\beta_8&0\\
\gamma_1 & \gamma_2& -\alpha_1&\gamma_4&\gamma_5&1+\alpha_2&1-\alpha_1-\beta_4&0&-\beta_8\end{array}\right),\] where 
$\beta_8\neq 0$,
\[A_{13}=\left(
\begin{array}{ccccccccc}
\alpha_1 & \alpha_2& 0&0&\alpha_5&1&0&-1&0\\
0 & 0& 0&\beta_4&\beta_5&0&0&0&0\\
\gamma_1 & \gamma_2& -\alpha_1&\gamma_4&0&1-\beta_5&1-\alpha_1-\beta_4&-\alpha_2-\beta_5&0\end{array}\right),\] where 
$1-\alpha_2-3\beta_5\neq 0$,
\[A_{14}=\left(
\begin{array}{ccccccccc}
\alpha_1 & 1-3\beta_5& 0&0&\alpha_5&1&0&-1&0\\
0 & 0& 0&\beta_4&\beta_5&0&0&0&0\\
\gamma_1 & \gamma_2& -\alpha_1&0&\gamma_5&1-\beta_5&1-\alpha_1-\beta_4&2\beta_5-1&0\end{array}\right),\] where $\alpha_1\neq 1-2\beta_4$,
\[A_{15}=\left(
\begin{array}{ccccccccc}
1-2\beta_4 & 1-3\beta_5& 0&0&\alpha_5&1&0&-1&0\\
0 & 0& 0&\beta_4&\beta_5&0&0&0&0\\
\gamma_1 & \gamma_2& -1+2\beta_4&\gamma_4&\gamma_5&1-\beta_5&\beta_4&2\beta_5-1&0\end{array}\right),\]
\[A_{16}=\left(
\begin{array}{ccccccccc}
\alpha_1 & \alpha_2& 1&0&\alpha_5&1&-1&-1&0\\
\beta_1 & \beta_2& -1&\beta_4&\beta_5&-1&1&1&0\\
\gamma_1 & \gamma_2& -\alpha_1-\beta_2&\gamma_4&\gamma_5&1-\beta_5&1-\alpha_1-\beta_4&-\alpha_2-\beta_5&0\end{array}\right),\]
\[A_{17}=\left(
\begin{array}{ccccccccc}
\alpha_1 & \alpha_2& \alpha_3&\alpha_4&\alpha_5&0&1+\alpha_3-\beta_8&\alpha_8&0\\
\beta_1 & \beta_2& \beta_3&0&0&1&\beta_7&\beta_8&0\\
\gamma_1 & \gamma_2& -\alpha_1-\beta_2&\gamma_4&\gamma_5&1-\alpha_4&1-\alpha_1&-\alpha_2&-\alpha_3-1\end{array}\right),\] where $\beta_8\neq -1$,
\[A_{18}=\left(
\begin{array}{ccccccccc}
\alpha_1 & \alpha_2& \alpha_3&\alpha_4&0&0&2+\alpha_3&\alpha_8&0\\
\beta_1 & \beta_2& \beta_3&0&\beta_5&1&\beta_7&-1&0\\
\gamma_1 & \gamma_2& -\alpha_1-\beta_2&\gamma_4&\gamma_5&1-\alpha_4-\beta_5&1-\alpha_1&-\alpha_2-\beta_5&-\alpha_3-1\end{array}\right),\] where $\alpha_8\neq 0$,
\[A_{19}=\left(
\begin{array}{ccccccccc}
\alpha_1 & \alpha_2& \alpha_3&0&\alpha_5&0&2+\alpha_3&0&0\\
\beta_1 & \beta_2& \beta_3&0&\beta_5&1&\beta_7&-1&0\\
\gamma_1 & \gamma_2& -\alpha_1-\beta_2&\gamma_4&\gamma_5&1-\beta_5&1-\alpha_1&-\alpha_2-\beta_5&-\alpha_3-1\end{array}\right),\] where $\alpha_3\neq -2$,
\[A_{20}=\left(
\begin{array}{ccccccccc}
\alpha_1 & \alpha_2& -2&\alpha_4&\alpha_5&0&0&0&0\\
\beta_1 & 0& \beta_3&0&\beta_5&1&\beta_7&-1&0\\
\gamma_1 & \gamma_2& -\alpha_1&\gamma_4&\gamma_5&1-\alpha_4-\beta_5&1-\alpha_1&-\alpha_2-\beta_5&1\end{array}\right),\] where $\beta_7(\beta_3-1)\neq 0$,
\[A_{21}=\left(
\begin{array}{ccccccccc}
\alpha_1 & 0& -2&\alpha_4&\alpha_5&0&0&0&0\\
\beta_1 & \beta_2& 1&0&\beta_5&1&\beta_7&-1&0\\
\gamma_1 & \gamma_2& -\alpha_1-\beta_2&\gamma_4&\gamma_5&1-\alpha_4-\beta_5&1-\alpha_1&-\beta_5&1\end{array}\right),\] where $\beta_7\neq 0$,
\[A_{22}=\left(
\begin{array}{ccccccccc}
\alpha_1 & \alpha_2& -2&\alpha_4&\alpha_5&0&0&0&0\\
\beta_1 & \beta_2& \beta_3&0&-\alpha_2&1&0&-1&0\\
\gamma_1 & \gamma_2& -\alpha_1-\beta_2&\gamma_4&\gamma_5&1-\alpha_4+\alpha_2&1-\alpha_1&0&1\end{array}\right),\] 
\[A_{23}=\left(
\begin{array}{ccccccccc}
\alpha_1 & \alpha_2& 1&\alpha_4&0&0&\alpha_7&\alpha_8&0\\
\beta_1 & -\alpha_1& \beta_3&\beta_4&\beta_5&0&\beta_7&1-\alpha_7&0\\
\gamma_1 & \gamma_2& 0&\gamma_4&\gamma_5&1-\alpha_4-\beta_5&1-\alpha_1-\beta_4&-\alpha_2-\beta_5&-1\end{array}\right),\] where $\alpha_8\neq 0$,
\[A_{24}=\left(
\begin{array}{ccccccccc}
\alpha_1 & \alpha_2& 1&0&\alpha_5&0&\alpha_7&0&0\\
\beta_1 & -\alpha_1& \beta_3&\beta_4&\beta_5&0&\beta_7&1-\alpha_7&0\\
\gamma_1 & \gamma_2& 0&\gamma_4&\gamma_5&1-\beta_5&1-\alpha_1-\beta_4&-\alpha_2-\beta_5&-1\end{array}\right),\] where $\alpha_7\neq 0$,
\[A_{25}=\left(
\begin{array}{ccccccccc}
\alpha_1 & \alpha_2& 1&\alpha_4&\alpha_5&0&0&0&0\\
\beta_1 & -\alpha_1& \beta_3&\beta_4&0&0&\beta_7&1&0\\
\gamma_1 & \gamma_2& 0&\gamma_4&\gamma_5&1-\alpha_4&1-\alpha_1-\beta_4&-\alpha_2&-1\end{array}\right),\]

\[A_{26}=\left(
\begin{array}{ccccccccc}
\alpha_1 & 0& 0&\alpha_4&\alpha_5&0&-\beta_8&\alpha_8&0\\
\beta_1 & -\alpha_1& 1&\beta_4&\beta_5&0&\beta_7&\beta_8&0\\
\gamma_1 & \gamma_2& 0&\gamma_4&\gamma_5&1-\alpha_4-\beta_5&1-\alpha_1-\beta_4&-\beta_5&0\end{array}\right),\] where $\alpha_8\neq 0$,
\[A_{27}=\left(
\begin{array}{ccccccccc}
0 & \alpha_2& 0&\alpha_4&\alpha_5&0&-\beta_8&0&0\\
\beta_1 & 0& 1&\beta_4&\beta_5&0&\beta_7&\beta_8&0\\
\gamma_1 & \gamma_2& 0&\gamma_4&\gamma_5&1-\alpha_4-\beta_5&1-\beta_4&-\beta_5&0\end{array}\right),\] where $\alpha_7=-\beta_8\neq 0$,
\[A_{28}=\left(
\begin{array}{ccccccccc}
\alpha_1 & \alpha_2& 0&\alpha_4&\alpha_5&0&0&0&0\\
\beta_1 & 0& 1&\beta_4&\beta_5&0&\beta_7&0&0\\
\gamma_1 & \gamma_2& -\alpha_1&\gamma_4&\gamma_5&1-\alpha_4-\beta_5&1-\alpha_1-\beta_4&-\alpha_2-\beta_5&0\end{array}\right),\]
\[A_{29}=\left(
\begin{array}{ccccccccc}
\alpha_1 & 0& 0&\alpha_4&0&0&\alpha_7&1&0\\
\beta_1 & \beta_2& 0&\beta_4&\beta_5&0&\beta_7&-\alpha_7&0\\
\gamma_1 & \gamma_2& -\alpha_1-\beta_2&\gamma_4&\gamma_5&1-\alpha_4-\beta_5&1-\alpha_1-\beta_4&-\beta_5&0\end{array}\right),\]
\[A_{30}=\left(
\begin{array}{ccccccccc}
0 & \alpha_2& 0&0&\alpha_5&0&1&0&0\\
\beta_1 & \beta_2& 0&\beta_4&\beta_5&0&\beta_7&-1&0\\
\gamma_1 & \gamma_2& -\beta_2&\gamma_4&\gamma_5&1-\beta_5&1-\beta_4&-\alpha_2-\beta_5&0\end{array}\right),\] 
\[A_{31}=\left(
\begin{array}{ccccccccc}
\alpha_1 & \alpha_2& 0&\alpha_4&\alpha_5&0&0&0&0\\
0 & \beta_2& 0&0&\beta_5&0&1&0&0\\
\gamma_1 & \gamma_2& -\alpha_1-\beta_2&\gamma_4&\gamma_5&1-\alpha_4-\beta_5&1-\alpha_1&-\alpha_2-\beta_5&0\end{array}\right),\] 
\[A_{32}=\left(
\begin{array}{ccccccccc}
\alpha_1 & \alpha_2& 0&\alpha_4&\alpha_5&0&0&0&0\\
\beta_1 & \beta_2& 0&\beta_4&\beta_5&0&0&0&0\\
0 & 0& -\alpha_1-\beta_2&0&0&1-\alpha_4-\beta_5&1-\alpha_1-\beta_4&-\alpha_2-\beta_5&0\end{array}\right),\]
\[A_{33}=\left(
\begin{array}{ccccccccc}
\alpha_1 & \alpha_2& 0&\alpha_4&\alpha_5&0&0&0&0\\
\beta_1 & \beta_2& 0&\beta_4&\beta_5&0&0&0&0\\
1 & 0& -\alpha_1-\beta_2&0&0&1-\alpha_4-\beta_5&1-\alpha_1-\beta_4&-\alpha_2-\beta_5&0\end{array}\right),\] where $rk(M_{(1,2,3,4)})=rk(M_{(2,3,4)})$,
\[A_{34}=\left(
\begin{array}{ccccccccc}
\alpha_1 & \alpha_2& 0&\alpha_4&\alpha_5&0&0&0&0\\
\beta_1 & \beta_2& 0&\beta_4&\beta_5&0&0&0&0\\
0 & 1& -\alpha_1-\beta_2&0&0&1-\alpha_4-\beta_5&1-\alpha_1-\beta_4&-\alpha_2-\beta_5&0\end{array}\right),\] where $rk(M_{(1,2,3,4)})=rk(M_{(1,3,4)})$, $rk(M_{(3,4)})=rk(M_{(2,3,4)})$,
\[A_{35}=\left(
\begin{array}{ccccccccc}
\alpha_1 & \alpha_2& 0&\alpha_4&\alpha_5&0&0&0&0\\
\beta_1 & \beta_2& 0&\beta_4&\beta_5&0&0&0&0\\
0 & 0& -\alpha_1-\beta_2&1&0&1-\alpha_4-\beta_5&1-\alpha_1-\beta_4&-\alpha_2-\beta_5&0\end{array}\right),\] where $rk(M_{(1,2,3,4)})=rk(M_{(1,2,4)})$, $rk(M_{(2,4)})=rk(M_{(2,3,4)})$, $rk(M_{(1,4)})=rk(M_{(1,3,4)})$,
\[A_{36}=\left(
\begin{array}{ccccccccc}
\alpha_1 & \alpha_2& 0&\alpha_4&\alpha_5&0&0&0&0\\
\beta_1 & \beta_2& 0&\beta_4&\beta_5&0&0&0&0\\
0 & 0& -\alpha_1-\beta_2&0&1&1-\alpha_4-\beta_5&1-\alpha_1-\beta_4&-\alpha_2-\beta_5&0\end{array}\right),\] where $rk(M_{(1,2,3,4)})=rk(M_{(1,2,3)})$, $rk(M_{(2,3)})=rk(M_{(2,3,4)})$, $rk(M_{(1,3)})=rk(M_{(1,3,4)})$, $rk(M_{(1,2)})=rk(M_{(1,2,4)})$,
\[A_{37}=\left(
\begin{array}{ccccccccc}
\alpha_1 & \alpha_2& 0&\alpha_4&\alpha_5&0&0&0&0\\
\beta_1 & \beta_2& 0&\beta_4&\beta_5&0&0&0&0\\
1 & \gamma_2& -\alpha_1-\beta_2&0&0&1-\alpha_4-\beta_5&1-\alpha_1-\beta_4&-\alpha_2-\beta_5&0\end{array}\right),\] where $rk(M_{(1,2,3,4)})=rk(M_{(2-\gamma_2\times1,3,4)})$, $rk(M_{(3,4)})=rk(M_{(2,3,4)})$, $rk(M_{(3,4)})=rk(M_{(1,3,4)})$, $rk(M_{(2-\gamma_2\times 1,4)})=rk(M_{(1,2,4)})$, 
$rk(M_{(2-\gamma_2\times 1,3)})=rk(M_{(1,2,3)})$ and $\gamma_2\neq 0$,
\[A_{38}=\left(
\begin{array}{ccccccccc}
\alpha_1 & \alpha_2& 0&\alpha_4&\alpha_5&0&0&0&0\\
\beta_1 & \beta_2& 0&\beta_4&\beta_5&0&0&0&0\\
1 & 0& -\alpha_1-\beta_2&\gamma_4&0&1-\alpha_4-\beta_5&1-\alpha_1-\beta_4&-\alpha_2-\beta_5&0\end{array}\right),\] where $rk(M_{(1,2,3,4)})=rk(M_{(2,3-\gamma_4\times 1,4)})$, $rk(M_{(2,4)})=rk(M_{(2,3,4)})$, $rk(M_{(3-\gamma_4\times 1,4)})=rk(M_{(1,3,4)})$, $rk(M_{(2,4)})=rk(M_{(1,2,4)})$, 
$rk(M_{(2,3-\gamma_4\times 1)})=rk(M_{(1,2,3)})$, $rk(M_{(4)})=rk(M_{(3,4)})$ and $\gamma_4\neq 0$,
\[A_{39}=\left(
\begin{array}{ccccccccc}
\alpha_1 & \alpha_2& 0&\alpha_4&\alpha_5&0&0&0&0\\
\beta_1 & \beta_2& 0&\beta_4&\beta_5&0&0&0&0\\
0 & 1& -\alpha_1-\beta_2&\gamma_4&0&1-\alpha_4-\beta_5&1-\alpha_1-\beta_4&-\alpha_2-\beta_5&0\end{array}\right),\] where $rk(M_{(1,2,3,4)})=rk(M_{(1,3-\gamma_4\times 2,4)})$, $rk(M_{(3-\gamma_4\times 2,4)})=rk(M_{(2,3,4)})$, $rk(M_{(1,4)})=rk(M_{(1,3,4)})$, $rk(M_{(1,4)})=rk(M_{(1,2,4)})$, 
$rk(M_{(1,3-\gamma_4\times 2)})=rk(M_{(1,2,3)})$, $rk(M_{(4)})=rk(M_{(3,4)})$,  $rk(M_{(4)})=rk(M_{(2,4)})$ and $\gamma_4\neq 0$,
\[A_{40}=\left(
\begin{array}{ccccccccc}
\alpha_1 & \alpha_2& 0&\alpha_4&\alpha_5&0&0&0&0\\
\beta_1 & \beta_2& 0&\beta_4&\beta_5&0&0&0&0\\
1 & 0& -\alpha_1-\beta_2&0&\gamma_5&1-\alpha_4-\beta_5&1-\alpha_1-\beta_4&-\alpha_2-\beta_5&0\end{array}\right),\] where $rk(M_{(1,2,3,4)})=rk(M_{(2,3,4-\gamma_5\times 1)})$, $rk(M_{(2,3)})=rk(M_{(2,3,4)})$, $rk(M_{(3,4-\gamma_5\times 1)})=rk(M_{(1,3,4})$, $rk(M_{(2,4-\gamma_5\times 1)})=rk(M_{(1,2,4)})$, 
$rk(M_{(2,3)})=rk(M_{(1,2,3)})$, $rk(M_{(3)})=rk(M_{(3,4)})$, $rk(M_{(2)})=rk(M_{(2,4)})$, $rk(M_{(4-\gamma_5\times 1)})=rk(M_{(1,4)})$ and $\gamma_5\neq 0$,
\[A_{41}=\left(
\begin{array}{ccccccccc}
\alpha_1 & \alpha_2& 0&\alpha_4&\alpha_5&0&0&0&0\\
\beta_1 & \beta_2& 0&\beta_4&\beta_5&0&0&0&0\\
0 & 1& -\alpha_1-\beta_2&0&\gamma_5&1-\alpha_4-\beta_5&1-\alpha_1-\beta_4&-\alpha_2-\beta_5&0\end{array}\right),\] where $rk(M_{(1,2,3,4)})=rk(M_{(1,3,4-\gamma_5\times 2)})$, $rk(M_{(3,4-\gamma_5\times 2)})=rk(M_{(2,3,4)})$, $rk(M_{(1,3)})=rk(M_{(1,3,4})$, $rk(M_{(1,4-\gamma_5\times 2)})=rk(M_{(1,2,4)})$, 
$rk(M_{(1,3)})=rk(M_{(1,2,3)})$, $rk(M_{(3)})=rk(M_{(3,4)})$, $rk(M_{(4-\gamma_5\times 2)})=rk(M_{(2,4)})$, $rk(M_{(1)})=rk(M_{(1,4)})$, $rk(M_{(3)})=rk(M_{(2,3)})$ and $\gamma_5\neq 0$,
\[A_{42}=\left(
\begin{array}{ccccccccc}
\alpha_1 & \alpha_2& 0&\alpha_4&\alpha_5&0&0&0&0\\
\beta_1 & \beta_2& 0&\beta_4&\beta_5&0&0&0&0\\
0 & 0& -\alpha_1-\beta_2&1&\gamma_5&1-\alpha_4-\beta_5&1-\alpha_1-\beta_4&-\alpha_2-\beta_5&0\end{array}\right),\] where $rk(M_{(1,2,3,4)})=rk(M_{(1,2,4-\gamma_5\times 3)})$, $rk(M_{(2,4-\gamma_5\times 3)})=rk(M_{(2,3,4)})$, $rk(M_{(1,4-\gamma_5\times 3)})=rk(M_{(1,3,4})$, $rk(M_{(1,2)})=rk(M_{(1,2,4)})$, 
$rk(M_{(1,2)})=rk(M_{(1,2,3)})$, $rk(M_{(4-\gamma_5\times 3)})=rk(M_{(3,4)})$, $rk(M_{(2)})=rk(M_{(2,4)})$, $rk(M_{(1)})=rk(M_{(1,4)})$, $rk(M_{(1)})=rk(M_{(1,3)})$ and $\gamma_5\neq 0$,
 \[A_{43}=\left(
\begin{array}{ccccccccc}
\alpha_1 & \alpha_2& 0&\alpha_4&\alpha_5&0&0&0&0\\
\beta_1 & \beta_2& 0&\beta_4&\beta_5&0&0&0&0\\
0 & 1& -\alpha_1-\beta_2&\gamma_4&\gamma_5&1-\alpha_4-\beta_5&1-\alpha_1-\beta_4&-\alpha_2-\beta_5&0\end{array}\right),\] where $rk(M_{(1,2,3,4)})=rk(M_{(1,3-\gamma_4\times 2,4-\gamma_5\times 2)})$, $rk(M_{(3-\gamma_4\times 2, 4-\gamma_5\times 2)})=rk(M_{(2,3,4)})$, $rk(M_{(1,\gamma_4\times 4-\gamma_5\times 3)})=rk(M_{(1,3,4})$, $rk(M_{(1,4-\gamma_5\times 2)})=rk(M_{(1,2,4)})$, 
$rk(M_{(1,3-\gamma_4\times 2)})=rk(M_{(1,2,3)})$, $rk(M_{(\gamma_4\times 4-\gamma_5\times 3)})=rk(M_{(3,4)})$, $rk(M_{(4-\gamma_5\times 2)})=rk(M_{(2,4)})$, $rk(M_{(1)})=rk(M_{(1,4)})$, $rk(M_{(1)})=rk(M_{(1,3)})$, $rk(M_{(1)})=rk(M_{(1,2)})$ and $\gamma_4\neq 0$, $\gamma_5\neq 0$,
\[A_{44}=\left(
\begin{array}{ccccccccc}
\alpha_1 & \alpha_2& 0&\alpha_4&\alpha_5&0&0&0&0\\
\beta_1 & \beta_2& 0&\beta_4&\beta_5&0&0&0&0\\
1 & 0& -\alpha_1-\beta_2&\gamma_4&\gamma_5&1-\alpha_4-\beta_5&1-\alpha_1-\beta_4&-\alpha_2-\beta_5&0\end{array}\right),\] where $rk(M_{(1,2,3,4)})=rk(M_{(2,3-\gamma_4\times 1,4-\gamma_5\times 1)})$, $rk(M_{(2,\gamma_4\times 3-\gamma_5\times 3)})=rk(M_{(2,3,4)})$, $rk(M_{(3-\gamma_4\times 13,4-\gamma_5\times 1)})=rk(M_{(1,3,4})$, $rk(M_{(2,4-\gamma_5\times 1)})=rk(M_{(1,2,4)})$, 
$rk(M_{(2,3-\gamma_4\times 1)})=rk(M_{(1,2,3)})$, $rk(M_{(\gamma_4\times 4-\gamma_5\times 3)})=rk(M_{(3,4)})$, $rk(M_{(2)})=rk(M_{(2,4)})$, $rk(M_{(4-\gamma_4\times 1)})=rk(M_{(1,4)})$, $rk(M_{(3-\gamma_5\times 1)})=rk(M_{(1,3)})$, $rk(M_{(2)})=rk(M_{(1,2)})$, $rk(M_{(1)})=0$ and $\gamma_4\neq 0$, $\gamma_5\neq 0$,
\[A_{45}=\left(
\begin{array}{ccccccccc}
\alpha_1 & \alpha_2& 0&\alpha_4&\alpha_5&0&0&0&0\\
0 & 0& 0&1-3\alpha_1&-\alpha_1-2\alpha_2&0&0&0&0\\
1 & \gamma_2& -\alpha_1&0&\gamma_5&1+\alpha_1+2\alpha_2-\alpha_4&2\alpha_1&\alpha_1+\alpha_2&0\end{array}\right),\] where $\gamma_2\neq 0$, $\gamma_5\neq 0$,
\[A_{46}=\left(
\begin{array}{ccccccccc}
\frac{1}{5} & \alpha_2& 0&\frac{3}{5}+\alpha_2&\alpha_5&0&0&0&0\\
0 & 0& 0&\frac{2}{5}&-\frac{1}{5}-2\alpha_2&0&0&0&0\\
1 & \gamma_2& -\frac{1}{5}&\gamma_4&0&\frac{3}{5}+\alpha_2&\frac{2}{5}&\frac{1}{5}+\alpha_2&0\end{array}\right),\] where $\gamma_2\neq 0$, $\gamma_4\neq 0$,
\[A_{47}=\left(
\begin{array}{ccccccccc}
\frac{1}{5} & -\frac{1}{4}& 0&\frac{7}{20}&0&0&0&0&0\\
0 & 0& 0&\frac{2}{5}&\frac{3}{10}&0&0&0&0\\
1 & \gamma_2& -\frac{1}{5} &\gamma_4&\gamma_5&\frac{7}{20} &\frac{2}{5} &-\frac{1}{20} &0\end{array}\right),\] where 
$\gamma_2\neq 0,\ \gamma_3\neq 0,\ \gamma_4\neq 0.$\end{theorem}

{\bf Proof.} 
 Due to the $3^{rd}$ column in (\ref{E4}) one can conclude:
  
 \underline{\textbf{Case 1.\ $\alpha_9\neq 0$}}. It is possible to make $\alpha'_9=1$ ($c^2\alpha_9=1$), $\alpha'_7=\alpha'_8=0$ to get  
\[A_1=\left(
\begin{array}{ccccccccc}
\alpha_1 & \alpha_2& \alpha_3&\alpha_4&\alpha_5&\alpha_6&0&0&1\\
\beta_1 & \beta_2& \beta_3&\beta_4&\beta_5&\beta_8-\alpha_3&\beta_7&\beta_8&\beta_9\\
\gamma_1 & \gamma_2& -\alpha_1-\beta_2&\gamma_4&\gamma_5&1-\alpha_4-\beta_5&1-\alpha_1-\beta_4&-\alpha_2-\beta_5&-\beta_8\end{array}\right)\cong \]
\[\left(
\begin{array}{ccccccccc}
	\alpha_1 & \alpha_2& -\alpha_3&\alpha_4&\alpha_5&-\alpha_6&0&0&1\\
	\beta_1 & \beta_2& -\beta_3&\beta_4&\beta_5&-\beta_8+\alpha_3&-\beta_7&-\beta_8&\beta_9\\
	-\gamma_1 & -\gamma_2& -\alpha_1-\beta_2&-\gamma_4&-\gamma_5&1-\alpha_4-\beta_5&1-\alpha_1-\beta_4&-\alpha_2-\beta_5&\beta_8\end{array}\right).\] 
	
 \underline{\textbf{Case 2.\ $\alpha_9= 0$}, $\beta_9\neq 0$}. One can make $\beta'_9=1$, $\beta'_7=\beta'_8=0$ to get  
\[A_2=\left(
\begin{array}{ccccccccc}
\alpha_1 & \alpha_2& \alpha_3&\alpha_4&\alpha_5&\alpha_6&\alpha_7&\alpha_8&0\\
\beta_1 & \beta_2&\beta_3&\beta_4&\beta_5&\alpha_7-\alpha_3&0&0&1\\
\gamma_1 & \gamma_2& -\alpha_1-\beta_2&\gamma_4&\gamma_5&1-\alpha_4-\beta_5&1-\alpha_1-\beta_4&-\alpha_2-\beta_5&-\alpha_7\end{array}\right)\cong\]
\[\left(
\begin{array}{ccccccccc}
\alpha_1 & \alpha_2& -\alpha_3&\alpha_4&\alpha_5&-\alpha_6&-\alpha_7&-\alpha_8&0\\
\beta_1 & \beta_2&-\beta_3&\beta_4&\beta_5&-\alpha_7+\alpha_3&0&0&1\\
-\gamma_1 & -\gamma_2& -\alpha_1-\beta_2&-\gamma_4&-\gamma_5&1-\alpha_4-\beta_5&1-\alpha_1-\beta_4&-\alpha_2-\beta_5&\alpha_7\end{array}\right).\]

If $\alpha_9=\beta_9= 0$, $\alpha_6\neq 0$ then due to the $3^{rd}$ column in (\ref{E3}) it is possible making $\alpha'_6=1$ by appropriate choose of $c$ to come to the case

\underline{\textbf{Case 3.\ $\alpha_9=\beta_9= 0$, $\alpha_6= 1 $}} with respect to $g$ of the form $g^{-1}=\left(
 \begin{array}{ccc}
 1 & 0& 0\\ 0 & 1& 0\\ a & b& 1\end{array}\right)$. If $a=-\alpha_4-b\alpha_7$ then $\alpha'_4=0$ and 
  $\mathcal{A}'_1=g(\mathcal{A}_1+b\mathcal{A}_3)g^{-1}=$ 
 \[\left( 
 \begin{array}{ccc}
 \alpha_1-\alpha_4(\alpha_3+\alpha_7)-b\alpha_7(\alpha_3+\alpha_7)&\alpha_2-\alpha_4(\alpha_3+\alpha_8)-b\alpha_7(\alpha_3+\alpha_8)&\alpha_3\\
 \beta_1-\alpha_4(\beta_3+\beta_7)-b\alpha_7(\beta_3+\beta_7)&\beta_2-\alpha_4(\beta_3+\beta_8)-b\alpha_7(\beta_3+\beta_8)&\beta_3\\
 *&*&*\end{array}\right),\]
  $\mathcal{A}'_2=g(\mathcal{A}_2+b\mathcal{A}_3)g^{-1}=$ 
 \[\left( 
 \begin{array}{ccc}
 0&\alpha_5+b(\alpha_8+1)&1\\
 \beta_4-\alpha_4+b(\beta_7-\alpha_7\beta_6)&\beta_5+b(\alpha_7+2\beta_8-\alpha_3)&\alpha_7+\beta_8-\alpha_3\\
 *&*&*\end{array}\right),\]
  \[\mathcal{A}'_3=g\mathcal{A}_3g^{-1}=\left( 
 \begin{array}{ccc}
 \alpha_7&\alpha_8&0\\
 \beta_7&\beta_8&0\\
 *&*&-\alpha_7-\beta_8\end{array}\right)\] which implies that:
  	 
 	1. If $\alpha_8\neq -1$ then one can make $\alpha'_5=0$ to get 
\[A_3=\left(
\begin{array}{ccccccccc}
\alpha_1 & \alpha_2& \alpha_3&0&0&1&\alpha_7&\alpha_8&0\\
\beta_1 & \beta_2& \beta_3&\beta_4&\beta_5&\alpha_7+\beta_8-\alpha_3&\beta_7&\beta_8&0\\
\gamma_1 & \gamma_2& -\alpha_1-\beta_2&\gamma_4&\gamma_5&1-\beta_5&1-\alpha_1-\beta_4&-\alpha_2-\beta_5&-\alpha_7-\beta_8\end{array}\right),\] where $\alpha_8\neq -1$.\\

2. If $\alpha_8=-1$ and $\alpha_7+2\beta_8-\alpha_3\neq 0$ then one can make  $\beta'_5=0$ to get
\[A_4=\left(
\begin{array}{ccccccccc}
\alpha_1 & \alpha_2& \alpha_3&0&\alpha_5&1&\alpha_7&-1&0\\
\beta_1 & \beta_2& \beta_3&\beta_4&0&\alpha_7+\beta_8-\alpha_3&\beta_7&\beta_8&0\\
\gamma_1 & \gamma_2& -\alpha_1-\beta_2&\gamma_4&\gamma_5&1&1-\alpha_1-\beta_4&-\alpha_2&-\alpha_7-\beta_8\end{array}\right),\] where $\alpha_7\neq \alpha_3-2\beta_8$.

3. If $\alpha_8=-1$, $\alpha_7=\alpha_3-2\beta_8$, that is $\beta_6=-\beta_8$, and $\beta_7\neq (2\beta_8-\alpha_3)\beta_8$ then one can make $\beta'_4=0$ to get
\[A_5=\left(
\begin{array}{ccccccccc}
\alpha_1 & \alpha_2& \alpha_3&0&\alpha_5&1&\alpha_3-2\beta_8&-1&0\\
\beta_1 & \beta_2& \beta_3&0&\beta_5&-\beta_8&\beta_7&\beta_8&0\\
\gamma_1 & \gamma_2& -\alpha_1-\beta_2&\gamma_4&\gamma_5&1-\beta_5&1-\alpha_1&-\alpha_2-\beta_5&-\alpha_3+\beta_8\end{array}\right),\]
where $\beta_7\neq (2\beta_8-\alpha_3)\beta_8$.

Note that until now there was no need for the assumption $Ch(\mathbb{F})\neq 2$ and from now one needs it.

4. If $\alpha_8=-1$, $\alpha_7=\alpha_3-2\beta_8$, $\beta_6=-\beta_8$, $\beta_7=(2\beta_8-\alpha_3)\beta_8$ and $\alpha_7(\alpha_3+\alpha_7)=2(\alpha_3-2\beta_8)(\alpha_3-\beta_8)\neq 0$ one can make $\alpha'_1=0$ to get

\[A_6=\left(
\begin{array}{ccccccccc}
0 & \alpha_2& \alpha_3&0&\alpha_5&1&\alpha_3-2\beta_8&-1&0\\
\beta_1 & \beta_2& \beta_3&\beta_4&\beta_5&-\beta_8&(2\beta_8-\alpha_3)\beta_8&\beta_8&0\\
\gamma_1 & \gamma_2& -\beta_2&\gamma_4&\gamma_5&1-\beta_5&1-\beta_4&-\alpha_2-\beta_5&-\alpha_3+\beta_8\end{array}\right),\] where 
$(\alpha_3-2\beta_8)(\alpha_3-\beta_8)\neq 0$.

5. If $\alpha_8=-1$, $\alpha_7=\alpha_3-2\beta_8$, $\beta_6=-\beta_8$, $\beta_7=(2\beta_8-\alpha_3)\beta_8$, $(\alpha_3-2\beta_8)(\alpha_3-\beta_8)= 0$, $(\alpha_3-2\beta_8)(\beta_3+(2\beta_8-\alpha_3)\beta_8)\neq 0$ then  $\alpha_8=-1$, $\alpha_7=-\alpha_3$, $\beta_6=-\alpha_3$, $\beta_8=\alpha_3$, $\beta_7=\alpha^2_3$ and one can make $\beta'_1=0$ to get

\[A_7=\left(
\begin{array}{ccccccccc}
\alpha_1 & \alpha_2& \alpha_3&0&\alpha_5&1&-\alpha_3&-1&0\\
0 & \beta_2& \beta_3&\beta_4&\beta_5&-\alpha_3&\alpha^2_3&\alpha_3&0\\
\gamma_1 & \gamma_2& -\alpha_1-\beta_2&\gamma_4&\gamma_5&1-\beta_5&1-\alpha_1-\beta_4&-\alpha_2-\beta_5&0\end{array}\right),\] where 
$\alpha_3(\alpha^2_3+\beta_3)\neq 0$.

6. If   $\alpha_8=-1$, $\alpha_7=\alpha_3-2\beta_8$, $\beta_6=-\beta_8$, $\beta_7=(2\beta_8-\alpha_3)\beta_8$, $(\alpha_3-2\beta_8)(\alpha_3-\beta_8)= 0$, $(\alpha_3-2\beta_8)(\beta_3+(2\beta_8-\alpha_3)\beta_8)= 0$ and $(\alpha_3-2\beta_8)(\alpha_3-1)\neq 0$ then $\beta_8=\alpha_3$, $\beta_3=-\alpha^2_3$, one can make $\alpha'_2=0$ to get

%$\alpha_8=-1$, $\alpha_7=-\alpha_3$, $\beta_6=-\alpha_3$, $\beta_8=\alpha_3$, $\beta_7=\alpha^2_3$, $\alpha_3(\alpha^2_3+\beta_3)= 0$ and $(\alpha_3-2\beta_8)(\alpha_3-1)\neq 0$ then $\beta_3=-\alpha^2_3$, one can make $\alpha'_2=0$ to get

%$\alpha_6+\alpha_8=0$, $\alpha_7+2\beta_8-\alpha_3= 0$, $\beta_7-\alpha_7\beta_6= 0$, $2(\alpha_3-2\beta_8)(\alpha_3-\beta_8)= 0$ $(\alpha_3-2\beta_8)(\beta_3+(2\beta_8-\alpha_3)\beta_8)= 0$ and $(\alpha_3-2\beta_8)(\alpha_3-1)\neq 0$ ( that is $\alpha_8=-1$, $\alpha_7=-\alpha_3$, $\beta_6=-\alpha_3$, $\beta_7=\alpha^2_3$, $\beta_3=-\alpha^2_3$ ) one can make $\alpha'_2=0$ to get

\[A_8=\left(
\begin{array}{ccccccccc}
\alpha_1 & 0& \alpha_3&0&\alpha_5&1&-\alpha_3&-1&0\\
\beta_1 & \beta_2& -\alpha^2_3&\beta_4&\beta_5&-\alpha_3&\alpha^2_3&\alpha_3&0\\
\gamma_1 & \gamma_2& -\alpha_1-\beta_2&\gamma_4&\gamma_5&1-\beta_5&1-\alpha_1-\beta_4&-\beta_5&0\end{array}\right),\] where 
$\alpha_3\neq 0, 1$.

7. In  $\alpha_8=-1$, $\alpha_7=\alpha_3-2\beta_8$, $\beta_6=-\beta_8$, $\beta_7=(2\beta_8-\alpha_3)\beta_8$, $(\alpha_3-2\beta_8)(\alpha_3-\beta_8)= 0$, $(\alpha_3-2\beta_8)(\beta_3+(2\beta_8-\alpha_3)\beta_8)= 0$, $(\alpha_3-2\beta_8)(\alpha_3-1)= 0$ case one has the following two possibilities. 

a) $\alpha_3=2\beta_8$ that is $\alpha_7=0$.

1. If $\beta_3\neq 0$ one can make $\gamma'_3=0$ to get 
\[A_{9}=\left(
\begin{array}{ccccccccc}
\alpha_1 & \alpha_2& 2\beta_8&0&\alpha_5&1&0&-1&0\\
\beta_1 & -\alpha_1& \beta_3&\beta_4&\beta_5&-\beta_8&0&\beta_8&0\\
\gamma_1 & \gamma_2& 0&\gamma_4&\gamma_5&1-\beta_5&1-\alpha_1-\beta_4&-\alpha_2-\beta_5&-\beta_8\end{array}\right),\] where 
$\beta_3\neq 0$.
 
2. If $\beta_3= 0$  and $\beta_2\neq \alpha_4\beta_8$ one can make $\gamma'_2=0$ to get 
\[A_{10}=\left(
\begin{array}{ccccccccc}
\alpha_1 & \alpha_2& 2\beta_8&0&\alpha_5&1&0&-1&0\\
\beta_1 & \beta_2& 0&\beta_4&\beta_5&-\beta_8&0&\beta_8&0\\
\gamma_1 & 0& -\alpha_1-\beta_2&\gamma_4&\gamma_5&1-\beta_5&1-\alpha_1-\beta_4&-\alpha_2-\beta_5&-\beta_8\end{array}\right),\] where 
$\beta_2\neq 0$.

3. If $\beta_3= 0$, $\beta_2= \alpha_4\beta_8$ and $\beta_1\neq 0$ one can make $\gamma'_1=0$ to get 
\[A_{11}=\left(
\begin{array}{ccccccccc}
\alpha_1 & \alpha_2& 2\beta_8&0&\alpha_5&1&0&-1&0\\
\beta_1 & 0& 0&\beta_4&\beta_5&-\beta_8&0&\beta_8&0\\
0 & \gamma_2& -\alpha_1&\gamma_4&\gamma_5&1-\beta_5&1-\alpha_1-\beta_4&-\alpha_2-\beta_5&-\beta_8\end{array}\right),\] where 
$\beta_1\neq 0$.

4. If $\beta_3= 0$, $\beta_2= \alpha_4\beta_8$, $\beta_1= 0$ and $2\beta_8\neq 0$ one can make $\gamma'_8=0$ to get 
\[A_{12}=\left(
\begin{array}{ccccccccc}
\alpha_1 & \alpha_2& 2\beta_8&0&\alpha_5&1&0&-1&0\\
0 & 0& 0&\beta_4&-\alpha_2&-\beta_8&0&\beta_8&0\\
\gamma_1 & \gamma_2& -\alpha_1&\gamma_4&\gamma_5&1+\alpha_2&1-\alpha_1-\beta_4&0&-\beta_8\end{array}\right),\] where 
$\beta_8\neq 0$.

% 5. If $\beta_3= 0$, $\beta_2= \alpha_4\beta_8$, $\beta_1= 0$, $\beta_8= 0$ and $\alpha_3\neq 0$ one can make $\gamma'_6=0$ to get 
% \[A_{13}=\left(
 %\begin{array}{ccccccccc}
 %\alpha_1 & \alpha_2& \alpha_3&0&\alpha_5&1&0&-1&0\\
 %0 & 0& 0&\beta_4&\beta_5&0&0&0&0\\
 %\gamma_1 & \gamma_2& -\alpha_1&\gamma_4&\gamma_5&1-\beta_5&1-\alpha_1-\beta_4&-\alpha_2-\beta_5&-\alpha_3\end{array}\right),\] where 
 %$\alpha_3\neq 0$.
 
 5. If $\beta_3= 0$, $\beta_2= \alpha_4\beta_8$, $\beta_1= 0$, $\beta_8= 0$ and $\alpha_4\neq 1-\alpha_2-3\beta_5$ one can make $\gamma'_5=0$ to get 
\[A_{13}=\left(
\begin{array}{ccccccccc}
\alpha_1 & \alpha_2& 0&0&\alpha_5&1&0&-1&0\\
0 & 0& 0&\beta_4&\beta_5&0&0&0&0\\
\gamma_1 & \gamma_2& -\alpha_1&\gamma_4&0&1-\beta_5&1-\alpha_1-\beta_4&-\alpha_2-\beta_5&0\end{array}\right),\] where 
$1-\alpha_2-3\beta_5\neq 0$.
 
7. If $\beta_3= 0$, $\beta_2= \alpha_4\beta_8$, $\beta_1= 0$, $\beta_8= 0$, $\alpha_3= 0$, $\alpha_4= 1-\alpha_2-3\beta_5$  and $\alpha_1\neq 1-2\beta_4$ one can make $\gamma'_4=0$ to get 
\[A_{14}=\left(
\begin{array}{ccccccccc}
\alpha_1 & 1-3\beta_5& 0&0&\alpha_5&1&0&-1&0\\
0 & 0& 0&\beta_4&\beta_5&0&0&0&0\\
\gamma_1 & \gamma_2& -\alpha_1&0&\gamma_5&1-\beta_5&1-\alpha_1-\beta_4&2\beta_5-1&0\end{array}\right),\] where $\alpha_1\neq 1-2\beta_4$.

8. If $\beta_3= 0$, $\beta_2= \alpha_4\beta_8$, $\beta_1= 0$, $\beta_8= 0$, $\alpha_3= 0$, $\alpha_4= 1-\alpha_2-3\beta_5$, $\alpha_1= 1-2\beta_4$ one gets 
\[A_{15}=\left(
\begin{array}{ccccccccc}
1-2\beta_4 & 1-3\beta_5& 0&0&\alpha_5&1&0&-1&0\\
0 & 0& 0&\beta_4&\beta_5&0&0&0&0\\
\gamma_1 & \gamma_2& -1+2\beta_4&\gamma_4&\gamma_5&1-\beta_5&\beta_4&2\beta_5-1&0\end{array}\right).\]

%b) $\alpha_3\neq 2\beta_8$. In this case $\alpha_3=\beta_8=1$, $\beta_3=-1$, $\alpha_7=\beta_6=-1$, $\beta_7=1$ and therefore

%\[A_{17}=\left(
%\begin{array}{ccccccccc}
%1-2\beta_4 & 1-3\beta_5& 0&0&\alpha_5&1&0&-1&0\\
%0 & 0& 0&\beta_4&\beta_5&0&0&0&0\\
%\gamma_1 & \gamma_2& -1+2\beta_4&0&\gamma_5&1-\beta_5&\beta_4&2\beta_5-1&0\end{array}\right),\] where 
%$1-\alpha_2-3\beta_5\neq 0$.

b) $\alpha_3\neq 2\beta_8$. In this case $\alpha_3=\beta_8=1$, $\beta_3=-1$, $\alpha_7=\beta_6=-1$, $\beta_7=1$ and therefore

\[A_{16}=\left(
\begin{array}{ccccccccc}
\alpha_1 & \alpha_2& 1&0&\alpha_5&1&-1&-1&0\\
\beta_1 & \beta_2& -1&\beta_4&\beta_5&-1&1&1&0\\
\gamma_1 & \gamma_2& -\alpha_1-\beta_2&\gamma_4&\gamma_5&1-\beta_5&1-\alpha_1-\beta_4&-\alpha_2-\beta_5&0\end{array}\right).\]
%\end{document}

If $\alpha_9=\beta_9=\alpha_6= 0$, $\alpha_7+\beta_8-\alpha_3\neq 0$ then it is possible making $\beta'_6=1$ by appropriate choose of $c$ to come to the case:
  
 \underline{\textbf{Case 4.\ $\alpha_9=\beta_9=\alpha_6= 0$, $\alpha_7+\beta_8-\alpha_3= 1$},} with respect to $g$ with the $c=1$ and if $a=-\beta_4-b\beta_7$ then $\beta'_4=0$,  
  $\mathcal{A}'_1=g(\mathcal{A}_1+b\mathcal{A}_3)g^{-1}=$ 
  \[\left( 
  \begin{array}{ccc}
  \alpha_1-\beta_4(\alpha_3+\alpha_7)-b\beta_7(\alpha_3+\alpha_7)&\alpha_2-\beta_4(\alpha_3+\alpha_8)-b\beta_7(\alpha_3+\alpha_8)&\alpha_3\\
  \beta_1-\beta_4(\beta_3+\beta_7)-b\beta_7(\beta_3+\beta_7)&\beta_2-\beta_4(\beta_3+\beta_8)-b\beta_7(\beta_3+\beta_8)&\beta_3\\
  *&*&*\end{array}\right),\]
  \[\mathcal{A}'_2=g(\mathcal{A}_2+b\mathcal{A}_3)g^{-1}= 
  \left( 
  \begin{array}{ccc}
  \alpha_4+b\alpha_7&\alpha_5+b\alpha_8&0\\
  0&\beta_5+b(1+\beta_8)&1\\
  *&*&*\end{array}\right),\]
  \[\mathcal{A}'_3=g\mathcal{A}_3g^{-1}=\left( 
  \begin{array}{ccc}
  \alpha_7&\alpha_8&0\\
  \beta_7&\beta_8&0\\
  *&*&-1-\alpha_3\end{array}\right).\] 
  Therefore:\\
1. If $\beta_8\neq -1$  one can make  $\beta'_5=0$ to get
\[A_{17}=\left(
\begin{array}{ccccccccc}
\alpha_1 & \alpha_2& \alpha_3&\alpha_4&\alpha_5&0&1+\alpha_3-\beta_8&\alpha_8&0\\
\beta_1 & \beta_2& \beta_3&0&0&1&\beta_7&\beta_8&0\\
\gamma_1 & \gamma_2& -\alpha_1-\beta_2&\gamma_4&\gamma_5&1-\alpha_4&1-\alpha_1&-\alpha_2&-\alpha_3-1\end{array}\right),\] where $\beta_8\neq -1$.

2. If $\beta_8=-1$ and $\alpha_8\neq 0$ one can make  $\alpha'_5=0$ to get
\[A_{18}=\left(
\begin{array}{ccccccccc}
\alpha_1 & \alpha_2& \alpha_3&\alpha_4&0&0&2+\alpha_3&\alpha_8&0\\
\beta_1 & \beta_2& \beta_3&0&\beta_5&1&\beta_7&-1&0\\
\gamma_1 & \gamma_2& -\alpha_1-\beta_2&\gamma_4&\gamma_5&1-\alpha_4-\beta_5&1-\alpha_1&-\alpha_2-\beta_5&-\alpha_3-1\end{array}\right),\] where $\alpha_8\neq 0$.

3. If $\beta_8=-1$, $\alpha_8= 0$ and $\alpha_7=2+\alpha_3\neq 0$ one can make  $\alpha'_4=0$ to get
\[A_{19}=\left(
\begin{array}{ccccccccc}
\alpha_1 & \alpha_2& \alpha_3&0&\alpha_5&0&2+\alpha_3&0&0\\
\beta_1 & \beta_2& \beta_3&0&\beta_5&1&\beta_7&-1&0\\
\gamma_1 & \gamma_2& -\alpha_1-\beta_2&\gamma_4&\gamma_5&1-\beta_5&1-\alpha_1&-\alpha_2-\beta_5&-\alpha_3-1\end{array}\right),\] where $\alpha_3\neq -2$.

4. If $\beta_8=-1$, $\alpha_8= 0$, $\alpha_3=-2$ and $\beta_7(\beta_3-1)\neq 0$ one can make  $\beta'_2=0$ to get
\[A_{20}=\left(
\begin{array}{ccccccccc}
\alpha_1 & \alpha_2& -2&\alpha_4&\alpha_5&0&0&0&0\\
\beta_1 & 0& \beta_3&0&\beta_5&1&\beta_7&-1&0\\
\gamma_1 & \gamma_2& -\alpha_1&\gamma_4&\gamma_5&1-\alpha_4-\beta_5&1-\alpha_1&-\alpha_2-\beta_5&1\end{array}\right),\] where $\beta_7(\beta_3-1)\neq 0$.

5. If $\beta_8=-1$, $\alpha_8= 0$, $\alpha_3=-2$, $\beta_7(\beta_3-1)= 0$, $\beta_7(\alpha_3+\alpha_8)\neq 0$ then $\beta_3=1$, $\beta_7\neq 0$ and one can make  $\alpha'_2=0$ to get
\[A_{21}=\left(
\begin{array}{ccccccccc}
\alpha_1 & 0& -2&\alpha_4&\alpha_5&0&0&0&0\\
\beta_1 & \beta_2& 1&0&\beta_5&1&\beta_7&-1&0\\
\gamma_1 & \gamma_2& -\alpha_1-\beta_2&\gamma_4&\gamma_5&1-\alpha_4-\beta_5&1-\alpha_1&-\beta_5&1\end{array}\right),\] where $\beta_7\neq 0$.

6. If $\beta_8=-1$, $\alpha_8= 0$, $\alpha_3=-2$, $\beta_7(\beta_3-1)= 0$, $\beta_7(\alpha_3+\alpha_8)= 0$ then $\beta_7= 0$ and one can make $\gamma
'_8=0$ to get
\[A_{22}=\left(
\begin{array}{ccccccccc}
\alpha_1 & \alpha_2& -2&\alpha_4&\alpha_5&0&0&0&0\\
\beta_1 & \beta_2& \beta_3&0&-\alpha_2&1&0&-1&0\\
\gamma_1 & \gamma_2& -\alpha_1-\beta_2&\gamma_4&\gamma_5&1-\alpha_4+\alpha_2&1-\alpha_1&0&1\end{array}\right).\] 

If $\alpha_9=\beta_9=\alpha_6= 0$, $\alpha_7+\beta_8-\alpha_3=0$, $\alpha_3\neq 0$ one can make $\alpha'_3=1$ to come to the case:

\underline{\textbf{Case 5.\ $\alpha_9=\beta_9=\alpha_6= 0$, $\alpha_7+\beta_8-\alpha_3=0$, $\alpha_3= 1$}} with respect to $g$ with the $c=1$. If $a=\frac{-1}{2}(\alpha_1+\beta_2+b\beta_3)$ then $\gamma'_3=0$ and 
$\mathcal{A}'_1=g(\mathcal{A}_1+a\mathcal{A}_3)g^{-1}=$ 
\[ \left(\begin{array}{ccc}
\alpha_1-\frac{1}{2}(\alpha_1+\beta_2)(\alpha_3+\alpha_7)-\frac{b}{2}\beta_3(\alpha_3+\alpha_7)&\alpha_2-\frac{1}{2}(\alpha_1+\beta_2)(\alpha_3+\alpha_8)-\frac{b}{2}\beta_3(\alpha_3+\alpha_8)&1\\
\beta_1-\frac{1}{2}(\alpha_1+\beta_2)(\beta_3+\beta_7)-\frac{b}{2}\beta_3(\beta_3+\beta_7)&\beta_2-\frac{1}{2}(\alpha_1+\beta_2)(\beta_3+\beta_8)-\frac{b}{2}\beta_3(\beta_3+\beta_8)&\beta_3\\
*&*&0\end{array}\right), \]
\[\mathcal{A}'_2=g(\mathcal{A}_2+b\mathcal{A}_3)g^{-1}= 
\left( 
\begin{array}{ccc}
\alpha_4+b\alpha_7&\alpha_5+b\alpha_8&0\\
\beta_4+b\beta_7&\beta_5+b\beta_8&0\\
*&*&*\end{array}\right),\] 
\[\mathcal{A}'_3=g\mathcal{A}_3g^{-1}=\left( 
\begin{array}{ccc}
\alpha_7&\alpha_8&0\\
\beta_7&1-\alpha_7&0\\
*&*&-1\end{array}\right).\] 

Therefore if $\alpha_8\neq 0$ one can make $\alpha'_5= 0$ to get
\[A_{23}=\left(
\begin{array}{ccccccccc}
\alpha_1 & \alpha_2& 1&\alpha_4&0&0&\alpha_7&\alpha_8&0\\
\beta_1 & -\alpha_1& \beta_3&\beta_4&\beta_5&0&\beta_7&1-\alpha_7&0\\
\gamma_1 & \gamma_2& 0&\gamma_4&\gamma_5&1-\alpha_4-\beta_5&1-\alpha_1-\beta_4&-\alpha_2-\beta_5&-1\end{array}\right),\] where $\alpha_8\neq 0$, if $\alpha_8= 0$ and $\alpha_7\neq 0$ one can make $\alpha'_4= 0$ to get
\[A_{24}=\left(
\begin{array}{ccccccccc}
\alpha_1 & \alpha_2& 1&0&\alpha_5&0&\alpha_7&0&0\\
\beta_1 & -\alpha_1& \beta_3&\beta_4&\beta_5&0&\beta_7&1-\alpha_7&0\\
\gamma_1 & \gamma_2& 0&\gamma_4&\gamma_5&1-\beta_5&1-\alpha_1-\beta_4&-\alpha_2-\beta_5&-1\end{array}\right),\] where $\alpha_7\neq 0$, if $\alpha_8= 0$, $\alpha_7= 0$ then $\beta_8=1$ and one can make $\beta'_5=0$ to get
\[A_{25}=\left(
\begin{array}{ccccccccc}
\alpha_1 & \alpha_2& 1&\alpha_4&\alpha_5&0&0&0&0\\
\beta_1 & -\alpha_1& \beta_3&\beta_4&0&0&\beta_7&1&0\\
\gamma_1 & \gamma_2& 0&\gamma_4&\gamma_5&1-\alpha_4&1-\alpha_1-\beta_4&-\alpha_2&-1\end{array}\right).\]  ( note that $\alpha_7+\beta_8=\alpha_3=1$ and in $\alpha_7=0$ case $\beta_8=1$, and condition $\beta_8=0$ has no meaning).

If $\alpha_9=\beta_9=\alpha_6= \alpha_7+\beta_8-\alpha_3=\alpha_3= 0, \beta_3\neq 0$ then one can make $\beta'_3=1$ by appropriate choose of $c$ to come to:

\underline{\textbf{Case 6.\ $\alpha_9=\beta_9=\alpha_6= \alpha_7+\beta_8-\alpha_3=\alpha_3= 0, \beta_3= 1$}} with respect to $g$ with $c=1$. If $b=-(\alpha_1+\beta_2)$ then $\gamma'_3=0$ and
\[\mathcal{A}'_1=g(\mathcal{A}_1+b\mathcal{A}_3)g^{-1}= 
\left(\begin{array}{ccc}
\alpha_1+a\alpha_7&\alpha_2+a\alpha_8&0\\
\beta_1+a(1+\beta_7)&\beta_2+a(1+\beta_8)&1\\
*&*&0\end{array}\right),\]
\[\mathcal{A}'_2=g(\mathcal{A}_2+b\mathcal{A}_3)g^{-1}= 
\left( 
\begin{array}{ccc}
\alpha_4+b\alpha_7&\alpha_5+b\alpha_8&0\\
\beta_4+b\beta_7&\beta_5+b\beta_8&0\\
*&*&*\end{array}\right),\]
\[\mathcal{A}'_3=g\mathcal{A}_3g^{-1}=\left( 
\begin{array}{ccc}
-\beta_8&\alpha_8&0\\
\beta_7&\beta_8&0\\
*&*&0\end{array}\right).\] 

 Therefore if $\alpha_8\neq 0$ one can make $\alpha'_2= 0$ to get
 \[A_{26}=\left(
\begin{array}{ccccccccc}
\alpha_1 & 0& 0&\alpha_4&\alpha_5&0&-\beta_8&\alpha_8&0\\
\beta_1 & -\alpha_1& 1&\beta_4&\beta_5&0&\beta_7&\beta_8&0\\
\gamma_1 & \gamma_2& 0&\gamma_4&\gamma_5&1-\alpha_4-\beta_5&1-\alpha_1-\beta_4&-\beta_5&0\end{array}\right),\] where $\alpha_8\neq 0$, if $\alpha_8= 0$, $\alpha_7\neq 0$ one can make $\alpha'_1= 0$ to get
\[A_{27}=\left(
\begin{array}{ccccccccc}
0 & \alpha_2& 0&\alpha_4&\alpha_5&0&-\beta_8&0&0\\
\beta_1 & 0& 1&\beta_4&\beta_5&0&\beta_7&\beta_8&0\\
\gamma_1 & \gamma_2& 0&\gamma_4&\gamma_5&1-\alpha_4-\beta_5&1-\beta_4&-\alpha_2-\beta_5&0\end{array}\right),\] where $\alpha_7=-\beta_8\neq 0$ and 
if $\alpha_8=\alpha_7=-\beta_8= 0$ then on can make $\beta'_2=0$ to get
\[A_{28}=\left(
\begin{array}{ccccccccc}
\alpha_1 & \alpha_2& 0&\alpha_4&\alpha_5&0&0&0&0\\
\beta_1 & 0& 1&\beta_4&\beta_5&0&\beta_7&0&0\\
\gamma_1 & \gamma_2& -\alpha_1&\gamma_4&\gamma_5&1-\alpha_4-\beta_5&1-\alpha_1-\beta_4&-\alpha_2-\beta_5&0\end{array}\right).\]

\underline{\textbf{Case 7.\ $\alpha_9=\beta_9=\alpha_6= \alpha_7+\beta_8-\alpha_3=\alpha_3= \beta_3= 0$}.} In this case
$\alpha'_9=\beta'_9=\alpha_6'= \beta'_6=\alpha'_3=\beta'_3=0$ and \[A'_1=g(A_1+aA_3)g^{-1}= 
 \left(\begin{array}{ccc}
\alpha_1+a\alpha_7&\alpha_2+a\alpha_8&0\\
\beta_1+a\beta_7&\beta_2+a\beta_8)&0\\
*&*&*\end{array}\right), \]
\[A'_2=g(A_2+aA_3)g^{-1}= 
 \left(\begin{array}{ccc}
\alpha_4+b\alpha_7&\alpha_5+b\alpha_8&0\\
\beta_4+b\beta_7&\beta_5+b\beta_8)&0\\
*&*&*\end{array}\right), \] 
\[A'_3=cgA_3g^{-1}= 
\left(\begin{array}{ccc}
c\alpha_7&c\alpha_8&0\\
c\beta_7&c\beta_8&0\\
*&*&0\end{array}\right). \]
Therefore if  $\alpha_8\neq 0$ one can make $\alpha'_2=\alpha'_5=0, \alpha'_8=1$ to get
\[A_{29}=\left(
\begin{array}{ccccccccc}
\alpha_1 & 0& 0&\alpha_4&0&0&\alpha_7&1&0\\
\beta_1 & \beta_2& 0&\beta_4&\beta_5&0&\beta_7&-\alpha_7&0\\
\gamma_1 & \gamma_2& -\alpha_1-\beta_2&\gamma_4&\gamma_5&1-\alpha_4-\beta_5&1-\alpha_1-\beta_4&-\beta_5&0\end{array}\right),\] if $\alpha_8= 0$, $\alpha_7\neq 0$ one can make $\alpha'_1=\alpha'_4=0, \alpha'_7=1$ to get

\[A_{30}=\left(
\begin{array}{ccccccccc}
0 & \alpha_2& 0&0&\alpha_5&0&1&0&0\\
\beta_1 & \beta_2& 0&\beta_4&\beta_5&0&\beta_7&-1&0\\
\gamma_1 & \gamma_2& -\beta_2&\gamma_4&\gamma_5&1-\beta_5&1-\beta_4&-\alpha_2-\beta_5&0\end{array}\right),\] if  $\alpha_7=\alpha_8= 0$ then $\beta_8=0$ and if $\beta_7\neq 0$ one can make $\beta'_1=\beta'_4=0, \beta'_7=1$ to get
\[A_{31}=\left(
\begin{array}{ccccccccc}
\alpha_1 & \alpha_2& 0&\alpha_4&\alpha_5&0&0&0&0\\
0 & \beta_2& 0&0&\beta_5&0&1&0&0\\
\gamma_1 & \gamma_2& -\alpha_1-\beta_2&\gamma_4&\gamma_5&1-\alpha_4-\beta_5&1-\alpha_1&-\alpha_2-\beta_5&0\end{array}\right).\] 

%\underline{\textbf{Case 8.\
In $\alpha_9=\alpha_8=\alpha_7=\alpha_6=\alpha_3=\beta_9=\beta_8=\beta_7=\beta_6=\beta_3=0$.case
\[A'_1=g(A_1+aA_3)g^{-1}= 
\left(\begin{array}{ccc}
\alpha_1&\alpha_2&0\\
\beta_1&\beta_2&0\\
\gamma'_1&\gamma'_2&-\alpha_1-\beta_2\end{array}\right), \]
\[A'_2=g(A_2+aA_3)g^{-1}= 
\left(\begin{array}{ccc}
\alpha_4&\alpha_5&0\\
\beta_4&\beta_5&0\\
\gamma'_4&\gamma'_5&1-\alpha_4-\beta_5\end{array}\right), \] 
\[A'_3=cgA_3g^{-1}= 
\left(\begin{array}{ccc}
0&0&0\\
0&0&0\\
1-\alpha_1-\beta_4&-\alpha_2-\beta_5&0\end{array}\right),\] where \\
$\gamma'_1=c^{-1}[-a(3\alpha_1+\beta_2+\beta_4-1)-b\beta_1+\gamma_1]$,\\
$\gamma'_2=c^{-1}[-a(2\alpha_2+\alpha_1+\beta_5+\beta_2)-b\beta_2+\gamma_2]$,\\
$\gamma'_4=c^{-1}[-a(2\alpha_4+\beta_5-1)-b(\alpha_1+2\beta_4-1)+\gamma_4]$,\\
$\gamma'_5=c^{-1}[-a\alpha_5-b(\alpha_2+\alpha_4+3\beta_5-1)+\gamma_5]$.

 Consider row vectors \[ \begin{array}{ll} 1.\  (3\alpha_1+\beta_2+\beta_4-1,\beta_1), & 1'.\  (3\alpha_1+\beta_2+\beta_4-1,\beta_1, \gamma_1),\\
   2.\ (2\alpha_2+\alpha_1+\beta_5+\beta_2, \beta_2),& 2'.\  (2\alpha_2+\alpha_1+\beta_5+\beta_2, \beta_2, \gamma_2),\\
       3.\ (2\alpha_4+\beta_5-1, \alpha_1+2\beta_4-1),&  3'.\
(2\alpha_4+\beta_5-1, \alpha_1+2\beta_4-1, \gamma_4),\\
4.\ (\alpha_5, \alpha_2+\alpha_4+3\beta_5-1),& 4'.\
(\alpha_5, \alpha_2+\alpha_4+3\beta_5-1, \gamma_5),\end{array}\]
 and let
$M_{(i,j,k,l)}(A)=M_{(i,j,k,l)}$ (respectively, $M'_{(i,j,k,l)}(A)=M'_{(i,j,k,l)}$), where $i,j,k,l=1,2,3,4$, stand for matrix with the above defined rows $i.,j.,k.,l.$ (respectively, $i'.,j'.,k'.,l'.$), $\lambda\times i+ \mu \times j$, where $\lambda, \mu \in F$, stand for a linear combination of the rows $i.,j$, $rk(A)$ stand for the rank of a matrix $A$.

\underline{\textbf{Case 8. $rk(M_{(1,2,3,4)})=rk(M'_{(1,2,3,4)})$}.} In this case one can make $\gamma'_1=\gamma'_2=\gamma'_4=\gamma'_5=0$ to get
\[A_{32}=\left(
\begin{array}{ccccccccc}
\alpha_1 & \alpha_2& 0&\alpha_4&\alpha_5&0&0&0&0\\
\beta_1 & \beta_2& 0&\beta_4&\beta_5&0&0&0&0\\
0 & 0& -\alpha_1-\beta_2&0&0&1-\alpha_4-\beta_5&1-\alpha_1-\beta_4&-\alpha_2-\beta_5&0\end{array}\right).\] 

\underline{\textbf{Case 9. $rk(M_{(1,2,3,4)})\neq rk(M'_{(1,2,3,4)})$}.}

1. If  $rk(M_{(2,3,4)})=rk(M'_{(2,3,4)})$ one can make $\gamma'_2=\gamma'_4=\gamma'_5=0$ and $\gamma'_1=1$ to get
\[A_{33}=\left(
\begin{array}{ccccccccc}
\alpha_1 & \alpha_2& 0&\alpha_4&\alpha_5&0&0&0&0\\
\beta_1 & \beta_2& 0&\beta_4&\beta_5&0&0&0&0\\
1 & 0& -\alpha_1-\beta_2&0&0&1-\alpha_4-\beta_5&1-\alpha_1-\beta_4&-\alpha_2-\beta_5&0\end{array}\right),\] where $rk(M_{(1,2,3,4)}(A_{33}))=rk(M_{(2,3,4)}(A_{33}))$.

2. If $rk(M_{(2,3,4)})<rk(M'_{(2,3,4)})$, $rk(M_{(1,3,4)})=rk(M'_{(1,3,4)})$ one can make $\gamma'_1=\gamma'_4=\gamma'_5=0$ and $\gamma'_2=1$ to get
\[A_{34}=\left(
\begin{array}{ccccccccc}
\alpha_1 & \alpha_2& 0&\alpha_4&\alpha_5&0&0&0&0\\
\beta_1 & \beta_2& 0&\beta_4&\beta_5&0&0&0&0\\
0 & 1& -\alpha_1-\beta_2&0&0&1-\alpha_4-\beta_5&1-\alpha_1-\beta_4&-\alpha_2-\beta_5&0\end{array}\right),\] where $rk(M_{(1,2,3,4)}(A_{34}))=rk(M_{(1,3,4)}(A_{34}))$, $rk(M_{(3,4)}(A_{34}))=rk(M_{(2,3,4)}(A_{34}))$.

3. If $rk(M_{(2,3,4)})<rk(M'_{(2,3,4)})$, $rk(M_{(1,3,4)})<rk(M'_{(1,3,4)})$, $rk(M_{(1,2,4)})=rk(M'_{(1,2,4)})$ one can make $\gamma'_1=\gamma'_2=\gamma'_5=0$ and $\gamma'_4=1$ to get
\[A_{35}=\left(
\begin{array}{ccccccccc}
\alpha_1 & \alpha_2& 0&\alpha_4&\alpha_5&0&0&0&0\\
\beta_1 & \beta_2& 0&\beta_4&\beta_5&0&0&0&0\\
0 & 0& -\alpha_1-\beta_2&1&0&1-\alpha_4-\beta_5&1-\alpha_1-\beta_4&-\alpha_2-\beta_5&0\end{array}\right),\] where $rk(M_{(1,2,3,4)}(A_{35}))=rk(M_{(1,2,4)}(A_{35}))$, $rk(M_{(2,4)}(A_{35}))=rk(M_{(2,3,4)}(A_{35}))$, $rk(M_{(1,4)}(A_{35}))=rk(M_{(1,3,4)}(A_{35}))$.

4. If $rk(M_{(2,3,4)})<rk(M'_{(2,3,4)})$, $rk(M_{(1,3,4)})<rk(M'_{(1,3,4)})$, $rk(M_{(1,2,4)})<rk(M'_{(1,2,4)})$, $rk(M_{(1,2,3)})=rk(M'_{(1,2,3)})$ . In this case one can make $\gamma'_1=\gamma'_2=\gamma'_4=0$ and $\gamma'_5=1$ to get
\[A_{36}=\left(
\begin{array}{ccccccccc}
\alpha_1 & \alpha_2& 0&\alpha_4&\alpha_5&0&0&0&0\\
\beta_1 & \beta_2& 0&\beta_4&\beta_5&0&0&0&0\\
0 & 0& -\alpha_1-\beta_2&0&1&1-\alpha_4-\beta_5&1-\alpha_1-\beta_4&-\alpha_2-\beta_5&0\end{array}\right),\] where $rk(M_{(1,2,3,4)}(A_{36}))=rk(M_{(1,2,3)}(A_{36}))$, $rk(M_{(2,3)}(A_{36}))=rk(M_{(2,3,4)}(A_{36}))$, $rk(M_{(1,3)}(A_{36}))=rk(M_{(1,3,4)}(A_{36}))$, $rk(M_{(1,2)}(A_{36}))=rk(M_{(1,2,4)}(A_{36}))$.

\underline{\textbf{Subcase 9-1. $rk(M_{(2,3,4)})<rk(M'_{(2,3,4)})$, $rk(M_{(1,3,4)})<rk(M'_{(1,3,4)})$},}\\ \underline{$rk(M_{(1,2,4)})<rk(M'_{(1,2,4)})$, $rk(M_{(1,2,3)})<rk(M'_{(1,2,3)})$.}

1. If $rk(M_{(3,4)})=rk(M'_{(3,4)})$ one can make $\gamma'_4=\gamma'_5=0$, $\gamma'_1=1$ and $\gamma'_2\neq 0$ to get
\[A_{37}=\left(
\begin{array}{ccccccccc}
\alpha_1 & \alpha_2& 0&\alpha_4&\alpha_5&0&0&0&0\\
\beta_1 & \beta_2& 0&\beta_4&\beta_5&0&0&0&0\\
1 & \gamma_2& -\alpha_1-\beta_2&0&0&1-\alpha_4-\beta_5&1-\alpha_1-\beta_4&-\alpha_2-\beta_5&0\end{array}\right),\] where $rk(M_{(1,2,3,4)}(A_{37}))=rk(M_{(2-\gamma_2\times1,3,4)}(A_{37}))$, $rk(M_{(3,4)}(A_{37}))=rk(M_{(2,3,4)}(A_{37}))$, $rk(M_{(3,4)}(A_{37}))=rk(M_{(1,3,4)}(A_{37}))$, $rk(M_{(2-\gamma_2\times 1,4)}(A_{37}))=rk(M_{(1,2,4)}(A_{37}))$, 
$rk(M_{(2-\gamma_2\times 1,3)}(A_{37}))=rk(M_{(1,2,3)}(A_{37}))$ and $\gamma_2\neq 0$.

2. If $rk(M_{(3,4)})<rk(M'_{(3,4)})$, $rk(M_{(2,4)})=rk(M'_{(2,4)})$ one can make $\gamma'_2=\gamma'_5=0$, $\gamma'_1=1$ and $\gamma'_4\neq 0$ to get
\[A_{38}=\left(
\begin{array}{ccccccccc}
\alpha_1 & \alpha_2& 0&\alpha_4&\alpha_5&0&0&0&0\\
\beta_1 & \beta_2& 0&\beta_4&\beta_5&0&0&0&0\\
1 & 0& -\alpha_1-\beta_2&\gamma_4&0&1-\alpha_4-\beta_5&1-\alpha_1-\beta_4&-\alpha_2-\beta_5&0\end{array}\right),\] where $rk(M_{(1,2,3,4)}(A_{38}))=rk(M_{(2,3-\gamma_4\times 1,4)}(A_{38}))$, $rk(M_{(2,4)}(A_{38}))=rk(M_{(2,3,4)}(A_{38}))$, $rk(M_{(3-\gamma_4\times 1,4)}(A_{38}))=rk(M_{(1,3,4)}(A_{38}))$, $rk(M_{(2,4)}(A_{38}))=rk(M_{(1,2,4)}(A_{38}))$, 
$rk(M_{(2,3-\gamma_4\times 1)}(A_{38}))=rk(M_{(1,2,3)}(A_{38}))$, $rk(M_{(4)}(A_{38}))=rk(M_{(3,4)}(A_{38}))$ and $\gamma_4\neq 0$.

3. If $rk(M_{(3,4)})<rk(M'_{(3,4)})$, $rk(M_{(2,4)})<rk(M'_{(2,4)})$, $rk(M_{(1,4)})=rk(M'_{(1,4)})$ one can make $\gamma'_1=\gamma'_5=0$, $\gamma'_2=1$ and $\gamma'_4\neq 0$ to get
\[A_{39}=\left(
\begin{array}{ccccccccc}
\alpha_1 & \alpha_2& 0&\alpha_4&\alpha_5&0&0&0&0\\
\beta_1 & \beta_2& 0&\beta_4&\beta_5&0&0&0&0\\
0 & 1& -\alpha_1-\beta_2&\gamma_4&0&1-\alpha_4-\beta_5&1-\alpha_1-\beta_4&-\alpha_2-\beta_5&0\end{array}\right),\] where $rk(M_{(1,2,3,4)}(A_{39}))=rk(M_{(1,3-\gamma_4\times 2,4)}(A_{39}))$, $rk(M_{(3-\gamma_4\times 2,4)}(A_{39}))=rk(M_{(2,3,4)}(A_{39}))$, $rk(M_{(1,4)}(A_{39}))=rk(M_{(1,3,4)}(A_{39}))$, $rk(M_{(1,4)}(A_{39}))=rk(M_{(1,2,4)}(A_{39}))$, 
$rk(M_{(1,3-\gamma_4\times 2)}(A_{39}))=rk(M_{(1,2,3)}(A_{39}))$, $rk(M_{(4)}(A_{39}))=rk(M_{(3,4)}(A_{39}))$,  $rk(M_{(4)}(A_{39}))=rk(M_{(2,4)}(A_{39}))$ and $\gamma_4\neq 0$.

4. If $rk(M_{(3,4)})<rk(M'_{(3,4)})$, $rk(M_{(2,4)})<rk(M'_{(2,4)})$, $rk(M_{(1,4)})<rk(M'_{(1,4)})$ , $rk(M_{(2,3)})=rk(M'_{(2,3)})$ one can make $\gamma'_2=\gamma'_4=0$, $\gamma'_1=1$ and $\gamma'_5\neq 0$ to get
\[A_{40}=\left(
\begin{array}{ccccccccc}
\alpha_1 & \alpha_2& 0&\alpha_4&\alpha_5&0&0&0&0\\
\beta_1 & \beta_2& 0&\beta_4&\beta_5&0&0&0&0\\
1 & 0& -\alpha_1-\beta_2&0&\gamma_5&1-\alpha_4-\beta_5&1-\alpha_1-\beta_4&-\alpha_2-\beta_5&0\end{array}\right),\] where $rk(M_{(1,2,3,4)}(A_{40}))=rk(M_{(2,3,4-\gamma_5\times 1)}(A_{40}))$, $rk(M_{(2,3)}(A_{40}))=rk(M_{(2,3,4)}(A_{40}))$, $rk(M_{(3,4-\gamma_5\times 1)}(A_{40}))=rk(M_{(1,3,4)}(A_{40}))$, $rk(M_{(2,4-\gamma_5\times 1)}(A_{40}))=rk(M_{(1,2,4)}(A_{40}))$, 
$rk(M_{(2,3)}(A_{40}))=rk(M_{(1,2,3)}(A_{40}))$, $rk(M_{(3)}(A_{40}))=rk(M_{(3,4)}(A_{40}))$, $rk(M_{(2)}(A_{40}))=rk(M_{(2,4)}(A_{40}))$, $rk(M_{(4-\gamma_5\times 1)}(A_{40}))=rk(M_{(1,4)}(A_{40}))$ and $\gamma_5\neq 0$.

5. If $rk(M_{(3,4)})<rk(M'_{(3,4)})$, $rk(M_{(2,4)})<rk(M'_{(2,4)})$, $rk(M_{(1,4)})<rk(M'_{(1,4)})$ , $rk(M_{(2,3)})<rk(M'_{(2,3)})$, $rk(M_{(1,3)})=rk(M'_{(1,3)})$ one can make $\gamma'_1=\gamma'_4=0$, $\gamma'_2=1$ and $\gamma'_5\neq 0$ to get
\[A_{41}=\left(
\begin{array}{ccccccccc}
\alpha_1 & \alpha_2& 0&\alpha_4&\alpha_5&0&0&0&0\\
\beta_1 & \beta_2& 0&\beta_4&\beta_5&0&0&0&0\\
0 & 1& -\alpha_1-\beta_2&0&\gamma_5&1-\alpha_4-\beta_5&1-\alpha_1-\beta_4&-\alpha_2-\beta_5&0\end{array}\right),\] where $rk(M_{(1,2,3,4)}(A_{41}))=rk(M_{(1,3,4-\gamma_5\times 2)}(A_{41}))$, $rk(M_{(3,4-\gamma_5\times 2)}(A_{41}))=rk(M_{(2,3,4)}(A_{41}))$, $rk(M_{(1,3)}(A_{41}))=rk(M_{(1,3,4}(A_{41}))$, $rk(M_{(1,4-\gamma_5\times 2)}(A_{41}))=rk(M_{(1,2,4)}(A_{41}))$, 
$rk(M_{(1,3)}(A_{41}))=rk(M_{(1,2,3)}(A_{41}))$, $rk(M_{(3)}(A_{41}))=rk(M_{(3,4)}(A_{41}))$, $rk(M_{(4-\gamma_5\times 2)}(A_{41}))=rk(M_{(2,4)}(A_{41}))$, $rk(M_{(1)}(A_{41}))=rk(M_{(1,4)}(A_{41}))$, $rk(M_{(3)}(A_{41}))=rk(M_{(2,3)}(A_{41}))$ and $\gamma_5\neq 0$.

6. If $rk(M_{(3,4)})<rk(M'_{(3,4)})$, $rk(M_{(2,4)})<rk(M'_{(2,4)})$, $rk(M_{(1,4)})<rk(M'_{(1,4)})$ , $rk(M_{(2,3)})<rk(M'_{(2,3)})$, $rk(M_{(1,3)})<rk(M'_{(1,3)})$, $rk(M_{(1,2)})=rk(M'_{(1,2)})$
 one can make $\gamma'_1=\gamma'_2=0$, $\gamma'_4=1$ and $\gamma'_5\neq 0$ to get
\[A_{42}=\left(
\begin{array}{ccccccccc}
\alpha_1 & \alpha_2& 0&\alpha_4&\alpha_5&0&0&0&0\\
\beta_1 & \beta_2& 0&\beta_4&\beta_5&0&0&0&0\\
0 & 0& -\alpha_1-\beta_2&1&\gamma_5&1-\alpha_4-\beta_5&1-\alpha_1-\beta_4&-\alpha_2-\beta_5&0\end{array}\right),\] where $rk(M_{(1,2,3,4)}(A_{42}))=rk(M_{(1,2,4-\gamma_5\times 3)}(A_{42}))$, $rk(M_{(2,4-\gamma_5\times 3)}(A_{42}))=rk(M_{(2,3,4)}(A_{42}))$, $rk(M_{(1,4-\gamma_5\times 3)}(A_{42}))=rk(M_{(1,3,4}(A_{42}))$, $rk(M_{(1,2)}(A_{42}))=rk(M_{(1,2,4)}(A_{42}))$, 
$rk(M_{(1,2)}(A_{42}))=rk(M_{(1,2,3)}(A_{42}))$, $rk(M_{(4-\gamma_5\times 3)}(A_{42}))=rk(M_{(3,4)}(A_{42}))$, $rk(M_{(2)}(A_{42}))=rk(M_{(2,4)}(A_{42}))$, $rk(M_{(1)}(A_{42}))=rk(M_{(1,4)}(A_{42}))$, $rk(M_{(1)}(A_{42}))=rk(M_{(1,3)}(A_{42}))$ and $\gamma_5\neq 0$.

\underline{\textbf{Subsubcase 9-1-1. $rk(M_{(2,3,4)})<rk(M'_{(2,3,4)})$, $rk(M_{(1,3,4)})<rk(M'_{(1,3,4)})$},}\\ \underline{$rk(M_{(1,2,4)})<rk(M'_{(1,2,4)})$, $rk(M_{(1,2,3)})<rk(M'_{(1,2,3)})$,  $rk(M_{(3,4)})<rk(M'_{(3,4)})$,}\\ \underline{$rk(M_{(2,4)})<rk(M'_{(2,4)})$, $rk(M_{(1,4)})<rk(M'_{(1,4)})$ , $rk(M_{(2,3)})<rk(M'_{(2,3)})$,}\\ \underline{$rk(M_{(1,3)})<rk(M'_{(1,3)})$, $rk(M_{(1,2)})<rk(M'_{(1,2)})$.}

1. If $rk(M_{(1)})=rk(M'_{(1)})$ one can make $\gamma'_1=0$, $\gamma'_2=1$, $\gamma'_4\neq 0$, $\gamma'_5\neq 0$ to get
 \[A_{43}=\left(
 \begin{array}{ccccccccc}
 \alpha_1 & \alpha_2& 0&\alpha_4&\alpha_5&0&0&0&0\\
 \beta_1 & \beta_2& 0&\beta_4&\beta_5&0&0&0&0\\
 0 & 1& -\alpha_1-\beta_2&\gamma_4&\gamma_5&1-\alpha_4-\beta_5&1-\alpha_1-\beta_4&-\alpha_2-\beta_5&0\end{array}\right),\] where $rk(M_{(1,2,3,4)})=rk(M_{(1,3-\gamma_4\times 2,4-\gamma_5\times 2)})$, $rk(M_{(3-\gamma_4\times 2, 4-\gamma_5\times 2)}(A_{43}))=rk(M_{(2,3,4)}(A_{43}))$, $rk(M_{(1,\gamma_4\times 4-\gamma_5\times 3)}(A_{43}))=rk(M_{(1,3,4}(A_{43}))$, $rk(M_{(1,4-\gamma_5\times 2)}(A_{43}))=rk(M_{(1,2,4)}(A_{43}))$, 
 $rk(M_{(1,3-\gamma_4\times 2)}(A_{43}))=rk(M_{(1,2,3)}(A_{43}))$, $rk(M_{(\gamma_4\times 4-\gamma_5\times 3)}(A_{43}))=rk(M_{(3,4)}(A_{43}))$, $rk(M_{(4-\gamma_5\times 2)}(A_{43}))=rk(M_{(2,4)}(A_{43}))$, $rk(M_{(1)}(A_{43}))=rk(M_{(1,4)}(A_{43}))$, $rk(M_{(1)}(A_{43}))=rk(M_{(1,3)}(A_{43}))$, $rk(M_{(1)}(A_{43}))=rk(M_{(1,2)}(A_{43}))$ and $\gamma_4\neq 0$, $\gamma_5\neq 0$.

2. If $rk(M_{(1)})<rk(M'_{(1)})$, $rk(M_{(2)})=rk(M'_{(2)})$ one can make $\gamma'_1=1$, $\gamma'_2=0$, $\gamma'_4\neq 0$, $\gamma'_5\neq 0$ to get
 \[A_{44}=\left(
 \begin{array}{ccccccccc}
 \alpha_1 & \alpha_2& 0&\alpha_4&\alpha_5&0&0&0&0\\
 \beta_1 & \beta_2& 0&\beta_4&\beta_5&0&0&0&0\\
 1 & 0& -\alpha_1-\beta_2&\gamma_4&\gamma_5&1-\alpha_4-\beta_5&1-\alpha_1-\beta_4&-\alpha_2-\beta_5&0\end{array}\right),\] where $rk(M_{(1,2,3,4)}(A_{44}))=rk(M_{(2,3-\gamma_4\times 1,4-\gamma_5\times 1)}(A_{44}))$, $rk(M_{(2,\gamma_4\times 3-\gamma_5\times 3)}(A_{44}))=rk(M_{(2,3,4)}(A_{44}))$, $rk(M_{(3-\gamma_4\times 1,4-\gamma_5\times 1)}(A_{44}))=rk(M_{(1,3,4}(A_{44}))$, $rk(M_{(2,4-\gamma_5\times 1)}(A_{44}))=rk(M_{(1,2,4)}(A_{44}))$, 
 $rk(M_{(2,3-\gamma_4\times 1)}(A_{44}))=rk(M_{(1,2,3)}(A_{44}))$, $rk(M_{(\gamma_4\times 4-\gamma_5\times 3)}(A_{44}))=rk(M_{(3,4)}(A_{44}))$, $rk(M_{(2)}(A_{44}))=rk(M_{(2,4)}(A_{44}))$, $rk(M_{(4-\gamma_4\times 1)}(A_{44}))=rk(M_{(1,4)}(A_{44}))$, $rk(M_{(3-\gamma_5\times 1)}(A_{44}))=rk(M_{(1,3)}(A_{44}))$, $rk(M_{(2)}(A_{44}))=rk(M_{(1,2)}(A_{44}))$, $rk(M_{(1)}(A_{44}))=0$ and $\gamma_4\neq 0$, $\gamma_5\neq 0$.
 
  3. If $rk(M_{(1)})<rk(M'_{(1)})$, $rk(M_{(2)})<rk(M'_{(2)})$, $rk(M_{(3)})=rk(M'_{(3)})$ one can make $\gamma'_1=1$, $\gamma'_2\neq 0$, $\gamma'_4= 0$, $\gamma'_5\neq 0$ to get MSC
  \[\left(
  \begin{array}{ccccccccc}
  \alpha_1 & \alpha_2& 0&\alpha_4&\alpha_5&0&0&0&0\\
  \beta_1 & \beta_2& 0&\beta_4&\beta_5&0&0&0&0\\
  1 & \gamma_2& -\alpha_1-\beta_2&0&\gamma_5&1-\alpha_4-\beta_5&1-\alpha_1-\beta_4&-\alpha_2-\beta_5&0\end{array}\right)\] for which $rk(M_{(1,2,3,4)})=rk(M_{(2-\gamma_2\times 1,3,4-\gamma_5\times 1)})$, $rk(M_{(3,\gamma_2\times 4-\gamma_5\times 2)})=rk(M_{(2,3,4)})$, $rk(M_{(3,4-\gamma_5\times 1)})=rk(M_{(1,3,4})$, $rk(M_{(2-\gamma_2\times 1,4-\gamma_5\times 1)})=rk(M_{(1,2,4)})$, 
  $rk(M_{(2-\gamma_2\times 1),3})=rk(M_{(1,2,3)})$, $rk(M_{(3)})=rk(M_{(3,4)})$, $rk(M_{(\gamma_2\times 4-\gamma_5\times 2)})=rk(M_{(2,4)})$, $rk(M_{(4-\gamma_5\times 1)})=rk(M_{(1,4)})$, $rk(M_{(3)})=rk(M_{(1,3)})$, $rk(M_{(1)})=0$, $rk(M_{(2)})=0$ and $\gamma_2\neq 0$, $\gamma_5\neq 0$.
  
  The system $rk(M_{(1)})=0$, $rk(M_{(2)})=0$ implies that $\beta_1=\beta_2=0, \beta_4=1-3\alpha_1, \beta_5=-\alpha_1-2\alpha_2$ and one can make $\gamma'_1= 1$ to get \[A_{45}=\left(
  \begin{array}{ccccccccc}
  \alpha_1 & \alpha_2& 0&\alpha_4&\alpha_5&0&0&0&0\\
  0 & 0& 0&1-3\alpha_1&-\alpha_1-2\alpha_2&0&0&0&0\\
  1 & \gamma_2& -\alpha_1&0&\gamma_5&1+\alpha_1+2\alpha_2-\alpha_4&2\alpha_1&\alpha_1+\alpha_2&0\end{array}\right),\] where $\gamma_2\neq 0$, $\gamma_5\neq 0$.

4. If $rk(M_{(1)})<rk(M'_{(1)})$, $rk(M_{(2)})<rk(M'_{(2)})$,
 $rk(M_{(3)})<rk(M'_{(3)})$, $rk(M_{(4)})=rk(M'_{(4)})$
 one can make $\gamma'_1=1$, $\gamma'_2\neq 0$, $\gamma'_4\neq 0$, $\gamma'_5=0$ to get
MSC \[\left(
\begin{array}{ccccccccc}
\alpha_1 & \alpha_2& 0&\alpha_4&\alpha_5&0&0&0&0\\
\beta_1 & \beta_2& 0&\beta_4&\beta_5&0&0&0&0\\
1 & \gamma_2& -\alpha_1-\beta_2&\gamma_4&0&1-\alpha_4-\beta_5&1-\alpha_1-\beta_4&-\alpha_2-\beta_5&0\end{array}\right)\] for which $rk(M_{(1,2,3,4)})=rk(M_{(2-\gamma_2\times 1,3-\gamma_4\times 1, 4)})$, $rk(M_{(\gamma_2\times 3-\gamma_4\times 2,4)})=rk(M_{(2,3,4)})$, $rk(M_{(3-\gamma_4\times 1,4)})=rk(M_{(1,3,4})$, $rk(M_{(2-\gamma_2\times 1,4)})=rk(M_{(1,2,4)})$, 
$rk(M_{(2-\gamma_2\times 1),3-\gamma_4\times 1})=rk(M_{(1,2,3)})$, $rk(M_{(4)})=rk(M_{(3,4)})$, $rk(M_{(4)})=rk(M_{(2,4)})$, $rk(M_{(4)})=rk(M_{(1,4)})$, $rk(M_{(3-\gamma_4\times 1)})=rk(M_{(1,3)})$, $rk(M_{(1)})=0$, $rk(M_{(2)})=0$, $rk(M_{(3)})=0$ and $\gamma_2\neq 0$, $\gamma_4\neq 0$.
 
 But the system of equations $rk(M_{(1)})=0$, $rk(M_{(2)})=0$, $rk(M_{(3)})=0$ has no solution in $Ch.(F)=5$ case and otherwise $\alpha_1=\frac{1}{5}, \alpha_4=\frac{3}{5}+\alpha_2, \beta_1=0, \beta_2=0, \beta_4=\frac{2}{5}, \beta_5= -\frac{1}{5}-2\alpha_2$ and one can make $\gamma'_1= 1$ to have
\[A_{46}=\left(
\begin{array}{ccccccccc}
\frac{1}{5} & \alpha_2& 0&\frac{3}{5}+\alpha_2&\alpha_5&0&0&0&0\\
0 & 0& 0&\frac{2}{5}&-\frac{1}{5}-2\alpha_2&0&0&0&0\\
1 & \gamma_2& -\frac{1}{5}&\gamma_4&0&\frac{3}{5}+\alpha_2&\frac{2}{5}&\frac{1}{5}+\alpha_2&0\end{array}\right),\] where $\gamma_2\neq 0$, $\gamma_4\neq 0$.

5. The last one is $rk(M_{(1)})<rk(M'_{(1)})$, $rk(M_{(2)})<rk(M'_{(2)})$,
$rk(M_{(3)})<rk(M'_{(3)})$, $rk(M_{(4)})<rk(M'_{(4)})$, that is  $rk(M_{(1,2,3,4)})=0$,  $\gamma_1\neq 0$, $\gamma_2\neq 0$, $\gamma_4\neq 0$, $\gamma_5\neq 0$, case. 
The equation $rk(M_{(1,2,3,4)})=0$ has no solution in $Ch(\mathbb{F})=5$ case and otherwise ($Ch(F)\neq 2,5$) it has unique solution  $\alpha_1=\frac{1}{5},\ \alpha_2=-\frac{1}{4},\ \alpha_4=\frac{7}{20},\  \alpha_5=0,\  \beta_1=\beta_2=0,\ \beta_3= \frac{2}{5},\ \beta_5=\frac{3}{10}$ and one can make $\gamma'_1= 1$ to get

\[A_{47}=\left(
\begin{array}{ccccccccc}
\frac{1}{5} & -\frac{1}{4}& 0&\frac{7}{20}&0&0&0&0&0\\
0 & 0& 0&\frac{2}{5}&\frac{3}{10}&0&0&0&0\\
1 & \gamma_2& -\frac{1}{5} &\gamma_4&\gamma_5&\frac{7}{20} &\frac{2}{5} &-\frac{1}{20} &0\end{array}\right),\] where 
$\gamma_2\neq 0,\ \gamma_4\neq 0,\ \gamma_5\neq 0.$

\section{\bf Classification of $3$-dimensional algebras in $Ch.(\mathbb{F})=2$ case}
Now it is assumed that $Ch.(\mathbb{F})=2$ and we prove the following result.
\begin{theorem}. Let $Ch.(\mathbb{F})= 2$ be a field over which any  quadratic equation is solvable. In this case any $3$-dimensional algebra $\mathbb{A}$ over $\mathbb{F}$, for which $\{\mathbf{Tr_1}(A),\mathbf{Tr_2}(A)\}$ is linear independent, is isomorphic to only one algebra listed below by their (canonical) MSCs : 
	\[A_{1,2}=\left(
	\begin{array}{ccccccccc}
	\alpha_1 & \alpha_2& \alpha_3&\alpha_4&\alpha_5&\alpha_6&0&0&1\\
	\beta_1 & \beta_2& \beta_3&\beta_4&\beta_5&\beta_8+\alpha_3&\beta_7&\beta_9&\beta_9\\
	\gamma_1 & \gamma_2& \alpha_1+\beta_2&\gamma_4&\gamma_5&1+\alpha_4+\beta_5&1+\alpha_1+\beta_4&\alpha_2+\beta_5&-\beta_8\end{array}\right),\]
	\[A_{2,2}=\left(
	\begin{array}{ccccccccc}
	\alpha_1 & \alpha_2& \alpha_3&\alpha_4&\alpha_5&\alpha_6&\alpha_7&\alpha_8&0\\
	\beta_1 & \beta_2&\beta_3&\beta_4&\beta_5&\alpha_7+\alpha_3&0&0&1\\
	\gamma_1 & \gamma_2& \alpha_1+\beta_2&\gamma_4&\gamma_5&1+\alpha_4+\beta_5&1+\alpha_1+\beta_4&\alpha_2+\beta_5&-\alpha_7\end{array}\right),\]
	\[A_{3,2}=\left(
	\begin{array}{ccccccccc}
	\alpha_1 & \alpha_2& \alpha_3&0&0&1&\alpha_7&\alpha_8&0\\
	\beta_1 & \beta_2& \beta_3&\beta_4&\beta_5&\alpha_7+\beta_8+\alpha_3&\beta_7&\beta_8&0\\
	\gamma_1 & \gamma_2& \alpha_1+\beta_2&\gamma_4&\gamma_5&1+\beta_5&1+\alpha_1+\beta_4&\alpha_2+\beta_5&\alpha_7+\beta_8\end{array}\right),\] where $\alpha_8\neq 1$,
	\[A_{4,2}=\left(
	\begin{array}{ccccccccc}
	\alpha_1 & \alpha_2& \alpha_3&0&\alpha_5&1&\alpha_7&1&0\\
	\beta_1 & \beta_2& \beta_3&\beta_4&0&\alpha_7+\beta_8+\alpha_3&\beta_7&\beta_8&0\\
	\gamma_1 & \gamma_2& \alpha_1+\beta_2&\gamma_4&\gamma_5&1&1+\alpha_1+\beta_4&\alpha_2&\alpha_7+\beta_8\end{array}\right),\] where $\alpha_7+2\beta_8+\alpha_3\neq 0$,
	\[A_{5,2}=\left(
	\begin{array}{ccccccccc}
	\alpha_1 & \alpha_2& \alpha_3&0&\alpha_5&1&\alpha_3&1&0\\
	\beta_1 & \beta_2& \beta_3&0&\beta_5&\beta_8&\beta_7&\beta_8&0\\
	\gamma_1 & \gamma_2& \alpha_1+\beta_2&\gamma_4&\gamma_5&1+\beta_5&1+\alpha_1&\alpha_2+\beta_5&\alpha_3+\beta_8\end{array}\right),\]
	where $\beta_7\neq \alpha_3\beta_8$,
	\[A_{6,2}=\left(
	\begin{array}{ccccccccc}
	0 & \alpha_2& \alpha_3&0&\alpha_5&1&\alpha_3&1&0\\
	\beta_1 & \beta_2& \beta_3&\beta_4&\beta_5&\beta_8&\alpha_3\beta_8&\beta_8&0\\
	\gamma_1 & \gamma_2& \beta_2&\gamma_4&\gamma_5&1+\beta_5&1+\beta_4&\alpha_2+\beta_5&\alpha_3+\beta_8\end{array}\right),\] where 
	$\alpha_3(\beta_3+\alpha_3\beta_8)\neq 0$,
	\[A_{7,2}=\left(
	\begin{array}{ccccccccc}
	\alpha_1 & 0& \alpha_3&0&\alpha_5&1&\alpha_3&1&0\\
	\beta_1 & \beta_2& \alpha_3\beta_8&\beta_4&\beta_5&\beta_8&\alpha_3\beta_8&\beta_8&0\\
	\gamma_1 & \gamma_2& \alpha_1+\beta_2&\gamma_4&\gamma_5&1+\beta_5&1+\alpha_1+\beta_4&\beta_5&\alpha_3+\beta_8\end{array}\right),\] where 
	$\alpha_3\neq 0, 1$,
	\[A_{8,2}=\left(
	\begin{array}{ccccccccc}
	\alpha_1 & \alpha_2& 0&0&\alpha_5&1&0&1&0\\
	\beta_1 & \alpha_1& \beta_3&\beta_4&\beta_5&\beta_8&0&\beta_8&0\\
	\gamma_1 & \gamma_2& 0&\gamma_4&\gamma_5&1+\beta_5&1+\alpha_1+\beta_4&\alpha_2+\beta_5&\beta_8\end{array}\right),\] where 
	$\beta_3\neq 0$,
	\[A_{9,2}=\left(
	\begin{array}{ccccccccc}
	\alpha_1 & \alpha_2& 0&0&\alpha_5&1&0&1&0\\
	\beta_1 & \beta_2& 0&\beta_4&\beta_5&\beta_8&0&\beta_8&0\\
	\gamma_1 & 0& \alpha_1+\beta_2&\gamma_4&\gamma_5&1+\beta_5&1+\alpha_1+\beta_4&\alpha_2+\beta_5&\beta_8\end{array}\right),\] where 
	$\beta_2\neq 0$,
	\[A_{10,2}=\left(
	\begin{array}{ccccccccc}
	\alpha_1 & \alpha_2& 0&0&\alpha_5&1&0&1&0\\
	\beta_1 & 0& 0&\beta_4&\beta_5&\beta_8&0&\beta_8&0\\
	0 & \gamma_2& \alpha_1&\gamma_4&\gamma_5&1+\beta_5&1+\alpha_1+\beta_4&\alpha_2+\beta_5&\beta_8\end{array}\right),\] where 
	$\beta_1\neq 0$,
	\[A_{11,2}=\left(
	\begin{array}{ccccccccc}
	\alpha_1 & \alpha_2& 0&0&\alpha_5&1&0&1&0\\
	0 & 0& 0&\beta_4&\beta_5&\beta_8&0&\beta_8&0\\
	\gamma_1 & \gamma_2& \alpha_1&0&\gamma_5&11\beta_5&1+\alpha_1+\beta_4&\alpha_2+\beta_5&\beta_8\end{array}\right),\] where 
	$\alpha_1\neq 1$,
	\[A_{12,2}=\left(
	\begin{array}{ccccccccc}
	1 & \alpha_2& 0&0&\alpha_5&1&0&1&0\\
	0 & 0& 0&\beta_4&\beta_5&\beta_8&0&\beta_8&0\\
	\gamma_1 & \gamma_2& 1&\gamma_4&0&1+\beta_5&\beta_4&\alpha_2+\beta_5&0\end{array}\right),\] where $\beta_8\neq 0$ and if $\beta_8=0$ then $1+\alpha_2+\beta_5\neq 0,$
	\[A_{13,2}=\left(
	\begin{array}{ccccccccc}
	1 & \alpha_2& 0&0&\alpha_5&1&0&1&0\\
	0 & 0& 0&\beta_4&1+\alpha_2&0&0&0&0\\
	\gamma_1 & \gamma_2& 1&\gamma_4&\gamma_5&\alpha_2&\beta_4&1&0\end{array}\right),\] 
	\[A_{14,2}=\left(
	\begin{array}{ccccccccc}
	\alpha_1 & \alpha_2& 1&0&\alpha_5&1&1&1&0\\
	\beta_1 & \beta_2& 1&\beta_4&\beta_5&1&1&1&0\\
	0 & \gamma_2& \alpha_1+\beta_2&\gamma_4&\gamma_5&1+\beta_5&1+\alpha_1+\beta_4&\alpha_2+\beta_5&0\end{array}\right),\] 
	where $1+\alpha_1+\beta_1+\beta_2+\beta_4\neq 0 $,
	\[A_{15,2}=\left(
	\begin{array}{ccccccccc}
	\alpha_1 & \alpha_2& 1&0&\alpha_5&1&1&1&0\\
	\beta_1 & \beta_2& 1&\beta_4&\beta_5&1&1&1&0\\
	\gamma_1 & 0& \alpha_1+\beta_2&\gamma_4&\gamma_5&1+\beta_5&\beta_1+\beta_2&\alpha_2+\beta_5&0\end{array}\right),\] 
	where $\beta_4=1+\alpha_1+\beta_1+\beta_2$, $\alpha_1+\beta_1+\beta_2+\beta_5=1+\beta_4+\beta_5\neq 0 $,
	\[A_{16,2}=\left(
	\begin{array}{ccccccccc}
	\alpha_1 & \alpha_2& 1&0&\alpha_5&1&1&1&0\\
	\beta_1 & \beta_2& 1&\beta_4&\beta_5&1&1&1&0\\
	\gamma_1 & \gamma_2& \alpha_1+\beta_2&0&\gamma_5&1+\beta_5&\beta_1+\beta_2&\alpha_2+\beta_5&0\end{array}\right),\] 
	where $\beta_5=\alpha_1+\beta_1+\beta_2$, $\beta_4=1+\alpha_1+\beta_1+\beta_2$, $\alpha_1+\beta_5\neq 0 $,
	\[A_{17,2}=\left(
	\begin{array}{ccccccccc}
	\alpha_1 & \alpha_2& 1&0&\alpha_5&1&1&1&0\\
	\beta_1 & \beta_1& 1&1+\alpha_1&\alpha_1&1&1&1&0\\
	\gamma_1 & \gamma_2& \alpha_1+\beta_1&\gamma_4&0&1+\alpha_1&0&\alpha_1+\alpha_2&0\end{array}\right),\] 
	where  $\alpha_1+\alpha_2+\alpha_5\neq 1 $,
	\[A_{18,2}=\left(
	\begin{array}{ccccccccc}
	\alpha_1 & \alpha_2& 1&0&1+\alpha_1+\alpha_2&1&1&1&0\\
	\beta_1 & \beta_1& 1&1+\alpha_1&\alpha_1&1&1&1&0\\
	\gamma_1 & \gamma_2& \alpha_1+\beta_1&\gamma_4&\gamma_5&1+\alpha_1&0&\alpha_1+\alpha_2&0\end{array}\right).\] 
	\[A_{19,2}=\left(
	\begin{array}{ccccccccc}
	\alpha_1 & \alpha_2& 1&0&\alpha_5&1&1&1&0\\
	\beta_1 & \beta_2& \beta_8&\beta_4&\beta_5&\beta_8&\beta_8&\beta_8&0\\
	\gamma_1 & \gamma_1& \alpha_1+\beta_2&\gamma_4&\gamma_5&1+\beta_5&1+\alpha_1+\beta_4&\alpha_2+\beta_5&1+\beta_8\end{array}\right),\] 
	where $\beta_8\neq 1$, $\beta_4+\beta_5\neq 1 $,
	\[A_{20,2}=\left(
	\begin{array}{ccccccccc}
	\alpha_1 & \alpha_2& 1&0&\alpha_5&1&1&1&0\\
	\beta_1 & \beta_2& \beta_8&\beta_4&1+\beta_4&\beta_8&\beta_8&\beta_8&0\\
	\gamma_1 & \gamma_2& \alpha_1+\beta_2&\gamma_1&\gamma_5&\beta_4&1+\alpha_1+\beta_4&1+\alpha_2+\beta_4&1+\beta_8\end{array}\right),\] 
	where $\beta_8\neq 1$, $\beta_2\neq \beta_1 $,
	\[A_{21,2}=\left(
	\begin{array}{ccccccccc}
	\alpha_1 & \alpha_2& 1&0&\alpha_5&1&1&1&0\\
	\beta_1 & \beta_1& \beta_8&\beta_4&1+\beta_4&\beta_8&\beta_8&\beta_8&0\\
	\gamma_1 & \gamma_2& \alpha_1+\beta_2&\gamma_4&\gamma_1&\beta_4&1+\alpha_1+\beta_4&1+\alpha_2+\beta_4&1+\beta_8\end{array}\right),\] 
	where $\beta_8\neq 1$, $\alpha_1+\alpha_2+\alpha_5\neq 1$,
	\[A_{22,2}=\left(
	\begin{array}{ccccccccc}
	\alpha_1 & \alpha_2& 1&0&1+\alpha_1+\alpha_2&1&1&1&0\\
	\beta_1 & \beta_1& \beta_8&\beta_4&1+\beta_4&\beta_8&\beta_8&\beta_8&0\\
	0 & \gamma_2& \alpha_1+\beta_2&\gamma_4&\gamma_5&\beta_4&1+\alpha_1+\beta_4&1+\alpha_2+\beta_4&1+\beta_8\end{array}\right),\] 
	where $\beta_8\neq 1$,
	\[A_{23,2}=\left(
	\begin{array}{ccccccccc}
	\alpha_1 & \alpha_2& \alpha_3&\alpha_4&\alpha_5&0&1+\alpha_3+\beta_8&\alpha_8&0\\
	\beta_1 & \beta_2& \beta_3&0&0&1&\beta_7&\beta_8&0\\
	\gamma_1 & \gamma_2& \alpha_1+\beta_2&\gamma_4&\gamma_5&1+\alpha_4&1+\alpha_1&\alpha_2&\alpha_3+1\end{array}\right),\] where $\beta_8\neq 1$,
	\[A_{24,2}=\left(
	\begin{array}{ccccccccc}
	\alpha_1 & \alpha_2& \alpha_3&\alpha_4&0&0&2+\alpha_3&\alpha_8&0\\
	\beta_1 & \beta_2& \beta_3&0&\beta_5&1&\beta_7&1&0\\
	\gamma_1 & \gamma_2& \alpha_1+\beta_2&\gamma_4&\gamma_5&1+\alpha_4+\beta_5&1+\alpha_1&\alpha_2+\beta_5&\alpha_3+1\end{array}\right),\] where $\alpha_8\neq 0$,
	\[A_{25,2}=\left(
	\begin{array}{ccccccccc}
	\alpha_1 & \alpha_2& \alpha_3&0&\alpha_5&0&\alpha_3&0&0\\
	\beta_1 & \beta_2& \beta_3&0&\beta_5&1&\beta_7&1&0\\
	\gamma_1 & \gamma_2& \alpha_1+\beta_2&\gamma_4&\gamma_5&1+\beta_5&1+\alpha_1&\alpha_2+\beta_5&\alpha_3+1\end{array}\right),\] where $\alpha_3\neq 0$,
	\[A_{26,2}=\left(
	\begin{array}{ccccccccc}
	\alpha_1 & \alpha_2& 0&\alpha_4&\alpha_5&0&0&0&0\\
	\beta_1 & 0& \beta_3&0&\beta_5&1&\beta_7&1&0\\
	\gamma_1 & \gamma_2& \alpha_1&\gamma_4&\gamma_5&1+\alpha_4+\beta_5&1+\alpha_1&\alpha_2+\beta_5&1\end{array}\right),\] where $\beta_7(\beta_3+1)\neq 0$,
	\[A_{27,2}=\left(
	\begin{array}{ccccccccc}
	\alpha_1 & \alpha_2& 0&\alpha_4&\alpha_5&0&0&0&0\\
	\beta_1 & \beta_2& \beta_3&0&\beta_5&1&0&1&0\\
	0 & \gamma_2& \alpha_1+\beta_2&\gamma_4&\gamma_5&1+\alpha_4+\beta_5&1+\alpha_1&\alpha_2+\beta_5&1\end{array}\right),\] where $\beta_1\neq 0$,
	\[A_{28,2}=\left(
	\begin{array}{ccccccccc}
	\alpha_1 & \alpha_2& 0&\alpha_4&\alpha_5&0&0&0&0\\
	0 & \beta_2& \beta_3&0&\beta_5&1&0&1&0\\
	\gamma_1 & 0& \alpha_1+\beta_2&\gamma_4&\gamma_5&1+\alpha_4+\beta_5&1+\alpha_1&\alpha_2+\beta_5&1\end{array}\right),\] where $\beta_2\neq 0$,
	\[A_{29,2}=\left(
	\begin{array}{ccccccccc}
	\alpha_1 & \alpha_2& 0&\alpha_4&\alpha_5&0&0&0&0\\
	0 & 0& \beta_3&0&\beta_5&1&0&1&0\\
	\gamma_1 & \gamma_2& \alpha_1+\beta_2&0&\gamma_5&1+\alpha_4+\beta_5&1+\alpha_1&\alpha_2+\beta_5&1\end{array}\right),\] where $\alpha_1\neq 1$,
	\[A_{30,2}=\left(
	\begin{array}{ccccccccc}
	1 & \alpha_2& 0&\alpha_4&\alpha_5&0&0&0&0\\
	0 & 0& \beta_3&0&\beta_5&1&0&1&0\\
	\gamma_1 & \gamma_2& 1&\gamma_4&0&1+\alpha_4+\beta_5&0&\alpha_2+\beta_5&1\end{array}\right),\] 
	\[A_{31,2}=\left(
	\begin{array}{ccccccccc}
	\alpha_1 & \alpha_2& 0&\alpha_4&\alpha_5&0&0&0&0\\
	0 & \beta_2& 1&0&\beta_5&1&\beta_7&1&0\\
	\gamma_1 & \gamma_2& \alpha_1+\beta_2&\gamma_4&\gamma_5&1+\alpha_4+\beta_5&1+\alpha_1&\alpha_2+\beta_5&1\end{array}\right),\] where $\beta_7\neq 0,1$,
	\[A_{32,2}=\left(
	\begin{array}{ccccccccc}
	\alpha_1 & \alpha_2& 0&\alpha_4&\alpha_5&0&0&0&0\\
	\beta_1 & \beta_2& 1&0&\beta_5&1&1&1&0\\
	\gamma_1 & \gamma_1& \alpha_1+\beta_2&\gamma_4&\gamma_5&1+\alpha_4+\beta_5&1+\alpha_1&\alpha_2+\beta_5&1\end{array}\right),\] where $\beta_5\neq \beta_1+\beta_2$,
	\[A_{33,2}=\left(
	\begin{array}{ccccccccc}
	\alpha_1 & \alpha_2& 0&\alpha_4&\alpha_5&0&0&0&0\\
	\beta_1 & \beta_2& 1&0&\beta_5&1&1&1&0\\
	\gamma_1 & \gamma_2& \alpha_1+\beta_2&\gamma_1&\gamma_5&1+\alpha_4+\beta_5&1+\alpha_1&\alpha_2+\beta_5&1\end{array}\right),\] where $\beta_5= \beta_1+\beta_2$, $\alpha_4\neq 1$,
	\[A_{34,2}=\left(
	\begin{array}{ccccccccc}
	\alpha_1 & \alpha_2& 0&1&\alpha_5&0&0&0&0\\
	\beta_1 & \beta_2& 1&0&\beta_5&1&1&1&0\\
	\gamma_1 & \gamma_2& \alpha_1+\beta_2&\gamma_4&\gamma_1&\beta_5&1+\alpha_1&\alpha_2+\beta_5&1\end{array}\right),\] where $\beta_5= \beta_1+\beta_2$, $\alpha_5\neq 1+\alpha_1+\alpha_2$,
	\[A_{35,2}=\left(
	\begin{array}{ccccccccc}
	\alpha_1 & \alpha_2& 0&1&\alpha_5&0&0&0&0\\
	\beta_1 & \beta_2& 1&0&\beta_5&1&1&1&0\\
	0 & \gamma_2& \alpha_1+\beta_2&\gamma_1&\gamma_5&\beta_5&1+\alpha_1&\alpha_2+\beta_5&1\end{array}\right),\] where $\alpha_5=1+\alpha_1+\alpha_2$, $\beta_5= \beta_1+\beta_2$,
	\[A_{36,2}=\left(
	\begin{array}{ccccccccc}
	0 & \alpha_2& 1&\alpha_4&\alpha_5&0&\alpha_7&\alpha_8&0\\
	\beta_1 & \beta_2& \beta_3&\beta_4&1+\alpha_4&0&\beta_7&1+\alpha_7&0\\
	\gamma_1 & \gamma_2& \beta_2&\gamma_4&\gamma_5&0&1+\beta_4&1+\alpha_2+\alpha_4&1\end{array}\right),\] where $\alpha_7\neq 1$,
	\[A_{37,2}=\left(
	\begin{array}{ccccccccc}
	\alpha_1 & \alpha_2& 1&\alpha_4&\alpha_5&0&1&\alpha_8&0\\
	\beta_1 & 0& \beta_3&\beta_4&1+\alpha_4&0&\beta_7&0&0\\
	\gamma_1 & \gamma_2& \alpha_1&\gamma_4&\gamma_5&0&1+\alpha_1+\beta_4&1+\alpha_2+\alpha_4&1\end{array}\right),\] where $\beta_3\neq 0$,
	\[A_{38,2}=\left(
	\begin{array}{ccccccccc}
	\alpha_1 & \alpha_2& 1&\alpha_4&\alpha_5&0&1&\alpha_8&0\\
	0 & \beta_2& 0&\beta_4&1+\alpha_4&0&\beta_7&0&0\\
	\gamma_1 & \gamma_2& \alpha_1+\beta_2&\gamma_4&\gamma_5&0&1+\alpha_1+\beta_4&1+\alpha_2+\alpha_4&1\end{array}\right),\] where $\beta_7\neq 0$,
	\[A_{39,2}=\left(
	\begin{array}{ccccccccc}
	\alpha_1 & 0& 1&\alpha_4&\alpha_5&0&1&\alpha_8&0\\
	\beta_1 & \beta_2& 0&\beta_4&1+\alpha_4&0&0&0&0\\
	\gamma_1 & \gamma_2& \alpha_1+\beta_2&\gamma_4&\gamma_5&0&1+\alpha_1+\beta_4&1+\alpha_4&1\end{array}\right),\] where $\alpha_8\neq 1$,
	\[A_{40,2}=\left(
	\begin{array}{ccccccccc}
	\alpha_1 &1+\alpha_4& 1&\alpha_4&\alpha_5&0&1&1&0\\
	\beta_1 & \beta_2& 0&\beta_4&1+\alpha_4&0&0&0&0\\
	\gamma_1 & \gamma_2& \alpha_1+\beta_2&\gamma_4&\gamma_5&0&1+\alpha_1+\beta_4&0&1\end{array}\right),\]
	\[A_{41,2}=\left(
	\begin{array}{ccccccccc}
	\alpha_1 & 0& 0&\alpha_4&\alpha_5&0&\beta_8&\alpha_8&0\\
	\beta_1 & \alpha_1& 1&\beta_4&\beta_5&0&\beta_7&\beta_8&0\\
	\gamma_1 & \gamma_2& 0&\gamma_4&\gamma_5&1+\alpha_4+\beta_5&1+\alpha_1+\beta_4&\beta_5&0\end{array}\right),\] where $\alpha_8\neq 0$, 
	\[A_{42,2}=\left(
	\begin{array}{ccccccccc}
	0 & \alpha_2& 0&\alpha_4&\alpha_5&0&\beta_8&0&0\\
	\beta_1 & 0& 1&\beta_4&\beta_5&0&\beta_7&\beta_8&0\\
	\gamma_1 & \gamma_2& 0&\gamma_4&\gamma_5&1+\alpha_4+\beta_5&1+\beta_4&\alpha_2+\beta_5&0\end{array}\right),\] where $\alpha_7=\beta_8\neq 0$,
	\[A_{43,2}=\left(
	\begin{array}{ccccccccc}
	\alpha_1 & \alpha_2& 0&\alpha_4&\alpha_5&0&0&0&0\\
	\beta_1 & 0& 1&\beta_4&\beta_5&0&\beta_7&0&0\\
	\gamma_1 & \gamma_2& \alpha_1&\gamma_4&\gamma_5&1+\alpha_4+\beta_5&1+\alpha_1+\beta_4&\alpha_2=\beta_5&0\end{array}\right),\]
	\[A_{44,2}=\left(
	\begin{array}{ccccccccc}
	\alpha_1 & 0& 0&\alpha_4&0&0&\alpha_7&1&0\\
	\beta_1 & \beta_2& 0&\beta_4&\beta_5&0&\beta_7&\alpha_7&0\\
	\gamma_1 & \gamma_2& \alpha_1+\beta_2&\gamma_4&\gamma_5&1+\alpha_4+\beta_5&1+\alpha_1+\beta_4&\beta_5&0\end{array}\right),\] 
	\[A_{45,2}=\left(
	\begin{array}{ccccccccc}
	0 & \alpha_2& 0&0&\alpha_5&0&1&0&0\\
	\beta_1 & \beta_2& 0&\beta_4&\beta_5&0&\beta_7&1&0\\
	\gamma_1 & \gamma_2& \beta_2&\gamma_4&\gamma_5&1+\beta_5&1+\beta_4&\alpha_2+\beta_5&0\end{array}\right),\] \[A_{46,2}=\left(
	\begin{array}{ccccccccc}
	\alpha_1 & \alpha_2& 0&\alpha_4&\alpha_5&0&0&0&0\\
	0 & \beta_2& 0&0&\beta_5&0&1&0&0\\
	\gamma_1 & \gamma_2& \alpha_1+\beta_2&\gamma_4&\gamma_5&1+\alpha_4+\beta_5&1+\alpha_1&\alpha_2+\beta_5&0\end{array}\right),\] 
	\[A_{47,2}=\left(
	\begin{array}{ccccccccc}
	\alpha_1 & \alpha_2& 0&\alpha_4&\alpha_5&0&0&0&0\\
	\beta_1 & \beta_2& 0&\beta_4&\beta_5&0&0&0&0\\
	0 & 0& \alpha_1+\beta_2&0&0&1+\alpha_4+\beta_5&1+\alpha_1+\beta_4&\alpha_2+\beta_5&0\end{array}\right),\] 
	\[A_{48,2}=\left(
	\begin{array}{ccccccccc}
	\alpha_1 & \alpha_2& 0&\alpha_4&\alpha_5&0&0&0&0\\
	\beta_1 & \beta_2& 0&\beta_4&\beta_5&0&0&0&0\\
	1 & 0& \alpha_1+\beta_2&0&0&1+\alpha_4+\beta_5&1+\alpha_1+\beta_4&\alpha_2+\beta_5&0\end{array}\right),\] where $rk(M_{(1,2,3,4)})=rk(M_{(2,3,4)})$.
	\[A_{49,2}=\left(
	\begin{array}{ccccccccc}
	\alpha_1 & \alpha_2& 0&\alpha_4&\alpha_5&0&0&0&0\\
	\beta_1 & \beta_2& 0&\beta_4&\beta_5&0&0&0&0\\
	0 & 1& \alpha_1+\beta_2&0&0&1+\alpha_4+\beta_5&1+\alpha_1+\beta_4&\alpha_2+\beta_5&0\end{array}\right),\] where $rk(M_{(1,2,3,4)})=rk(M_{(1,3,4)})$, $rk(M_{(3,4)})=rk(M_{(2,3,4)})$,
	\[A_{50,2}=\left(
	\begin{array}{ccccccccc}
	\alpha_1 & \alpha_2& 0&\alpha_4&\alpha_5&0&0&0&0\\
	\beta_1 & \beta_2& 0&\beta_4&\beta_5&0&0&0&0\\
	0 & 0& \alpha_1+\beta_2&1&0&1+\alpha_4+\beta_5&1+\alpha_1+\beta_4&\alpha_2+\beta_5&0\end{array}\right),\] where $rk(M_{(1,2,3,4)})=rk(M_{(1,2,4)})$, $rk(M_{(2,4)})=rk(M_{(2,3,4)})$, $rk(M_{(1,4)})=rk(M_{(1,3,4)})$,
	\[A_{51,2}=\left(
	\begin{array}{ccccccccc}
	\alpha_1 & \alpha_2& 0&\alpha_4&\alpha_5&0&0&0&0\\
	\beta_1 & \beta_2& 0&\beta_4&\beta_5&0&0&0&0\\
	0 & 0& \alpha_1+\beta_2&0&1&1+\alpha_4+\beta_5&1+\alpha_1+\beta_4&\alpha_2+\beta_5&0\end{array}\right),\] where $rk(M_{(1,2,3,4)})=rk(M_{(1,2,3)})$, $rk(M_{(2,3)})=rk(M_{(2,3,4)})$, $rk(M_{(1,3)})=rk(M_{(1,3,4)})$, $rk(M_{(1,2)})=rk(M_{(1,2,4)})$,
	\[A_{52,2}=\left(
	\begin{array}{ccccccccc}
	\alpha_1 & \alpha_2& 0&\alpha_4&\alpha_5&0&0&0&0\\
	\beta_1 & \beta_2& 0&\beta_4&\beta_5&0&0&0&0\\
	1 & \gamma_2& \alpha_1+\beta_2&0&0&1+\alpha_4+\beta_5&1+\alpha_1+\beta_4&\alpha_2+\beta_5&0\end{array}\right),\] where $rk(M_{(1,2,3,4)})=rk(M_{(2-\gamma_2\times1,3,4)})$, $rk(M_{(3,4)})=rk(M_{(2,3,4)})$, $rk(M_{(3,4)})=rk(M_{(1,3,4)})$, $rk(M_{(2-\gamma_2\times 1,4)})=rk(M_{(1,2,4)})$, 
	$rk(M_{(2-\gamma_2\times 1,3)})=rk(M_{(1,2,3)})$ and $\gamma_2\neq 0$,
	\[A_{53,2}=\left(
	\begin{array}{ccccccccc}
	\alpha_1 & \alpha_2& 0&\alpha_4&\alpha_5&0&0&0&0\\
	\beta_1 & \beta_2& 0&\beta_4&\beta_5&0&0&0&0\\
	1 & 0& \alpha_1+\beta_2&\gamma_4&0&1+\alpha_4+\beta_5&1+\alpha_1+\beta_4&\alpha_2+\beta_5&0\end{array}\right),\] where $rk(M_{(1,2,3,4)})=rk(M_{(2,3-\gamma_4\times 1,4)})$, $rk(M_{(2,4)})=rk(M_{(2,3,4)})$, $rk(M_{(3-\gamma_4\times 1,4)})=rk(M_{(1,3,4)})$, $rk(M_{(2,4)})=rk(M_{(1,2,4)})$, 
	$rk(M_{(2,3-\gamma_4\times 1)})=rk(M_{(1,2,3)})$, $rk(M_{(4)})=rk(M_{(3,4)})$ and $\gamma_4\neq 0$,
	\[A_{54,2}=\left(
	\begin{array}{ccccccccc}
	\alpha_1 & \alpha_2& 0&\alpha_4&\alpha_5&0&0&0&0\\
	\beta_1 & \beta_2& 0&\beta_4&\beta_5&0&0&0&0\\
	0 & 1& \alpha_1+\beta_2&\gamma_4&0&1+\alpha_4+\beta_5&1+\alpha_1+\beta_4&\alpha_2+\beta_5&0\end{array}\right),\] where $rk(M_{(1,2,3,4)})=rk(M_{(1,3-\gamma_4\times 2,4)})$, $rk(M_{(3-\gamma_4\times 2,4)})=rk(M_{(2,3,4)})$, $rk(M_{(1,4)})=rk(M_{(1,3,4)})$, $rk(M_{(1,4)})=rk(M_{(1,2,4)})$, 
	$rk(M_{(1,3-\gamma_4\times 2)})=rk(M_{(1,2,3)})$, $rk(M_{(4)})=rk(M_{(3,4)})$,  $rk(M_{(4)})=rk(M_{(2,4)})$ and $\gamma_4\neq 0$,
	\[A_{55,2}=\left(
	\begin{array}{ccccccccc}
	\alpha_1 & \alpha_2& 0&\alpha_4&\alpha_5&0&0&0&0\\
	\beta_1 & \beta_2& 0&\beta_4&\beta_5&0&0&0&0\\
	1 & 0& \alpha_1+\beta_2&0&\gamma_5&1+\alpha_4+\beta_5&1+\alpha_1+\beta_4&\alpha_2+\beta_5&0\end{array}\right),\] where $rk(M_{(1,2,3,4)})=rk(M_{(2,3,4-\gamma_5\times 1)})$, $rk(M_{(2,3)})=rk(M_{(2,3,4)})$, $rk(M_{(3,4-\gamma_5\times 1)})=rk(M_{(1,3,4})$, $rk(M_{(2,4-\gamma_5\times 1)})=rk(M_{(1,2,4)})$, 
	$rk(M_{(2,3)})=rk(M_{(1,2,3)})$, $rk(M_{(3)})=rk(M_{(3,4)})$, $rk(M_{(2)})=rk(M_{(2,4)})$, $rk(M_{(4-\gamma_5\times 1)})=rk(M_{(1,4)})$ and $\gamma_5\neq 0$,
	\[A_{56,2}=\left(
	\begin{array}{ccccccccc}
	\alpha_1 & \alpha_2& 0&\alpha_4&\alpha_5&0&0&0&0\\
	\beta_1 & \beta_2& 0&\beta_4&\beta_5&0&0&0&0\\
	0 & 1& \alpha_1+\beta_2&0&\gamma_5&1+\alpha_4+\beta_5&1+\alpha_1+\beta_4&\alpha_2+\beta_5&0\end{array}\right),\] where $rk(M_{(1,2,3,4)})=rk(M_{(1,3,4-\gamma_5\times 2)})$, $rk(M_{(3,4-\gamma_5\times 2)})=rk(M_{(2,3,4)})$, $rk(M_{(1,3)})=rk(M_{(1,3,4})$, $rk(M_{(1,4-\gamma_5\times 2)})=rk(M_{(1,2,4)})$, 
	$rk(M_{(1,3)})=rk(M_{(1,2,3)})$, $rk(M_{(3)})=rk(M_{(3,4)})$, $rk(M_{(4-\gamma_5\times 2)})=rk(M_{(2,4)})$, $rk(M_{(1)})=rk(M_{(1,4)})$, $rk(M_{(3)})=rk(M_{(2,3)})$ and $\gamma_5\neq 0$,
	\[A_{57,2}=\left(
	\begin{array}{ccccccccc}
	\alpha_1 & \alpha_2& 0&\alpha_4&\alpha_5&0&0&0&0\\
	\beta_1 & \beta_2& 0&\beta_4&\beta_5&0&0&0&0\\
	0 & 0& \alpha_1+\beta_2&1&\gamma_5&1+\alpha_4+\beta_5&1+\alpha_1+\beta_4&\alpha_2+\beta_5&0\end{array}\right),\] where $rk(M_{(1,2,3,4)})=rk(M_{(1,2,4-\gamma_5\times 3)})$, $rk(M_{(2,4-\gamma_5\times 3)})=rk(M_{(2,3,4)})$, $rk(M_{(1,4-\gamma_5\times 3)})=rk(M_{(1,3,4})$, $rk(M_{(1,2)})=rk(M_{(1,2,4)})$, 
	$rk(M_{(1,2)})=rk(M_{(1,2,3)})$, $rk(M_{(4-\gamma_5\times 3)})=rk(M_{(3,4)})$, $rk(M_{(2)})=rk(M_{(2,4)})$, $rk(M_{(1)})=rk(M_{(1,4)})$, $rk(M_{(1)})=rk(M_{(1,3)})$ and $\gamma_5\neq 0$,
	\[A_{58,2}=\left(
	\begin{array}{ccccccccc}
	\alpha_1 & \alpha_2& 0&\alpha_4&\alpha_5&0&0&0&0\\
	\beta_1 & \beta_2& 0&\beta_4&\beta_5&0&0&0&0\\
	0 & 1& \alpha_1+\beta_2&\gamma_4&\gamma_5&1+\alpha_4+\beta_5&1+\alpha_1+\beta_4&\alpha_2+\beta_5&0\end{array}\right),\] where $rk(M_{(1,2,3,4)})=rk(M_{(1,3-\gamma_4\times 2,4-\gamma_5\times 2)})$, $rk(M_{(3-\gamma_4\times 2, 4-\gamma_5\times 2)})=rk(M_{(2,3,4)})$, $rk(M_{(1,\gamma_4\times 4-\gamma_5\times 3)})=rk(M_{(1,3,4})$, $rk(M_{(1,4-\gamma_5\times 2)})=rk(M_{(1,2,4)})$, 
	$rk(M_{(1,3-\gamma_4\times 2)})=rk(M_{(1,2,3)})$, $rk(M_{(\gamma_4\times 4-\gamma_5\times 3)})=rk(M_{(3,4)})$, $rk(M_{(4-\gamma_5\times 2)})=rk(M_{(2,4)})$, $rk(M_{(1)})=rk(M_{(1,4)})$, $rk(M_{(1)})=rk(M_{(1,3)})$, $rk(M_{(1)})=rk(M_{(1,2)})$ and $\gamma_4\neq 0$, $\gamma_5\neq 0$,
	\[A_{59,2}=\left(
	\begin{array}{ccccccccc}
	\alpha_1 & \alpha_2& 0&\alpha_4&\alpha_5&0&0&0&0\\
	\beta_1 & \beta_2& 0&\beta_4&\beta_5&0&0&0&0\\
	1 & 0& \alpha_1+\beta_2&\gamma_4&\gamma_5&1+\alpha_4+\beta_5&1+\alpha_1+\beta_4&\alpha_2+\beta_5&0\end{array}\right),\] where $rk(M_{(1,2,3,4)})=rk(M_{(2,3-\gamma_4\times 1,4-\gamma_5\times 1)})$, $rk(M_{(2,\gamma_4\times 3-\gamma_5\times 3)})=rk(M_{(2,3,4)})$, $rk(M_{(3-\gamma_4\times 1,4-\gamma_5\times 1)})=rk(M_{(1,3,4})$, $rk(M_{(2,4-\gamma_5\times 1)})=rk(M_{(1,2,4)})$, 
	$rk(M_{(2,3-\gamma_4\times 1)})=rk(M_{(1,2,3)})$, $rk(M_{(\gamma_4\times 4-\gamma_5\times 3)})=rk(M_{(3,4)})$, $rk(M_{(2)})=rk(M_{(2,4)})$, $rk(M_{(4-\gamma_4\times 1)})=rk(M_{(1,4)})$, $rk(M_{(3-\gamma_5\times 1)})=rk(M_{(1,3)})$, $rk(M_{(2)})=rk(M_{(1,2)})$, $rk(M_{(1)})=0$ and $\gamma_4\neq 0$, $\gamma_5\neq 0$,
	\[A_{60,2}=\left(
	\begin{array}{ccccccccc}
	\alpha_1 & \alpha_2& 0&\alpha_4&\alpha_5&0&0&0&0\\
	0 & 0& 0&1+\alpha_1&\alpha_1&0&0&0&0\\
	1 & \gamma_2& \alpha_1&0&\gamma_5&1+\alpha_1+\alpha_4&0&\alpha_1+\alpha_2&0\end{array}\right),\] where $\gamma_2\neq 0$, $\gamma_5\neq 0$,
	\[A_{61,2}=\left(
	\begin{array}{ccccccccc}
	1 & \alpha_2& 0&\alpha_4&\alpha_5&0&0&0&0\\
	0 & 0& 0&0&1&0&0&0&0\\
	1 & \gamma_2& 1&\gamma_4&0&\alpha_4&0&1+\alpha_2&0\end{array}\right),\] where $\gamma_2\neq 0$, $\gamma_4\neq 0$,
	\[A_{62,2}=\left(
	\begin{array}{ccccccccc}
	1 & \alpha_2& 0&\alpha_2&0&0&0&0&0\\
	0 & 0& 0&0&1&0&0&0&0\\
	1 & \gamma_2& 1&\gamma_4&\gamma_5&\alpha_2 &0 &1+\alpha_2&0\end{array}\right),\] where 
	$\gamma_2\neq 0,\ \gamma_4\neq 0,\ \gamma_5\neq 0.$
\end{theorem}

{\bf Proof.} 
Due to the $3^{rd}$ column in (\ref{E4}) one can conclude:

\underline{\textbf{Case 1.\ $\alpha_9\neq 0$}}. It is possible to make $\alpha'_9=1$ ($c^2\alpha_9=1$), $\alpha'_7=\alpha'_8=0$ to get  
\[A_{1,2}=\left(
\begin{array}{ccccccccc}
\alpha_1 & \alpha_2& \alpha_3&\alpha_4&\alpha_5&\alpha_6&0&0&1\\
\beta_1 & \beta_2& \beta_3&\beta_4&\beta_5&\beta_8-\alpha_3&\beta_7&\beta_9&\beta_9\\
\gamma_1 & \gamma_2& -\alpha_1-\beta_2&\gamma_4&\gamma_5&1-\alpha_4-\beta_5&1-\alpha_1-\beta_4&-\alpha_2-\beta_5&-\beta_8\end{array}\right). \]
\underline{\textbf{Case 2.\ $\alpha_9= 0$}, $\beta_9\neq 0$}. One can make $\beta'_9=1$, $\beta'_7=\beta'_8=0$ to get  
\[A_{2,2}=\left(
\begin{array}{ccccccccc}
\alpha_1 & \alpha_2& \alpha_3&\alpha_4&\alpha_5&\alpha_6&\alpha_7&\alpha_8&0\\
\beta_1 & \beta_2&\beta_3&\beta_4&\beta_5&\alpha_7-\alpha_3&0&0&1\\
\gamma_1 & \gamma_2& -\alpha_1-\beta_2&\gamma_4&\gamma_5&1-\alpha_4-\beta_5&1-\alpha_1-\beta_4&-\alpha_2-\beta_5&-\alpha_7\end{array}\right).\]

If $\alpha_9=\beta_9= 0$, $\alpha_6\neq 0$ then due to the $3^{rd}$ column in (\ref{E3}) it is possible making $\alpha'_6=1$ by appropriate choose of $c$ to come to the case

\underline{\textbf{Case 3.\ $\alpha_9=\beta_9= 0$, $\alpha_6= 1 $}} with respect to $g$ of the form $g^{-1}=\left(
\begin{array}{ccc}
1 & 0& 0\\ 0 & 1& 0\\ a & b& 1\end{array}\right)$. If $a=-\alpha_4-b\alpha_7$ then $\alpha'_4=0$ and 
$\mathcal{A}'_1=g(\mathcal{A}_1+b\mathcal{A}_3)g^{-1}=$ 
\[\left( 
\begin{array}{ccc}
\alpha_1-\alpha_4(\alpha_3+\alpha_7)-b\alpha_7(\alpha_3+\alpha_7)&\alpha_2-\alpha_4(\alpha_3+\alpha_8)-b\alpha_7(\alpha_3+\alpha_8)&\alpha_3\\
\beta_1-\alpha_4(\beta_3+\beta_7)-b\alpha_7(\beta_3+\beta_7)&\beta_2-\alpha_4(\beta_3+\beta_8)-b\alpha_7(\beta_3+\beta_8)&\beta_3\\
*&*&*\end{array}\right),\]
$\mathcal{A}'_2=g(\mathcal{A}_2+b\mathcal{A}_3)g^{-1}=$ 
\[\left( 
\begin{array}{ccc}
0&\alpha_5+b(\alpha_8+1)&1\\
\beta_4-\alpha_4+b(\beta_7-\alpha_7\beta_6)&\beta_5+b(\alpha_7+2\beta_8-\alpha_3)&\alpha_7+\beta_8-\alpha_3\\
*&*&*\end{array}\right),\]
\[\mathcal{A}'_3=g\mathcal{A}_3g^{-1}=\left( 
\begin{array}{ccc}
\alpha_7&\alpha_8&0\\
\beta_7&\beta_8&0\\
*&*&-\alpha_7-\beta_8\end{array}\right)\] which implies that:

1. If $\alpha_8\neq 1$ then one can make $\alpha'_5=0$ to get 
\[A_{3,2}=\left(
\begin{array}{ccccccccc}
\alpha_1 & \alpha_2& \alpha_3&0&0&1&\alpha_7&\alpha_8&0\\
\beta_1 & \beta_2& \beta_3&\beta_4&\beta_5&\alpha_7+\beta_8+\alpha_3&\beta_7&\beta_8&0\\
\gamma_1 & \gamma_2& \alpha_1+\beta_2&\gamma_4&\gamma_5&1+\beta_5&1+\alpha_1+\beta_4&\alpha_2+\beta_5&\alpha_7+\beta_8\end{array}\right),\] where $\alpha_8\neq 1$.

2. If $\alpha_8=1$ and $\alpha_7+\alpha_3\neq 0$ then one can make  $\beta'_5=0$ to get
\[A_{4,2}=\left(
\begin{array}{ccccccccc}
\alpha_1 & \alpha_2& \alpha_3&0&\alpha_5&1&\alpha_7&1&0\\
\beta_1 & \beta_2& \beta_3&\beta_4&0&\alpha_7+\beta_8+\alpha_3&\beta_7&\beta_8&0\\
\gamma_1 & \gamma_2& \alpha_1+\beta_2&\gamma_4&\gamma_5&1&1+\alpha_1+\beta_4&\alpha_2&\alpha_7+\beta_8\end{array}\right),\] where $\alpha_7\neq \alpha_3$.

3. If $\alpha_8=1$, $\alpha_7=\alpha_3$, that is $\beta_6=\beta_8$, and $\beta_7\neq \alpha_3\beta_8$ then one can make $\beta'_4=0$ to get
\[A_{5,2}=\left(
\begin{array}{ccccccccc}
\alpha_1 & \alpha_2& \alpha_3&0&\alpha_5&1&\alpha_3&1&0\\
\beta_1 & \beta_2& \beta_3&0&\beta_5&\beta_8&\beta_7&\beta_8&0\\
\gamma_1 & \gamma_2& \alpha_1+\beta_2&\gamma_4&\gamma_5&1+\beta_5&1+\alpha_1&\alpha_2+\beta_5&\alpha_3+\beta_8\end{array}\right),\]
where $\beta_7\neq \alpha_3\beta_8$.

4. If $\alpha_8=1$, $\alpha_7=\alpha_3$, $\beta_6=\beta_8$, $\beta_7=\alpha_3\beta_8$, $\alpha_3(\beta_3+\alpha_3\beta_8)\neq 0$ then   and one can make $\beta'_1=0$ to get

\[A_{6,2}=\left(
\begin{array}{ccccccccc}
\alpha_1 & \alpha_2& \alpha_3&0&\alpha_5&1&\alpha_3&1&0\\
0 & \beta_2& \beta_3&\beta_4&\beta_5&\beta_8&\alpha_3\beta_8&\beta_8&0\\
\gamma_1 & \gamma_2& \alpha_1+\beta_2&\gamma_4&\gamma_5&1+\beta_5&1+\alpha_1+\beta_4&\alpha_2+\beta_5&\alpha_3+\beta_8\end{array}\right),\] where 
$\alpha_3(\beta_3+\alpha_3\beta_8)\neq 0$.

5. If   $\alpha_8=1$, $\alpha_7=\alpha_3$, $\beta_6=\beta_8$, $\beta_7=\alpha_3\beta_8$, $\alpha_3(\beta_3+\alpha_3\beta_8)= 0$ and $\alpha_3(\alpha_3+1)\neq 0$ then, one can make $\alpha'_2=0$ to get
\[A_{7,2}=\left(
\begin{array}{ccccccccc}
\alpha_1 & 0& \alpha_3&0&\alpha_5&1&\alpha_3&1&0\\
\beta_1 & \beta_2& \alpha_3\beta_8&\beta_4&\beta_5&\beta_8&\alpha_3\beta_8&\beta_8&0\\
\gamma_1 & \gamma_2& \alpha_1+\beta_2&\gamma_4&\gamma_5&1+\beta_5&1+\alpha_1+\beta_4&\beta_5&\alpha_3+\beta_8\end{array}\right),\] where 
$\alpha_3\neq 0, 1$.

In the next two sub cases we deal with $\alpha_3(\alpha_3+1)= 0$ case.

\underline{\textbf{Subcase 3-1.\ $\alpha_8=1$, $\alpha_7=0$, $\beta_6=\beta_8$, $\beta_7=0$,  $\alpha_3= 0$.}}

1. If  $\beta_3\neq 0$ one can make $\gamma'_3=0$ to get 
\[A_{8,2}=\left(
\begin{array}{ccccccccc}
\alpha_1 & \alpha_2& 0&0&\alpha_5&1&0&1&0\\
\beta_1 & \alpha_1& \beta_3&\beta_4&\beta_5&\beta_8&0&\beta_8&0\\
\gamma_1 & \gamma_2& 0&\gamma_4&\gamma_5&1+\beta_5&1+\alpha_1+\beta_4&\alpha_2+\beta_5&\beta_8\end{array}\right),\] where 
$\beta_3\neq 0$.

2. If $\beta_3= 0$  and $\beta_2\neq \alpha_4\beta_8$ one can make $\gamma'_2=0$ to get 
\[A_{9,2}=\left(
\begin{array}{ccccccccc}
\alpha_1 & \alpha_2& 0&0&\alpha_5&1&0&1&0\\
\beta_1 & \beta_2& 0&\beta_4&\beta_5&\beta_8&0&\beta_8&0\\
\gamma_1 & 0& \alpha_1+\beta_2&\gamma_4&\gamma_5&1+\beta_5&1+\alpha_1+\beta_4&\alpha_2+\beta_5&\beta_8\end{array}\right),\] where 
$\beta_2\neq 0$.

3. If $\beta_3= 0$, $\beta_2= \alpha_4\beta_8$ and $\beta_1\neq 0$ one can make $\gamma'_1=0$ to get 
\[A_{10,2}=\left(
\begin{array}{ccccccccc}
\alpha_1 & \alpha_2& 0&0&\alpha_5&1&0&1&0\\
\beta_1 & 0& 0&\beta_4&\beta_5&\beta_8&0&\beta_8&0\\
0 & \gamma_2& \alpha_1&\gamma_4&\gamma_5&1+\beta_5&1+\alpha_1+\beta_4&\alpha_2+\beta_5&\beta_8\end{array}\right),\] where 
$\beta_1\neq 0$.

4. If  $\beta_3= 0$, $\beta_2= \alpha_4\beta_8$, $\beta_1=0$ and $\alpha_1\neq 1$ one can make $\gamma'_4=0$ to get 
\[A_{11,2}=\left(
\begin{array}{ccccccccc}
\alpha_1 & \alpha_2& 0&0&\alpha_5&1&0&1&0\\
0 & 0& 0&\beta_4&\beta_5&\beta_8&0&\beta_8&0\\
\gamma_1 & \gamma_2& \alpha_1&0&\gamma_5&11\beta_5&1+\alpha_1+\beta_4&\alpha_2+\beta_5&\beta_8\end{array}\right),\] where 
$\alpha_1\neq 1$.

5. If  $\beta_3= 0$, $\beta_2= \alpha_4\beta_8$, $\beta_1=0$, $\alpha_1= 1$ and $\beta_8\neq 0$ or $\beta_8=0$, $1+\alpha_2+\alpha_4+\beta_5\neq 0$ one can make $\gamma'_5=0$ to get 
\[A_{12,2}=\left(
\begin{array}{ccccccccc}
1 & \alpha_2& 0&0&\alpha_5&1&0&1&0\\
0 & 0& 0&\beta_4&\beta_5&\beta_8&0&\beta_8&0\\
\gamma_1 & \gamma_2& 1&\gamma_4&0&1+\beta_5&\beta_4&\alpha_2+\beta_5&0\end{array}\right),\] where $\beta_8\neq 0$ and if $\beta_8=0$ then $1+\alpha_2+\beta_5\neq 0.$

6. If  $\beta_3= 0$, $\beta_2= \alpha_4\beta_8$, $\beta_1=0$, $\alpha_1= 1$ and $\beta_8=0$, $1+\alpha_2+\alpha_4+\beta_5= 0$ one gets 
\[A_{13,2}=\left(
\begin{array}{ccccccccc}
1 & \alpha_2& 0&0&\alpha_5&1&0&1&0\\
0 & 0& 0&\beta_4&1+\alpha_2&0&0&0&0\\
\gamma_1 & \gamma_2& 1&\gamma_4&\gamma_5&\alpha_2&\beta_4&1&0\end{array}\right).\] 

In $\alpha_8=1$, $\alpha_7=1$, $\beta_6=\beta_7=\beta_8$, $\beta_3=\beta_8$,  $\alpha_3= 1$
case due to $a=\alpha_4+b$ one has

$\gamma'_1=-a(\alpha_1+a\alpha_7+a(\alpha_3+a\alpha_9))-b(\beta_1+a\beta_7+a(\beta_3+a\beta_9))+\gamma_1+ a(1-\alpha_1-\beta_4)+a(-\alpha_1-\beta_2-a(\alpha_7+\beta_8))=a(-\alpha_1-\beta_2-\beta_4+1)-a^2(1+\beta_8)-b\beta_1+\gamma_1=(\alpha_4+b)(\alpha_1+\beta_2+\beta_4-1)-(\alpha_4+b)^2(1+\beta_8)-b\beta_1+\gamma_1=(1+\beta_8)b^2+(\alpha_1+\beta_2+\beta_4+1+\beta_1)b+*$,

$\gamma'_2=-a(\alpha_2+a\alpha_8+a(\alpha_3+a\alpha_9))-b(\beta_2+a\beta_8+a(\beta_3+a\beta_9))+\gamma_2- a(\alpha_2+\beta_5)+a(-\alpha_1-\beta_2-a(\alpha_7+\beta_8))=
a(\alpha_1+\beta_2+\beta_5)+a^2(1+\beta_8)+b\beta_2=(\alpha_4+b)(\alpha_1+\beta_2+\beta_5)+(\alpha_4+b)^2(1+\beta_8)+b\beta_2=(1+\beta_8)b^2+(\alpha_1+\beta_1+\beta_2+\beta_5)b+*$,

$\gamma'_4=-a(\alpha_4+b\alpha_7+a(\alpha_6+b\alpha_9))-b(\beta_4+b\beta_7+a(\alpha_7+\beta_8-\alpha_3+b\beta_9))+\gamma_4+ b(1-\alpha_1-\beta_4)+a(1-\alpha_4-\beta_5-b(\alpha_7+\beta_8))=(\alpha_4+b)(\alpha_1+\beta_2+\beta_5)+(\alpha_4+b)^2(1+\beta_8)+b\beta_2=(1+\beta_8)b^2+(\alpha_1+\beta_5)b+*$,

$\gamma'_5=-a(\alpha_5+b\alpha_8+b(\alpha_6+b\alpha_9))-b(\beta_5+b\beta_8+b(\alpha_7+\beta_8-\alpha_3+b\beta_9))+\gamma_5- b(\alpha_2+\beta_5)+b(1-\alpha_4-\beta_5-b(\alpha_7+\beta_8))=a\alpha_5-b\beta_5+\gamma_5+
b(\alpha_2+1-\alpha_4)-b^2(1+\beta_8)=(1+\beta_8)b^2+(1+\alpha_2+\alpha_4+\alpha_5+\beta_5)b+*$, so one has the following case.

\underline{\textbf{Subcase 3-2.\ $\alpha_8=1$, $\alpha_7=1$, $\beta_6=\beta_7=\beta_8$, $\beta_3=\beta_8$,  $\alpha_3= 1$.}}

1. If  $\beta_8=1$ and $1+\alpha_1+\beta_1+\beta_2+\beta_4\neq 0 $ one can make $\gamma'_1=0$ to get 
\[A_{14,2}=\left(
\begin{array}{ccccccccc}
\alpha_1 & \alpha_2& 1&0&\alpha_5&1&1&1&0\\
\beta_1 & \beta_2& 1&\beta_4&\beta_5&1&1&1&0\\
0 & \gamma_2& \alpha_1+\beta_2&\gamma_4&\gamma_5&1+\beta_5&1+\alpha_1+\beta_4&\alpha_2+\beta_5&0\end{array}\right),\] 
where $1+\alpha_1+\beta_1+\beta_2+\beta_4\neq 0 $.

2. If $\beta_8=1$, $\beta_4=1+\alpha_1+\beta_1+\beta_2$ and $\alpha_1+\beta_1+\beta_2+\beta_5\neq 0$ one can make $\gamma'_2=0$ to get 
\[A_{15,2}=\left(
\begin{array}{ccccccccc}
\alpha_1 & \alpha_2& 1&0&\alpha_5&1&1&1&0\\
\beta_1 & \beta_2& 1&\beta_4&\beta_5&1&1&1&0\\
\gamma_1 & 0& \alpha_1+\beta_2&\gamma_4&\gamma_5&1+\beta_5&\beta_1+\beta_2&\alpha_2+\beta_5&0\end{array}\right),\] 
where $\beta_4=1+\alpha_1+\beta_1+\beta_2$, $\alpha_1+\beta_1+\beta_2+\beta_5=1+\beta_4+\beta_5\neq 0 $.

3. If $\beta_8=1$, $\beta_4=1+\alpha_1+\beta_1+\beta_2$, $\beta_5=1+\beta_4$ and $\alpha_1+\beta_5\neq 0$ one can make $\gamma'_4=0$ to get 
\[A_{16,2}=\left(
\begin{array}{ccccccccc}
\alpha_1 & \alpha_2& 1&0&\alpha_5&1&1&1&0\\
\beta_1 & \beta_2& 1&\beta_4&\beta_5&1&1&1&0\\
\gamma_1 & \gamma_2& \alpha_1+\beta_2&0&\gamma_5&1+\beta_5&\beta_1+\beta_2&\alpha_2+\beta_5&0\end{array}\right),\] 
where $\beta_5=\alpha_1+\beta_1+\beta_2$, $\beta_4=1+\alpha_1+\beta_1+\beta_2$, $\alpha_1+\beta_5\neq 0 $.

4. If $\beta_8=1$, $\beta_4=1+\alpha_1+\beta_1+\beta_2$, $\beta_5=\alpha_1+\beta_1+\beta_2$, $\alpha_1+\beta_5= 0$ and $1+\alpha_2+\alpha_4+\alpha_5+\beta_5\neq 0$ one can make $\gamma'_5=0$ to get 
\[A_{17,2}=\left(
\begin{array}{ccccccccc}
\alpha_1 & \alpha_2& 1&0&\alpha_5&1&1&1&0\\
\beta_1 & \beta_1& 1&1+\alpha_1&\alpha_1&1&1&1&0\\
\gamma_1 & \gamma_2& \alpha_1+\beta_1&\gamma_4&0&1+\alpha_1&0&\alpha_1+\alpha_2&0\end{array}\right),\] 
where  $\alpha_1+\alpha_2+\alpha_5\neq 1 $.

5. If $\beta_8=1$, $\beta_5=\alpha_1$, $\beta_2=\beta_1$, $\beta_4=1+\alpha_1$, and $\alpha_5=1+\alpha_1+\alpha_2+\alpha_4$ one gets 
\[A_{18,2}=\left(
\begin{array}{ccccccccc}
\alpha_1 & \alpha_2& 1&0&1+\alpha_1+\alpha_2&1&1&1&0\\
\beta_1 & \beta_1& 1&1+\alpha_1&\alpha_1&1&1&1&0\\
\gamma_1 & \gamma_2& \alpha_1+\beta_1&\gamma_4&\gamma_5&1+\alpha_1&0&\alpha_1+\alpha_2&0\end{array}\right).\] 
6. If  $\beta_8\neq1$, $\beta_4+\beta_5\neq 1 $ one can make $\gamma'_1=\gamma'_2$ to get 
\[A_{19,2}=\left(
\begin{array}{ccccccccc}
\alpha_1 & \alpha_2& 1&0&\alpha_5&1&1&1&0\\
\beta_1 & \beta_2& \beta_8&\beta_4&\beta_5&\beta_8&\beta_8&\beta_8&0\\
\gamma_1 & \gamma_1& \alpha_1+\beta_2&\gamma_4&\gamma_5&1+\beta_5&1+\alpha_1+\beta_4&\alpha_2+\beta_5&1+\beta_8\end{array}\right),\] 
where $\beta_8\neq 1$, $\beta_4+\beta_5\neq 1 $.

7. If  $\beta_8\neq1$, $\beta_5=1+\beta_4 $ and $\beta_2\neq \beta_1 $ one can make $\gamma'_1=\gamma'_4$ to get 
\[A_{20,2}=\left(
\begin{array}{ccccccccc}
\alpha_1 & \alpha_2& 1&0&\alpha_5&1&1&1&0\\
\beta_1 & \beta_2& \beta_8&\beta_4&1+\beta_4&\beta_8&\beta_8&\beta_8&0\\
\gamma_1 & \gamma_2& \alpha_1+\beta_2&\gamma_1&\gamma_5&\beta_4&1+\alpha_1+\beta_4&1+\alpha_2+\beta_4&1+\beta_8\end{array}\right),\] 
where $\beta_8\neq 1$, $\beta_2\neq \beta_1 $.

8. If  $\beta_8\neq1$, $\beta_5=1+\beta_4 $, $\beta_2= \beta_1 $ and $\alpha_1+\alpha_2+\alpha_4+\alpha_5\neq 1$ one can make $\gamma'_1=\gamma'_5$ to get 
\[A_{21,2}=\left(
\begin{array}{ccccccccc}
\alpha_1 & \alpha_2& 1&0&\alpha_5&1&1&1&0\\
\beta_1 & \beta_1& \beta_8&\beta_4&1+\beta_4&\beta_8&\beta_8&\beta_8&0\\
\gamma_1 & \gamma_2& \alpha_1+\beta_2&\gamma_4&\gamma_1&\beta_4&1+\alpha_1+\beta_4&1+\alpha_2+\beta_4&1+\beta_8\end{array}\right),\] 
where $\beta_8\neq 1$, $\alpha_1+\alpha_2+\alpha_5\neq 1$.

9. If  $\beta_8\neq1$, $\beta_5=1+\beta_4 $, $\beta_2= \beta_1 $ and $\alpha_5=1+\alpha_1+\alpha_2+\alpha_4$ one can make $\gamma'_1=0$ to get 
\[A_{22,2}=\left(
\begin{array}{ccccccccc}
\alpha_1 & \alpha_2& 1&0&1+\alpha_1+\alpha_2&1&1&1&0\\
\beta_1 & \beta_1& \beta_8&\beta_4&1+\beta_4&\beta_8&\beta_8&\beta_8&0\\
0 & \gamma_2& \alpha_1+\beta_2&\gamma_4&\gamma_5&\beta_4&1+\alpha_1+\beta_4&1+\alpha_2+\beta_4&1+\beta_8\end{array}\right),\] 
where $\beta_8\neq 1$.

\underline{\textbf{Case 4.\ $\alpha_9=\beta_9=\alpha_6= 0$, $\alpha_7+\beta_8+\alpha_3= 1$},} with respect to $g$ with the $c=1$ and if $a=\beta_4+b\beta_7$ then $\beta'_4=0$,  
$\mathcal{A}'_1=g(\mathcal{A}_1+b\mathcal{A}_3)g^{-1}=$ 
\[\left( 
\begin{array}{ccc}
\alpha_1-\beta_4(\alpha_3+\alpha_7)-b\beta_7(\alpha_3+\alpha_7)&\alpha_2-\beta_4(\alpha_3+\alpha_8)-b\beta_7(\alpha_3+\alpha_8)&\alpha_3\\
\beta_1-\beta_4(\beta_3+\beta_7)-b\beta_7(\beta_3+\beta_7)&\beta_2-\beta_4(\beta_3+\beta_8)-b\beta_7(\beta_3+\beta_8)&\beta_3\\
*&*&*\end{array}\right),\]
\[\mathcal{A}'_2=g(\mathcal{A}_2+b\mathcal{A}_3)g^{-1}= 
\left( 
\begin{array}{ccc}
\alpha_4+b\alpha_7&\alpha_5+b\alpha_8&0\\
0&\beta_5+b(1+\beta_8)&1\\
*&*&*\end{array}\right),\]
\[\mathcal{A}'_3=g\mathcal{A}_3g^{-1}=\left( 
\begin{array}{ccc}
\alpha_7&\alpha_8&0\\
\beta_7&\beta_8&0\\
*&*&-1-\alpha_3\end{array}\right).\] 
Therefore:\\
1. If $\beta_8\neq 1$  one can make  $\beta'_5=0$ to get
\[A_{23,2}=\left(
\begin{array}{ccccccccc}
\alpha_1 & \alpha_2& \alpha_3&\alpha_4&\alpha_5&0&1+\alpha_3+\beta_8&\alpha_8&0\\
\beta_1 & \beta_2& \beta_3&0&0&1&\beta_7&\beta_8&0\\
\gamma_1 & \gamma_2& \alpha_1+\beta_2&\gamma_4&\gamma_5&1+\alpha_4&1+\alpha_1&\alpha_2&\alpha_3+1\end{array}\right),\] where $\beta_8\neq 1$.

2. If $\beta_8=1$ and $\alpha_8\neq 0$ one can make  $\alpha'_5=0$ to get
\[A_{24,2}=\left(
\begin{array}{ccccccccc}
\alpha_1 & \alpha_2& \alpha_3&\alpha_4&0&0&2+\alpha_3&\alpha_8&0\\
\beta_1 & \beta_2& \beta_3&0&\beta_5&1&\beta_7&1&0\\
\gamma_1 & \gamma_2& \alpha_1+\beta_2&\gamma_4&\gamma_5&1+\alpha_4+\beta_5&1+\alpha_1&\alpha_2+\beta_5&\alpha_3+1\end{array}\right),\] where $\alpha_8\neq 0$.

3. If $\beta_8=1$, $\alpha_8= 0$ and $\alpha_7=\alpha_3\neq 0$ one can make  $\alpha'_4=0$ to get
\[A_{25,2}=\left(
\begin{array}{ccccccccc}
\alpha_1 & \alpha_2& \alpha_3&0&\alpha_5&0&\alpha_3&0&0\\
\beta_1 & \beta_2& \beta_3&0&\beta_5&1&\beta_7&1&0\\
\gamma_1 & \gamma_2& \alpha_1+\beta_2&\gamma_4&\gamma_5&1+\beta_5&1+\alpha_1&\alpha_2+\beta_5&\alpha_3+1\end{array}\right),\] where $\alpha_3\neq 0$.

4. If $\beta_8=1$, $\alpha_8= 0$, $\alpha_3=0$ and $\beta_7(\beta_3+1)\neq 0$ one can make  $\beta'_2=0$ to get
\[A_{26,2}=\left(
\begin{array}{ccccccccc}
\alpha_1 & \alpha_2& 0&\alpha_4&\alpha_5&0&0&0&0\\
\beta_1 & 0& \beta_3&0&\beta_5&1&\beta_7&1&0\\
\gamma_1 & \gamma_2& \alpha_1&\gamma_4&\gamma_5&1+\alpha_4+\beta_5&1+\alpha_1&\alpha_2+\beta_5&1\end{array}\right),\] where $\beta_7(\beta_3+1)\neq 0$.

If $\beta_8=1$, $\alpha_8= 0$, $\alpha_3=0$, $\beta_7= 0$ then  $\alpha_7= 0$ and

$\gamma'_1=-a(\alpha_1+a\alpha_7+a(\alpha_3+a\alpha_9))-b(\beta_1+a\beta_7+a(\beta_3+a\beta_9))+\gamma_1+ a(1-\alpha_1-\beta_4)+a(-\alpha_1-\beta_2-a(\alpha_7+\beta_8))=(\beta_1+\beta_3\beta_4)b+*$,

$\gamma'_2=-a(\alpha_2+a\alpha_8+a(\alpha_3+a\alpha_9))-b(\beta_2+a\beta_8+a(\beta_3+a\beta_9))+\gamma_2- a(\alpha_2+\beta_5)+a(-\alpha_1-\beta_2-a(\alpha_7+\beta_8))=(\beta_2+\beta_4+\beta_3\beta_4)b+*$,

$\gamma'_4=-a(\alpha_4+b\alpha_7+a(\alpha_6+b\alpha_9))-b(\beta_4+b\beta_7+a(\alpha_7+\beta_8-\alpha_3+b\beta_9))+\gamma_4+ b(1-\alpha_1-\beta_4)+a(1-\alpha_4-\beta_5-b(\alpha_7+\beta_8))=(1+\alpha_1)b+*$,

$\gamma'_5=-a(\alpha_5+b\alpha_8+b(\alpha_6+b\alpha_9))-b(\beta_5+b\beta_8+b(\alpha_7+\beta_8-\alpha_3+b\beta_9))+\gamma_5- b(\alpha_2+\beta_5)+b(1-\alpha_4-\beta_5-b(\alpha_7+\beta_8))=b^2+(1+\alpha_2+\alpha_4+\beta_5)b+*$. Therefore one has the following:

5. If $\beta_8=1$, $\alpha_8= 0$, $\alpha_3=0$, $\beta_7= 0$ and $\beta_1\neq \beta_3\beta_4$ one can make  $\gamma'_1=0$ to get
\[A_{27,2}=\left(
\begin{array}{ccccccccc}
\alpha_1 & \alpha_2& 0&\alpha_4&\alpha_5&0&0&0&0\\
\beta_1 & \beta_2& \beta_3&0&\beta_5&1&0&1&0\\
0 & \gamma_2& \alpha_1+\beta_2&\gamma_4&\gamma_5&1+\alpha_4+\beta_5&1+\alpha_1&\alpha_2+\beta_5&1\end{array}\right),\] where $\beta_1\neq 0$.

6. If $\beta_8=1$, $\alpha_8= 0$, $\alpha_3=0$, $\beta_7= 0$, $\beta_1= \beta_3\beta_4$ and $\beta_2\neq \beta_3\beta_4+\beta_4$ one can make  $\gamma'_2=0$ to get
\[A_{28,2}=\left(
\begin{array}{ccccccccc}
\alpha_1 & \alpha_2& 0&\alpha_4&\alpha_5&0&0&0&0\\
0 & \beta_2& \beta_3&0&\beta_5&1&0&1&0\\
\gamma_1 & 0& \alpha_1+\beta_2&\gamma_4&\gamma_5&1+\alpha_4+\beta_5&1+\alpha_1&\alpha_2+\beta_5&1\end{array}\right),\] where $\beta_2\neq 0$.

7. If $\beta_8=1$, $\alpha_8= 0$, $\alpha_3=0$, $\beta_7= 0$, $\beta_1= \beta_3\beta_4$,  $\beta_2= \beta_3\beta_4+\beta_4$ and $\alpha_1\neq 1$ one can make  $\gamma'_4=0$ to get
\[A_{29,2}=\left(
\begin{array}{ccccccccc}
\alpha_1 & \alpha_2& 0&\alpha_4&\alpha_5&0&0&0&0\\
0 & 0& \beta_3&0&\beta_5&1&0&1&0\\
\gamma_1 & \gamma_2& \alpha_1+\beta_2&0&\gamma_5&1+\alpha_4+\beta_5&1+\alpha_1&\alpha_2+\beta_5&1\end{array}\right),\] where $\alpha_1\neq 1$.

8. If $\beta_8=1$, $\alpha_8= 0$, $\alpha_3=0$, $\beta_7= 0$, $\beta_1= \beta_3\beta_4$,  $\beta_2= \beta_3\beta_4+\beta_4$ and $\alpha_1=1 \beta_4$ one can make  $\gamma'_5=0$ to get
\[A_{30,2}=\left(
\begin{array}{ccccccccc}
1 & \alpha_2& 0&\alpha_4&\alpha_5&0&0&0&0\\
0 & 0& \beta_3&0&\beta_5&1&0&1&0\\
\gamma_1 & \gamma_2& 1&\gamma_4&0&1+\alpha_4+\beta_5&0&\alpha_2+\beta_5&1\end{array}\right).\] 

If $\beta_8=1$, $\alpha_8= 0$, $\alpha_3=0$, $\beta_7\neq  0$, $\beta_3=1$ and $\beta_7\neq  1$ then  one can make $\beta'_1= 0$ to get

\[A_{31,2}=\left(
\begin{array}{ccccccccc}
\alpha_1 & \alpha_2& 0&\alpha_4&\alpha_5&0&0&0&0\\
0 & \beta_2& 1&0&\beta_5&1&\beta_7&1&0\\
\gamma_1 & \gamma_2& \alpha_1+\beta_2&\gamma_4&\gamma_5&1+\alpha_4+\beta_5&1+\alpha_1&\alpha_2+\beta_5&1\end{array}\right),\] where $\beta_7\neq 0,1$.

If $\beta_8=1$, $\alpha_8= 0$, $\alpha_3=0$, $\beta_3=1$ and $\beta_7=  1$ then

$\gamma'_1=-a(\alpha_1+a\alpha_7+a(\alpha_3+a\alpha_9))-b(\beta_1+a\beta_7+a(\beta_3+a\beta_9))+\gamma_1+ a(1-\alpha_1-\beta_4)+a(-\alpha_1-\beta_2-a(\alpha_7+\beta_8))=b^2+(1+\alpha_1+\beta_1+\beta_2+\beta_4)b+*$,

$\gamma'_2=-a(\alpha_2+a\alpha_8+a(\alpha_3+a\alpha_9))-b(\beta_2+a\beta_8+a(\beta_3+a\beta_9))+\gamma_2- a(\alpha_2+\beta_5)+a(-\alpha_1-\beta_2-a(\alpha_7+\beta_8))=b^2+(1+\alpha_1+\beta_5)b+*$,

$\gamma'_4=-a(\alpha_4+b\alpha_7+a(\alpha_6+b\alpha_9))-b(\beta_4+b\beta_7+a(\alpha_7+\beta_8-\alpha_3+b\beta_9))+\gamma_4+ b(1-\alpha_1-\beta_4)+a(1-\alpha_4-\beta_5-b(\alpha_7+\beta_8))=b^2+(\alpha_1+\alpha_4+\beta_4+\beta_5)b+*$,

$\gamma'_5=-a(\alpha_5+b\alpha_8+b(\alpha_6+b\alpha_9))-b(\beta_5+b\beta_8+b(\alpha_7+\beta_8-\alpha_3+b\beta_9))+\gamma_5- b(\alpha_2+\beta_5)+b(1-\alpha_4-\beta_5-b(\alpha_7+\beta_8))=b^2+(1+\alpha_2+\alpha_4+\alpha_5+\beta_5)b+*$. Therefore one has the following:

9. If $\beta_8=1$, $\alpha_8= 0$, $\alpha_3=0$, $\beta_3=1$, $\beta_7= 1$ and $\beta_5\neq \beta_1+\beta_2+\beta_4$ then  one can make $\gamma'_1=\gamma'_2$ to get

\[A_{32,2}=\left(
\begin{array}{ccccccccc}
\alpha_1 & \alpha_2& 0&\alpha_4&\alpha_5&0&0&0&0\\
\beta_1 & \beta_2& 1&0&\beta_5&1&1&1&0\\
\gamma_1 & \gamma_1& \alpha_1+\beta_2&\gamma_4&\gamma_5&1+\alpha_4+\beta_5&1+\alpha_1&\alpha_2+\beta_5&1\end{array}\right),\] where $\beta_5\neq \beta_1+\beta_2$.

10. If $\beta_8=1$, $\alpha_8= 0$, $\alpha_3=0$, $\beta_3=1$, $\beta_7= 1$, $\beta_5= \beta_1+\beta_2+\beta_4$ and $\beta_4\neq 1+\alpha_4$ then  one can make $\gamma'_1=\gamma'_4$ to get

\[A_{33,2}=\left(
\begin{array}{ccccccccc}
\alpha_1 & \alpha_2& 0&\alpha_4&\alpha_5&0&0&0&0\\
\beta_1 & \beta_2& 1&0&\beta_5&1&1&1&0\\
\gamma_1 & \gamma_2& \alpha_1+\beta_2&\gamma_1&\gamma_5&1+\alpha_4+\beta_5&1+\alpha_1&\alpha_2+\beta_5&1\end{array}\right),\] where $\beta_5= \beta_1+\beta_2$, $\alpha_4\neq 1$.

11. If $\beta_8=1$, $\alpha_8= 0$, $\alpha_3=0$, $\beta_3=1$, $\beta_7= 1$, $\beta_5= \beta_1+\beta_2+\beta_4$, $\beta_4= 1+\alpha_4$ and $\alpha_5\neq\alpha_1+\alpha_2+\alpha_4$ then  one can make $\gamma'_1=\gamma'_5$ to get

\[A_{34,2}=\left(
\begin{array}{ccccccccc}
\alpha_1 & \alpha_2& 0&1&\alpha_5&0&0&0&0\\
\beta_1 & \beta_2& 1&0&\beta_5&1&1&1&0\\
\gamma_1 & \gamma_2& \alpha_1+\beta_2&\gamma_4&\gamma_1&\beta_5&1+\alpha_1&\alpha_2+\beta_5&1\end{array}\right),\] where $\beta_5= \beta_1+\beta_2$, $\alpha_5\neq 1+\alpha_1+\alpha_2$.

12. If $\beta_8=1$, $\alpha_8= 0$, $\alpha_3=0$, $\beta_3=1$, $\beta_7= 1$, $\beta_5= \beta_1+\beta_2+\beta_4$, $\beta_4= 1+\alpha_4$, $\alpha_5=\alpha_1+\alpha_2+\alpha_4$ then  one can make $\gamma'_1=0$ to get

\[A_{35,2}=\left(
\begin{array}{ccccccccc}
\alpha_1 & \alpha_2& 0&1&\alpha_5&0&0&0&0\\
\beta_1 & \beta_2& 1&0&\beta_5&1&1&1&0\\
0 & \gamma_2& \alpha_1+\beta_2&\gamma_1&\gamma_5&\beta_5&1+\alpha_1&\alpha_2+\beta_5&1\end{array}\right),\] where $\alpha_5=1+\alpha_1+\alpha_2$, $\beta_5= \beta_1+\beta_2$.

If $\alpha_9=\beta_9=\alpha_6= 0$, $\alpha_7+\beta_8-\alpha_3=0$, $\alpha_3\neq 0$ one can make $\alpha'_3=1$ to come to the case:
$\alpha_9=\beta_9=\alpha_6= 0$, $\alpha_7+\beta_8-\alpha_3=0$, $\alpha_3= 1$ with respect to $g$ with the $c=1$. In this case  
$\mathcal{A}'_1=g(\mathcal{A}_1+a\mathcal{A}_3)g^{-1}=$ 
\[ \left(\begin{array}{ccc}
\alpha_1+a\beta_8&\alpha_2+a(\alpha_3+\alpha_8)&1\\
\beta_1+a(\beta_3+\beta_7)&\beta_2+a(\beta_3+\beta_8)&\beta_3\\
*&*&*\end{array}\right), \]
\[\mathcal{A}'_2=g(\mathcal{A}_2+b\mathcal{A}_3)g^{-1}= 
\left( 
\begin{array}{ccc}
\alpha_4+b\alpha_7&\alpha_5+b\alpha_8&0\\
\beta_4+b\beta_7&\beta_5+b\beta_8&0\\
*&*&1+\alpha_4+\beta_5+b\end{array}\right),\] 
\[\mathcal{A}'_3=g\mathcal{A}_3g^{-1}=\left( 
\begin{array}{ccc}
\alpha_7&\alpha_8&0\\
\beta_7&1-\alpha_7&0\\
*&*&1\end{array}\right).\] 

\underline{\textbf{Case 5.\ $\alpha_9=\beta_9=\alpha_6= 0$, $\alpha_7+\beta_8-\alpha_3=0$, $\alpha_3= 1$.}} One can make $\gamma'_6=0$.

1. If $\beta_8\neq 0$ one can make $\alpha'_1= 0$ to get
\[A_{36,2}=\left(
\begin{array}{ccccccccc}
0 & \alpha_2& 1&\alpha_4&\alpha_5&0&\alpha_7&\alpha_8&0\\
\beta_1 & \beta_2& \beta_3&\beta_4&1+\alpha_4&0&\beta_7&1+\alpha_7&0\\
\gamma_1 & \gamma_2& \beta_2&\gamma_4&\gamma_5&0&1+\beta_4&1+\alpha_2+\alpha_4&1\end{array}\right),\] where $\alpha_7\neq 1$.

2. If $\beta_8= 0$ and $\beta_3\neq 0$ one can make $\beta'_2= 0$ to get
\[A_{37,2}=\left(
\begin{array}{ccccccccc}
\alpha_1 & \alpha_2& 1&\alpha_4&\alpha_5&0&1&\alpha_8&0\\
\beta_1 & 0& \beta_3&\beta_4&1+\alpha_4&0&\beta_7&0&0\\
\gamma_1 & \gamma_2& \alpha_1&\gamma_4&\gamma_5&0&1+\alpha_1+\beta_4&1+\alpha_2+\alpha_4&1\end{array}\right),\] where $\beta_3\neq 0$.

3. If $\beta_8= 0$, $\beta_3= 0$ and $\beta_7\neq 0$ one can make $\beta'_1= 0$ to get
\[A_{38,2}=\left(
\begin{array}{ccccccccc}
\alpha_1 & \alpha_2& 1&\alpha_4&\alpha_5&0&1&\alpha_8&0\\
0 & \beta_2& 0&\beta_4&1+\alpha_4&0&\beta_7&0&0\\
\gamma_1 & \gamma_2& \alpha_1+\beta_2&\gamma_4&\gamma_5&0&1+\alpha_1+\beta_4&1+\alpha_2+\alpha_4&1\end{array}\right),\] where $\beta_7\neq 0$.

4. If $\beta_8= 0$, $\beta_3= 0$, $\beta_7= 0$ and $\alpha_8\neq 1$ one can make $\alpha'_2= 0$ to get
\[A_{39,2}=\left(
\begin{array}{ccccccccc}
\alpha_1 & 0& 1&\alpha_4&\alpha_5&0&1&\alpha_8&0\\
\beta_1 & \beta_2& 0&\beta_4&1+\alpha_4&0&0&0&0\\
\gamma_1 & \gamma_2& \alpha_1+\beta_2&\gamma_4&\gamma_5&0&1+\alpha_1+\beta_4&1+\alpha_4&1\end{array}\right),\] where $\alpha_8\neq 1$.

5. If $\beta_8= 0$, $\beta_3= 0$, $\beta_7= 0$, $\alpha_8= 1$ one can make $\gamma'_8= 0$ to get
\[A_{40,2}=\left(
\begin{array}{ccccccccc}
\alpha_1 &1+\alpha_4& 1&\alpha_4&\alpha_5&0&1&1&0\\
\beta_1 & \beta_2& 0&\beta_4&1+\alpha_4&0&0&0&0\\
\gamma_1 & \gamma_2& \alpha_1+\beta_2&\gamma_4&\gamma_5&0&1+\alpha_1+\beta_4&0&1\end{array}\right).\]

If $\alpha_9=\beta_9=\alpha_6= \alpha_7+\beta_8-\alpha_3=\alpha_3= 0, \beta_3\neq 0$ then one can make $\beta'_3=1$ by appropriate choose of $c$ to come to:

\underline{\textbf{Case 6.\ $\alpha_9=\beta_9=\alpha_6= \alpha_7+\beta_8-\alpha_3=\alpha_3= 0, \beta_3= 1$}} with respect to $g$ with $c=1$. If $b=-(\alpha_1+\beta_2)$ then $\gamma'_3=0$ and
\[\mathcal{A}'_1=g(\mathcal{A}_1+b\mathcal{A}_3)g^{-1}= 
\left(\begin{array}{ccc}
\alpha_1+a\alpha_7&\alpha_2+a\alpha_8&0\\
\beta_1+a(1+\beta_7)&\beta_2+a(1+\beta_8)&1\\
*&*&0\end{array}\right),\]
\[\mathcal{A}'_2=g(\mathcal{A}_2+b\mathcal{A}_3)g^{-1}= 
\left( 
\begin{array}{ccc}
\alpha_4+b\alpha_7&\alpha_5+b\alpha_8&0\\
\beta_4+b\beta_7&\beta_5+b\beta_8&0\\
*&*&*\end{array}\right),\]
\[\mathcal{A}'_3=g\mathcal{A}_3g^{-1}=\left( 
\begin{array}{ccc}
-\beta_8&\alpha_8&0\\
\beta_7&\beta_8&0\\
*&*&0\end{array}\right).\] 

Therefore if $\alpha_8\neq 0$ one can make $\alpha'_2= 0$ to get
\[A_{41,2}=\left(
\begin{array}{ccccccccc}
\alpha_1 & 0& 0&\alpha_4&\alpha_5&0&\beta_8&\alpha_8&0\\
\beta_1 & \alpha_1& 1&\beta_4&\beta_5&0&\beta_7&\beta_8&0\\
\gamma_1 & \gamma_2& 0&\gamma_4&\gamma_5&1+\alpha_4+\beta_5&1+\alpha_1+\beta_4&\beta_5&0\end{array}\right),\] where $\alpha_8\neq 0$, if $\alpha_8= 0$, $\alpha_7\neq 0$ one can make $\alpha'_1= 0$ to get
\[A_{42,2}=\left(
\begin{array}{ccccccccc}
0 & \alpha_2& 0&\alpha_4&\alpha_5&0&\beta_8&0&0\\
\beta_1 & 0& 1&\beta_4&\beta_5&0&\beta_7&\beta_8&0\\
\gamma_1 & \gamma_2& 0&\gamma_4&\gamma_5&1+\alpha_4+\beta_5&1+\beta_4&\alpha_2+\beta_5&0\end{array}\right),\] where $\alpha_7=\beta_8\neq 0$ and 
if $\alpha_8=\alpha_7=\beta_8= 0$ then on can make $\beta'_2=0$ to get
\[A_{43,2}=\left(
\begin{array}{ccccccccc}
\alpha_1 & \alpha_2& 0&\alpha_4&\alpha_5&0&0&0&0\\
\beta_1 & 0& 1&\beta_4&\beta_5&0&\beta_7&0&0\\
\gamma_1 & \gamma_2& \alpha_1&\gamma_4&\gamma_5&1+\alpha_4+\beta_5&1+\alpha_1+\beta_4&\alpha_2=\beta_5&0\end{array}\right).\]

\underline{\textbf{Case 7.\ $\alpha_9=\beta_9=\alpha_6= \alpha_7+\beta_8-\alpha_3=\alpha_3= \beta_3= 0$}.} In this case
$\alpha'_9=\beta'_9=\alpha_6'= \beta'_6=\alpha'_3=\beta'_3=0$ and \[A'_1=g(A_1+aA_3)g^{-1}= 
\left(\begin{array}{ccc}
\alpha_1+a\alpha_7&\alpha_2+a\alpha_8&0\\
\beta_1+a\beta_7&\beta_2+a\beta_8)&0\\
*&*&*\end{array}\right), \]
\[A'_2=g(A_2+aA_3)g^{-1}= 
\left(\begin{array}{ccc}
\alpha_4+b\alpha_7&\alpha_5+b\alpha_8&0\\
\beta_4+b\beta_7&\beta_5+b\beta_8)&0\\
*&*&*\end{array}\right), \] 
\[A'_3=cgA_3g^{-1}= 
\left(\begin{array}{ccc}
c\alpha_7&c\alpha_8&0\\
c\beta_7&c\beta_8&0\\
*&*&0\end{array}\right). \]
Therefore if  $\alpha_8\neq 0$ one can make $\alpha'_2=\alpha'_5=0, \alpha'_8=1$ to get
\[A_{44,2}=\left(
\begin{array}{ccccccccc}
\alpha_1 & 0& 0&\alpha_4&0&0&\alpha_7&1&0\\
\beta_1 & \beta_2& 0&\beta_4&\beta_5&0&\beta_7&\alpha_7&0\\
\gamma_1 & \gamma_2& \alpha_1+\beta_2&\gamma_4&\gamma_5&1+\alpha_4+\beta_5&1+\alpha_1+\beta_4&\beta_5&0\end{array}\right),\] if $\alpha_8= 0$, $\alpha_7\neq 0$ one can make $\alpha'_1=\alpha'_4=0, \alpha'_7=1$ to get

\[A_{45,2}=\left(
\begin{array}{ccccccccc}
0 & \alpha_2& 0&0&\alpha_5&0&1&0&0\\
\beta_1 & \beta_2& 0&\beta_4&\beta_5&0&\beta_7&1&0\\
\gamma_1 & \gamma_2& \beta_2&\gamma_4&\gamma_5&1+\beta_5&1+\beta_4&\alpha_2+\beta_5&0\end{array}\right),\] if  $\alpha_7=\alpha_8= 0$ then $\beta_8=0$ and if $\beta_7\neq 0$ one can make $\beta'_1=\beta'_4=0, \beta'_7=1$ to get
\[A_{46,2}=\left(
\begin{array}{ccccccccc}
\alpha_1 & \alpha_2& 0&\alpha_4&\alpha_5&0&0&0&0\\
0 & \beta_2& 0&0&\beta_5&0&1&0&0\\
\gamma_1 & \gamma_2& \alpha_1+\beta_2&\gamma_4&\gamma_5&1+\alpha_4+\beta_5&1+\alpha_1&\alpha_2+\beta_5&0\end{array}\right).\] 

%\underline{\textbf{Case 8.\
In $\alpha_9=\alpha_8=\alpha_7=\alpha_6=\alpha_3=\beta_9=\beta_8=\beta_7=\beta_6=\beta_3=0$ case
\[A'_1=g(A_1+aA_3)g^{-1}= 
\left(\begin{array}{ccc}
\alpha_1&\alpha_2&0\\
\beta_1&\beta_2&0\\
*&*&-\alpha_1-\beta_2\end{array}\right), \]
\[A'_2=g(A_2+aA_3)g^{-1}= 
\left(\begin{array}{ccc}
\alpha_4&\alpha_5&0\\
\beta_4&\beta_5&0\\
*&*&1-\alpha_4-\beta_5\end{array}\right), \] 
\[A'_3=cgA_3g^{-1}= 
\left(\begin{array}{ccc}
0&0&0\\
0&0&0\\
1-\alpha_1-\beta_4&-\alpha_2-\beta_5&0\end{array}\right)\] and \\
$\gamma'_1=c^{-1}[-a(3\alpha_1+\beta_2+\beta_4-1)-b\beta_1+\gamma_1]$,\\
$\gamma'_2=c^{-1}[-a(2\alpha_2+\alpha_1+\beta_5+\beta_2)-b\beta_2+\gamma_2]$,\\
$\gamma'_4=c^{-1}[-a(2\alpha_4+\beta_5-1)-b(\alpha_1+2\beta_4-1)+\gamma_4]$,\\
$\gamma'_5=c^{-1}[-a\alpha_5-b(\alpha_2+\alpha_4+3\beta_5-1)+\gamma_5]$.

Consider row vectors \[ \begin{array}{ll} 1.\  (3\alpha_1+\beta_2+\beta_4-1,\beta_1), & 1'.\  (3\alpha_1+\beta_2+\beta_4-1,\beta_1, \gamma_1),\\
2.\ (2\alpha_2+\alpha_1+\beta_5+\beta_2, \beta_2),& 2'.\  (2\alpha_2+\alpha_1+\beta_5+\beta_2, \beta_2, \gamma_2),\\
3.\ (2\alpha_4+\beta_5-1, \alpha_1+2\beta_4-1),&  3'.\
(2\alpha_4+\beta_5-1, \alpha_1+2\beta_4-1, \gamma_4),\\
4.\ (\alpha_5, \alpha_2+\alpha_4+3\beta_5-1),& 4'.\
(\alpha_5, \alpha_2+\alpha_4+3\beta_5-1, \gamma_5),\end{array}\]
and let
$M_{(i,j,k,l)}$ (respectively, $M'_{(i,j,k,l)}$), where $i,j,k,l=1,2,3,4$, stand for matrix with the above defined rows $i.,j.,k.,l.$ (respectively, $i'.,j'.,k'.,l'.$), $\lambda\times i+ \mu \times j$, with $\lambda, \mu \in F$, stand for a linear combination of the rows $i.,j$, $rk(A)$ stand for the rank of a matrix $A$.

\underline{\textbf{Case 8. $rk(M_{(1,2,3,4)})=rk(M'_{(1,2,3,4)})$}.} In this case one can make $\gamma'_1=\gamma'_2=\gamma'_4=\gamma'_5=0$ to get
\[A_{47,2}=\left(
\begin{array}{ccccccccc}
\alpha_1 & \alpha_2& 0&\alpha_4&\alpha_5&0&0&0&0\\
\beta_1 & \beta_2& 0&\beta_4&\beta_5&0&0&0&0\\
0 & 0& \alpha_1+\beta_2&0&0&1+\alpha_4+\beta_5&1+\alpha_1+\beta_4&\alpha_2+\beta_5&0\end{array}\right).\] 

\underline{\textbf{Case 9. $rk(M_{(1,2,3,4)})\neq rk(M'_{(1,2,3,4)})$}.}

1. If  $rk(M_{(2,3,4)})=rk(M'_{(2,3,4)})$ one can make $\gamma'_2=\gamma'_4=\gamma'_5=0$ and $\gamma'_1=1$ to get
\[A_{48,2}=\left(
\begin{array}{ccccccccc}
\alpha_1 & \alpha_2& 0&\alpha_4&\alpha_5&0&0&0&0\\
\beta_1 & \beta_2& 0&\beta_4&\beta_5&0&0&0&0\\
1 & 0& \alpha_1+\beta_2&0&0&1+\alpha_4+\beta_5&1+\alpha_1+\beta_4&\alpha_2+\beta_5&0\end{array}\right),\] where $rk(M_{(1,2,3,4)})=rk(M_{(2,3,4)})$.

2. If $rk(M_{(2,3,4)})<rk(M'_{(2,3,4)})$, $rk(M_{(1,3,4)})=rk(M'_{(1,3,4)})$ one can make $\gamma'_1=\gamma'_4=\gamma'_5=0$ and $\gamma'_2=1$ to get
\[A_{49,2}=\left(
\begin{array}{ccccccccc}
\alpha_1 & \alpha_2& 0&\alpha_4&\alpha_5&0&0&0&0\\
\beta_1 & \beta_2& 0&\beta_4&\beta_5&0&0&0&0\\
0 & 1& \alpha_1+\beta_2&0&0&1+\alpha_4+\beta_5&1+\alpha_1+\beta_4&\alpha_2+\beta_5&0\end{array}\right),\] where $rk(M_{(1,2,3,4)})=rk(M_{(1,3,4)})$, $rk(M_{(3,4)})=rk(M_{(2,3,4)})$.

3. If $rk(M_{(2,3,4)})<rk(M'_{(2,3,4)})$, $rk(M_{(1,3,4)})<rk(M'_{(1,3,4)})$, $rk(M_{(1,2,4)})=rk(M'_{(1,2,4)})$ one can make $\gamma'_1=\gamma'_2=\gamma'_5=0$ and $\gamma'_4=1$ to get
\[A_{50,2}=\left(
\begin{array}{ccccccccc}
\alpha_1 & \alpha_2& 0&\alpha_4&\alpha_5&0&0&0&0\\
\beta_1 & \beta_2& 0&\beta_4&\beta_5&0&0&0&0\\
0 & 0& \alpha_1+\beta_2&1&0&1+\alpha_4+\beta_5&1+\alpha_1+\beta_4&\alpha_2+\beta_5&0\end{array}\right),\] where $rk(M_{(1,2,3,4)})=rk(M_{(1,2,4)})$, $rk(M_{(2,4)})=rk(M_{(2,3,4)})$, $rk(M_{(1,4)})=rk(M_{(1,3,4)})$.

4. If $rk(M_{(2,3,4)})<rk(M'_{(2,3,4)})$, $rk(M_{(1,3,4)})<rk(M'_{(1,3,4)})$, $rk(M_{(1,2,4)})<rk(M'_{(1,2,4)})$, $rk(M_{(1,2,3)})=rk(M'_{(1,2,3)})$ . In this case one can make $\gamma'_1=\gamma'_2=\gamma'_4=0$ and $\gamma'_5=1$ to get
\[A_{51,2}=\left(
\begin{array}{ccccccccc}
\alpha_1 & \alpha_2& 0&\alpha_4&\alpha_5&0&0&0&0\\
\beta_1 & \beta_2& 0&\beta_4&\beta_5&0&0&0&0\\
0 & 0& \alpha_1+\beta_2&0&1&1+\alpha_4+\beta_5&1+\alpha_1+\beta_4&\alpha_2+\beta_5&0\end{array}\right),\] where $rk(M_{(1,2,3,4)})=rk(M_{(1,2,3)})$, $rk(M_{(2,3)})=rk(M_{(2,3,4)})$, $rk(M_{(1,3)})=rk(M_{(1,3,4)})$, $rk(M_{(1,2)})=rk(M_{(1,2,4)})$.

\underline{\textbf{Subcase 9-1. $rk(M_{(2,3,4)})<rk(M'_{(2,3,4)})$, $rk(M_{(1,3,4)})<rk(M'_{(1,3,4)})$},}\\ \underline{$rk(M_{(1,2,4)})<rk(M'_{(1,2,4)})$, $rk(M_{(1,2,3)})<rk(M'_{(1,2,3)})$.}

1. If $rk(M_{(3,4)})=rk(M'_{(3,4)})$ one can make $\gamma'_4=\gamma'_5=0$, $\gamma'_1=1$ and $\gamma'_2\neq 0$ to get
\[A_{52,2}=\left(
\begin{array}{ccccccccc}
\alpha_1 & \alpha_2& 0&\alpha_4&\alpha_5&0&0&0&0\\
\beta_1 & \beta_2& 0&\beta_4&\beta_5&0&0&0&0\\
1 & \gamma_2& \alpha_1+\beta_2&0&0&1+\alpha_4+\beta_5&1+\alpha_1+\beta_4&\alpha_2+\beta_5&0\end{array}\right),\] where $rk(M_{(1,2,3,4)})=rk(M_{(2-\gamma_2\times1,3,4)})$, $rk(M_{(3,4)})=rk(M_{(2,3,4)})$, $rk(M_{(3,4)})=rk(M_{(1,3,4)})$, $rk(M_{(2-\gamma_2\times 1,4)})=rk(M_{(1,2,4)})$, 
$rk(M_{(2-\gamma_2\times 1,3)})=rk(M_{(1,2,3)})$ and $\gamma_2\neq 0$.

2. If $rk(M_{(3,4)})<rk(M'_{(3,4)})$, $rk(M_{(2,4)})=rk(M'_{(2,4)})$ one can make $\gamma'_2=\gamma'_5=0$, $\gamma'_1=1$ and $\gamma'_4\neq 0$ to get
\[A_{53,2}=\left(
\begin{array}{ccccccccc}
\alpha_1 & \alpha_2& 0&\alpha_4&\alpha_5&0&0&0&0\\
\beta_1 & \beta_2& 0&\beta_4&\beta_5&0&0&0&0\\
1 & 0& \alpha_1+\beta_2&\gamma_4&0&1+\alpha_4+\beta_5&1+\alpha_1+\beta_4&\alpha_2+\beta_5&0\end{array}\right),\] where $rk(M_{(1,2,3,4)})=rk(M_{(2,3-\gamma_4\times 1,4)})$, $rk(M_{(2,4)})=rk(M_{(2,3,4)})$, $rk(M_{(3-\gamma_4\times 1,4)})=rk(M_{(1,3,4)})$, $rk(M_{(2,4)})=rk(M_{(1,2,4)})$, 
$rk(M_{(2,3-\gamma_4\times 1)})=rk(M_{(1,2,3)})$, $rk(M_{(4)})=rk(M_{(3,4)})$ and $\gamma_4\neq 0$.

3. If $rk(M_{(3,4)})<rk(M'_{(3,4)})$, $rk(M_{(2,4)})<rk(M'_{(2,4)})$, $rk(M_{(1,4)})=rk(M'_{(1,4)})$ one can make $\gamma'_1=\gamma'_5=0$, $\gamma'_2=1$ and $\gamma'_4\neq 0$ to get
\[A_{54,2}=\left(
\begin{array}{ccccccccc}
\alpha_1 & \alpha_2& 0&\alpha_4&\alpha_5&0&0&0&0\\
\beta_1 & \beta_2& 0&\beta_4&\beta_5&0&0&0&0\\
0 & 1& \alpha_1+\beta_2&\gamma_4&0&1+\alpha_4+\beta_5&1+\alpha_1+\beta_4&\alpha_2+\beta_5&0\end{array}\right),\] where $rk(M_{(1,2,3,4)})=rk(M_{(1,3-\gamma_4\times 2,4)})$, $rk(M_{(3-\gamma_4\times 2,4)})=rk(M_{(2,3,4)})$, $rk(M_{(1,4)})=rk(M_{(1,3,4)})$, $rk(M_{(1,4)})=rk(M_{(1,2,4)})$, 
$rk(M_{(1,3-\gamma_4\times 2)})=rk(M_{(1,2,3)})$, $rk(M_{(4)})=rk(M_{(3,4)})$,  $rk(M_{(4)})=rk(M_{(2,4)})$ and $\gamma_4\neq 0$.

4. If $rk(M_{(3,4)})<rk(M'_{(3,4)})$, $rk(M_{(2,4)})<rk(M'_{(2,4)})$, $rk(M_{(1,4)})<rk(M'_{(1,4)})$ , $rk(M_{(2,3)})=rk(M'_{(2,3)})$ one can make $\gamma'_2=\gamma'_4=0$, $\gamma'_1=1$ and $\gamma'_5\neq 0$ to get
\[A_{55,2}=\left(
\begin{array}{ccccccccc}
\alpha_1 & \alpha_2& 0&\alpha_4&\alpha_5&0&0&0&0\\
\beta_1 & \beta_2& 0&\beta_4&\beta_5&0&0&0&0\\
1 & 0& \alpha_1+\beta_2&0&\gamma_5&1+\alpha_4+\beta_5&1+\alpha_1+\beta_4&\alpha_2+\beta_5&0\end{array}\right),\] where $rk(M_{(1,2,3,4)})=rk(M_{(2,3,4-\gamma_5\times 1)})$, $rk(M_{(2,3)})=rk(M_{(2,3,4)})$, $rk(M_{(3,4-\gamma_5\times 1)})=rk(M_{(1,3,4})$, $rk(M_{(2,4-\gamma_5\times 1)})=rk(M_{(1,2,4)})$, 
$rk(M_{(2,3)})=rk(M_{(1,2,3)})$, $rk(M_{(3)})=rk(M_{(3,4)})$, $rk(M_{(2)})=rk(M_{(2,4)})$, $rk(M_{(4-\gamma_5\times 1)})=rk(M_{(1,4)})$ and $\gamma_5\neq 0$.

5. If $rk(M_{(3,4)})<rk(M'_{(3,4)})$, $rk(M_{(2,4)})<rk(M'_{(2,4)})$, $rk(M_{(1,4)})<rk(M'_{(1,4)})$ , $rk(M_{(2,3)})<rk(M'_{(2,3)})$, $rk(M_{(1,3)})=rk(M'_{(1,3)})$ one can make $\gamma'_1=\gamma'_4=0$, $\gamma'_2=1$ and $\gamma'_5\neq 0$ to get
\[A_{56,2}=\left(
\begin{array}{ccccccccc}
\alpha_1 & \alpha_2& 0&\alpha_4&\alpha_5&0&0&0&0\\
\beta_1 & \beta_2& 0&\beta_4&\beta_5&0&0&0&0\\
0 & 1& \alpha_1+\beta_2&0&\gamma_5&1+\alpha_4+\beta_5&1+\alpha_1+\beta_4&\alpha_2+\beta_5&0\end{array}\right),\] where $rk(M_{(1,2,3,4)})=rk(M_{(1,3,4-\gamma_5\times 2)})$, $rk(M_{(3,4-\gamma_5\times 2)})=rk(M_{(2,3,4)})$, $rk(M_{(1,3)})=rk(M_{(1,3,4})$, $rk(M_{(1,4-\gamma_5\times 2)})=rk(M_{(1,2,4)})$, 
$rk(M_{(1,3)})=rk(M_{(1,2,3)})$, $rk(M_{(3)})=rk(M_{(3,4)})$, $rk(M_{(4-\gamma_5\times 2)})=rk(M_{(2,4)})$, $rk(M_{(1)})=rk(M_{(1,4)})$, $rk(M_{(3)})=rk(M_{(2,3)})$ and $\gamma_5\neq 0$.

6. If $rk(M_{(3,4)})<rk(M'_{(3,4)})$, $rk(M_{(2,4)})<rk(M'_{(2,4)})$, $rk(M_{(1,4)})<rk(M'_{(1,4)})$ , $rk(M_{(2,3)})<rk(M'_{(2,3)})$, $rk(M_{(1,3)})<rk(M'_{(1,3)})$, $rk(M_{(1,2)})=rk(M'_{(1,2)})$
one can make $\gamma'_1=\gamma'_2=0$, $\gamma'_4=1$ and $\gamma'_5\neq 0$ to get
\[A_{57,2}=\left(
\begin{array}{ccccccccc}
\alpha_1 & \alpha_2& 0&\alpha_4&\alpha_5&0&0&0&0\\
\beta_1 & \beta_2& 0&\beta_4&\beta_5&0&0&0&0\\
0 & 0& \alpha_1+\beta_2&1&\gamma_5&1+\alpha_4+\beta_5&1+\alpha_1+\beta_4&\alpha_2+\beta_5&0\end{array}\right),\] where $rk(M_{(1,2,3,4)})=rk(M_{(1,2,4-\gamma_5\times 3)})$, $rk(M_{(2,4-\gamma_5\times 3)})=rk(M_{(2,3,4)})$, $rk(M_{(1,4-\gamma_5\times 3)})=rk(M_{(1,3,4})$, $rk(M_{(1,2)})=rk(M_{(1,2,4)})$, 
$rk(M_{(1,2)})=rk(M_{(1,2,3)})$, $rk(M_{(4-\gamma_5\times 3)})=rk(M_{(3,4)})$, $rk(M_{(2)})=rk(M_{(2,4)})$, $rk(M_{(1)})=rk(M_{(1,4)})$, $rk(M_{(1)})=rk(M_{(1,3)})$ and $\gamma_5\neq 0$.

\underline{\textbf{Subsubcase 9-1-1. $rk(M_{(2,3,4)})<rk(M'_{(2,3,4)})$, $rk(M_{(1,3,4)})<rk(M'_{(1,3,4)})$},}\\ \underline{$rk(M_{(1,2,4)})<rk(M'_{(1,2,4)})$, $rk(M_{(1,2,3)})<rk(M'_{(1,2,3)})$,  $rk(M_{(3,4)})<rk(M'_{(3,4)})$,}\\ \underline{$rk(M_{(2,4)})<rk(M'_{(2,4)})$, $rk(M_{(1,4)})<rk(M'_{(1,4)})$ , $rk(M_{(2,3)})<rk(M'_{(2,3)})$,}\\ \underline{$rk(M_{(1,3)})<rk(M'_{(1,3)})$, $rk(M_{(1,2)})<rk(M'_{(1,2)})$.}

1. If $rk(M_{(1)})=rk(M'_{(1)})$ one can make $\gamma'_1=0$, $\gamma'_2=1$, $\gamma'_4\neq 0$, $\gamma'_5\neq 0$ to get
\[A_{58,2}=\left(
\begin{array}{ccccccccc}
\alpha_1 & \alpha_2& 0&\alpha_4&\alpha_5&0&0&0&0\\
\beta_1 & \beta_2& 0&\beta_4&\beta_5&0&0&0&0\\
0 & 1& \alpha_1+\beta_2&\gamma_4&\gamma_5&1+\alpha_4+\beta_5&1+\alpha_1+\beta_4&\alpha_2+\beta_5&0\end{array}\right),\] where $rk(M_{(1,2,3,4)})=rk(M_{(1,3-\gamma_4\times 2,4-\gamma_5\times 2)})$, $rk(M_{(3-\gamma_4\times 2, 4-\gamma_5\times 2)})=rk(M_{(2,3,4)})$, $rk(M_{(1,\gamma_4\times 4-\gamma_5\times 3)})=rk(M_{(1,3,4})$, $rk(M_{(1,4-\gamma_5\times 2)})=rk(M_{(1,2,4)})$, 
$rk(M_{(1,3-\gamma_4\times 2)})=rk(M_{(1,2,3)})$, $rk(M_{(\gamma_4\times 4-\gamma_5\times 3)})=rk(M_{(3,4)})$, $rk(M_{(4-\gamma_5\times 2)})=rk(M_{(2,4)})$, $rk(M_{(1)})=rk(M_{(1,4)})$, $rk(M_{(1)})=rk(M_{(1,3)})$, $rk(M_{(1)})=rk(M_{(1,2)})$ and $\gamma_4\neq 0$, $\gamma_5\neq 0$.

2. If $rk(M_{(1)})<rk(M'_{(1)})$, $rk(M_{(2)})=rk(M'_{(2)})$ one can make $\gamma'_1=1$, $\gamma'_2=0$, $\gamma'_4\neq 0$, $\gamma'_5\neq 0$ to get
\[A_{59,2}=\left(
\begin{array}{ccccccccc}
\alpha_1 & \alpha_2& 0&\alpha_4&\alpha_5&0&0&0&0\\
\beta_1 & \beta_2& 0&\beta_4&\beta_5&0&0&0&0\\
1 & 0& \alpha_1+\beta_2&\gamma_4&\gamma_5&1+\alpha_4+\beta_5&1+\alpha_1+\beta_4&\alpha_2+\beta_5&0\end{array}\right),\] where $rk(M_{(1,2,3,4)})=rk(M_{(2,3-\gamma_4\times 1,4-\gamma_5\times 1)})$, $rk(M_{(2,\gamma_4\times 3-\gamma_5\times 3)})=rk(M_{(2,3,4)})$, $rk(M_{(3-\gamma_4\times 1,4-\gamma_5\times 1)})=rk(M_{(1,3,4})$, $rk(M_{(2,4-\gamma_5\times 1)})=rk(M_{(1,2,4)})$, 
$rk(M_{(2,3-\gamma_4\times 1)})=rk(M_{(1,2,3)})$, $rk(M_{(\gamma_4\times 4-\gamma_5\times 3)})=rk(M_{(3,4)})$, $rk(M_{(2)})=rk(M_{(2,4)})$, $rk(M_{(4-\gamma_4\times 1)})=rk(M_{(1,4)})$, $rk(M_{(3-\gamma_5\times 1)})=rk(M_{(1,3)})$, $rk(M_{(2)})=rk(M_{(1,2)})$, $rk(M_{(1)})=0$ and $\gamma_4\neq 0$, $\gamma_5\neq 0$.

3. If $rk(M_{(1)})<rk(M'_{(1)})$, $rk(M_{(2)})<rk(M'_{(2)})$, $rk(M_{(3)})=rk(M'_{(3)})$ one can make $\gamma'_1=1$, $\gamma'_2\neq 0$, $\gamma'_4= 0$, $\gamma'_5\neq 0$ to get MSC
\[\left(
\begin{array}{ccccccccc}
\alpha_1 & \alpha_2& 0&\alpha_4&\alpha_5&0&0&0&0\\
\beta_1 & \beta_2& 0&\beta_4&\beta_5&0&0&0&0\\
1 & \gamma_2& -\alpha_1-\beta_2&0&\gamma_5&1-\alpha_4-\beta_5&1-\alpha_1-\beta_4&-\alpha_2-\beta_5&0\end{array}\right),\] where $rk(M_{(1,2,3,4)})=rk(M_{(2-\gamma_2\times 1,3,4-\gamma_5\times 1)})$, $rk(M_{(3,\gamma_2\times 4-\gamma_5\times 2)})=rk(M_{(2,3,4)})$, $rk(M_{(3,4-\gamma_5\times 1)})=rk(M_{(1,3,4})$, $rk(M_{(2-\gamma_2\times 1,4-\gamma_5\times 1)})=rk(M_{(1,2,4)})$, 
$rk(M_{(2-\gamma_2\times 1),3})=rk(M_{(1,2,3)})$, $rk(M_{(3)})=rk(M_{(3,4)})$, $rk(M_{(\gamma_2\times 4-\gamma_5\times 2)})=rk(M_{(2,4)})$, $rk(M_{(4-\gamma_5\times 1)})=rk(M_{(1,4)})$, $rk(M_{(3)})=rk(M_{(1,3)})$, $rk(M_{(1)})=0$, $rk(M_{(2)})=0$ and $\gamma_2\neq 0$, $\gamma_5\neq 0$.

The system $rk(M_{(1)})=0$, $rk(M_{(2)})=0$ implies that $\beta_1=\beta_2=0, \beta_4=1-3\alpha_1, \beta_5=-\alpha_1-2\alpha_2$ and one can make $\gamma'_1= 1$ to get \[A_{60,2}=\left(
\begin{array}{ccccccccc}
\alpha_1 & \alpha_2& 0&\alpha_4&\alpha_5&0&0&0&0\\
0 & 0& 0&1+\alpha_1&\alpha_1&0&0&0&0\\
1 & \gamma_2& \alpha_1&0&\gamma_5&1+\alpha_1+\alpha_4&0&\alpha_1+\alpha_2&0\end{array}\right),\] where $\gamma_2\neq 0$, $\gamma_5\neq 0$.

4. If $rk(M_{(1)})<rk(M'_{(1)})$, $rk(M_{(2)})<rk(M'_{(2)})$,
$rk(M_{(3)})<rk(M'_{(3)})$, $rk(M_{(4)})=rk(M'_{(4)})$
one can make $\gamma'_1=1$, $\gamma'_2\neq 0$, $\gamma'_4\neq 0$, $\gamma'_5=0$ to get
MSC \[\left(
\begin{array}{ccccccccc}
\alpha_1 & \alpha_2& 0&\alpha_4&\alpha_5&0&0&0&0\\
\beta_1 & \beta_2& 0&\beta_4&\beta_5&0&0&0&0\\
1 & \gamma_2& -\alpha_1-\beta_2&\gamma_4&0&1-\alpha_4-\beta_5&1-\alpha_1-\beta_4&-\alpha_2-\beta_5&0\end{array}\right),\] where $rk(M_{(1,2,3,4)})=rk(M_{(2-\gamma_2\times 1,3-\gamma_4\times 1, 4)})$, $rk(M_{(\gamma_2\times 3-\gamma_4\times 2,4)})=rk(M_{(2,3,4)})$, $rk(M_{(3-\gamma_4\times 1,4)})=rk(M_{(1,3,4})$, $rk(M_{(2-\gamma_2\times 1,4)})=rk(M_{(1,2,4)})$, 
$rk(M_{(2-\gamma_2\times 1),3-\gamma_4\times 1})=rk(M_{(1,2,3)})$, $rk(M_{(4)})=rk(M_{(3,4)})$, $rk(M_{(4)})=rk(M_{(2,4)})$, $rk(M_{(4)})=rk(M_{(1,4)})$, $rk(M_{(3-\gamma_4\times 1)})=rk(M_{(1,3)})$, $rk(M_{(1)})=0$, $rk(M_{(2)})=0$, $rk(M_{(3)})=0$ and $\gamma_2\neq 0$, $\gamma_4\neq 0$.

But the system of equations $rk(M_{(1)})=0$, $rk(M_{(2)})=0$, $rk(M_{(3)})=0$ has  solution $\alpha_1=1, \alpha_4=\alpha_2, \beta_1=0, \beta_2=0, \beta_4=0, \beta_5= 1$ and one can make $\gamma'_1= 1$ to have
\[A_{61,2}=\left(
\begin{array}{ccccccccc}
1 & \alpha_2& 0&\alpha_4&\alpha_5&0&0&0&0\\
0 & 0& 0&0&1&0&0&0&0\\
1 & \gamma_2& 1&\gamma_4&0&\alpha_4&0&1+\alpha_2&0\end{array}\right),\] where $\gamma_2\neq 0$, $\gamma_4\neq 0$.

5. The last one is $rk(M_{(1)})<rk(M'_{(1)})$, $rk(M_{(2)})<rk(M'_{(2)})$,
$rk(M_{(3)})<rk(M'_{(3)})$, $rk(M_{(4)})<rk(M'_{(4)})$, that is  $rk(M_{(1,2,3,4)})=0$,  $\gamma_1\neq 0$, $\gamma_2\neq 0$, $\gamma_4\neq 0$, $\gamma_5\neq 0$, case. 
The equation $rk(M_{(1,2,3,4)})=0$ has solution $\alpha_1=1, \alpha_4=\alpha_2, \alpha_5=0, \beta_1=0, \beta_2=0, \beta_4=0, \beta_5= 1$ and one can make $\gamma'_1= 1$ to get

\[A_{62,2}=\left(
\begin{array}{ccccccccc}
1 & \alpha_2& 0&\alpha_2&0&0&0&0&0\\
0 & 0& 0&0&1&0&0&0&0\\
1 & \gamma_2& 1&\gamma_4&\gamma_5&\alpha_2 &0 &1+\alpha_2&0\end{array}\right),\] where 
$\gamma_2\neq 0,\ \gamma_4\neq 0,\ \gamma_5\neq 0.$

\end{document}